\documentclass[12pt]{amsart}
\usepackage{indentfirst,amsmath,amssymb,amsthm,mathrsfs}
\usepackage[colorlinks,linkcolor=blue,citecolor=red]{hyperref}

\numberwithin{equation}{section}
\theoremstyle{plain}
\newtheorem{theorem}{Theorem}[section]
\newtheorem{lemma}{Lemma}[section]
\newtheorem{corollary}{Corollary}[section]
\newtheorem{proposition}{Proposition}[section]
\theoremstyle{definition}
\newtheorem{definition}{Definition}[section]
\newtheorem{example}{Example}[section]

\theoremstyle{remark}
\newtheorem{remark}{Remark}[section]
\newtheorem*{remark*}{Remark}

\providecommand{\IIref}[2]{#2}

\begin{document}
\date{}
\title[Shrinking Sasaki--Ricci Solitons I]{On the Structure of Gradient
Shrinking Sasaki--Ricci Solitons of Dimension Five I: Curvature Estimates and
Local Transverse Splitting}
\author{Shu-Cheng Chang}
\address{Shanghai Institute for Mathematics and Interdisciplinary Sciences,
Shanghai 200433, China}
\email{scchang@math.ntu.edu.tw}
\author{Chien Lin}
\address{Mathematical Science Research Center, Chongqing University of
Technology, Chongqing 400054, China}
\email{chienlin@cqut.edu.cn}
\author{Hongbing Qiu}
\address{School of Mathematics and Statistics, Wuhan University, Wuhan
430072, China}
\email{hbqiu@whu.edu.cn}

\begin{abstract}
This is the first of two papers on the structure of complete noncompact
gradient shrinking Sasaki--Ricci solitons of dimension five. Here we study
their local transverse geometry without any regularity assumption on the
Reeb foliation. We prove that an upper bound for the transverse scalar
curvature implies that the transverse curvature tensor is bounded by a
multiple of the transverse scalar curvature; this is the Sasakian analogue
of the curvature estimate of Munteanu--Wang. Using a transverse version of
Perelman's reduced geometry we show that, whenever a smooth leaf space is
available, the associated ancient Sasaki--Ricci flow is transversely
noncollapsed, and that at transverse infinity the developed local limits
split off a flat complex line, the remaining surface factor having
nonnegative curvature, as in the first splitting theorem of Naber. These
local conclusions do not produce a global complex surface; the
classification of the global quotients and the structure theory of
irregular Reeb foliations are carried out in the second paper.
\end{abstract}

\maketitle
\tableofcontents

\section{Introduction}

\label{sec:intro-conv}

The Sasaki--Ricci flow was introduced by Smoczyk, Wang and Zhang \cite{SWZ}
as the odd-dimensional counterpart of the K\"{a}hler--Ricci flow. On a
Sasakian manifold $(M^{2n+1},\xi ,\eta ,\Phi ,g)$ the Reeb field $\xi $ is
kept fixed and the transverse K\"{a}hler metric evolves by
\begin{equation*}
\frac{\partial }{\partial t}g^{T}(x,t)=-\mathrm{Ric}^{T}(x,t).
\end{equation*}
Smoczyk--Wang--Zhang proved the existence of Sasaki $\eta $-Einstein metrics
when the basic first Chern class satisfies $c_{1}^{B}(M)<0$ or
$c_{1}^{B}(M)=0$. The transverse Fano case $c_{1}^{B}(M)>0$ is still widely
open; see Collins \cite{Co1,Co2}, He \cite{He}, Collins--Jacob \cite{CJ} and
\cite{CHL,CHLW,CLLW,CCLW} for later developments.

Gradient Sasaki--Ricci solitons, also called transverse K\"{a}hler--Ricci
solitons, are the self-similar solutions of this flow. They generalize
$\eta $-Einstein metrics and model the singularities of the flow
\cite{FOW,CLL}; see also \cite{CHW}. Throughout this paper a
\emph{gradient Sasaki--Ricci soliton} is a Sasakian manifold
$(M,\xi ,\eta ,\Phi ,g)$ together with a basic function $f$ satisfying the
\emph{full real} transverse equation
\begin{equation}
\mathrm{Ric}^{T}+\operatorname{Hess}^{T}f=\lambda g^{T},\qquad \xi f=0.
\label{eq:intro-soliton}
\end{equation}
It is shrinking, steady or expanding according as $\lambda >0$, $\lambda =0$
or $\lambda <0$. In complex coordinates \eqref{eq:intro-soliton} is the
equation $R_{j\overline{k}}^{T}+f_{j\overline{k}}=\lambda g_{j\overline{k}}$
together with the vanishing of the $(2,0)+(0,2)$-part of
$\operatorname{Hess}^{T}f$.
Since $\mathrm{Ric}^{T}=\mathrm{Ric}_{g}|_{D}+2g^{T}$ on $D=\ker \eta $, it
reads $R_{j\overline{k}}+f_{j\overline{k}}=Ag_{j\overline{k}}$ with
$A=\lambda -2$ in terms of the Ricci tensor of $g$; for constant $f$ it is
the $\eta $-Einstein condition, and the Sasaki--Einstein condition when
$\lambda =2n+2$. A transverse $D$-homothety with parameter $\tau >0$ replaces
$\lambda $ by $\lambda /\tau $ (Remark~\ref{rem:Dhom}). For calculations involving a shrinking soliton we use the following normalizations, unless another scale is explicitly stated:
\begin{equation}
\lambda =2n+2\ \text{ in Sections~\ref{sec:identities}
and~\ref{sec:curvature}},\qquad \lambda =\tfrac{1}{2}\ \text{ in
Sections~\ref{sec:flows}, \ref{sec:splitting} and
Appendix~\ref{app:reduced}},  \label{eq:intro-normal}
\end{equation}
and in both cases the additive constant of the potential of a shrinking
soliton is fixed by the Hamilton identity $R^{T}+|\nabla ^{T}f|^{2}=2\lambda f$
(Lemma~\ref{lem:nonneg}); such solitons are called \emph{normalized}.

The present paper and its sequel \cite{PartII} study complete noncompact
gradient shrinking Sasaki--Ricci solitons of dimension five, that is, of
transverse complex dimension two, without any regularity assumption on the
Reeb foliation. If the Reeb foliation is regular, the leaf space is a
complete gradient K\"{a}hler--Ricci shrinker surface, and such surfaces are
now well understood: their curvature is bounded by Li--Wang \cite{LW}, and,
by the work of Conlon--Deruelle--Sun \cite{CDS}, Cifarelli--Conlon--Deruelle
\cite{CCD} and Bamler--Cifarelli--Conlon--Deruelle \cite{BCCD}, the
noncompact ones are, up to biholomorphism, the Gaussian shrinker on
$\mathbb{C}^{2}$, the Feldman--Ilmanen--Knopf (FIK) shrinker \cite{FIK} on
the blow-up of $\mathbb{C}^{2}$ at one point, the product of the Gaussian
and the round shrinker on $\mathbb{C}\times \mathbb{CP}^{1}$, and the
Bamler--Cifarelli--Conlon--Deruelle (BCCD) shrinker \cite{BCCD} on the
blow-up of $\mathbb{C}\times \mathbb{CP}^{1}$ at one point (see also
\cite{MW3,LWheat,LZ}). For quasi-regular, and above all for irregular, Reeb
foliations the leaf space is an orbifold or is not Hausdorff, and these
results do not apply directly. The present paper develops the local
transverse theory, which holds for every Reeb type; the global consequences
are derived in \cite{PartII}.

Our first result is the Sasakian analogue of the curvature estimate of
Munteanu--Wang \cite{MW1} for four-dimensional shrinking Ricci solitons. We
write $R_{g}$ and $R^{T}$ for the ambient and the transverse scalar
curvature, and $\mathrm{Rm}^{T}$, $\mathrm{Ric}^{T}$ for the transverse
curvature and Ricci tensors; in dimension five $R^{T}=R_{g}+4$.

\begin{theorem}
\label{thm:intro-curvature}Let $(M^{5},\xi ,\eta ,\Phi ,g,f)$ be a
connected, complete, noncompact gradient shrinking Sasaki--Ricci soliton,
with $f$ additively normalized by
$R^{T}+|\nabla^{T}f|^{2}=2\lambda f$. If the transverse scalar curvature
$R^{T}$ is bounded from above, then the
transverse curvature tensor $\mathrm{Rm}^{T}$ and the curvature tensor
$\mathrm{Rm}_{g}$ of $g$ are bounded. More precisely, either

\begin{enumerate}
\item $R^{T}\equiv 0$ and $\mathrm{Rm}^{T}\equiv 0$, or

\item $R^{T}>0$ on $M$, and there are constants $C,c>0$ such that
\begin{equation*}
R^{T}f\geq c\qquad \text{and}\qquad |\mathrm{Rm}^{T}|+|\mathrm{Ric}^{T}|\leq
CR^{T}.
\end{equation*}
\end{enumerate}
In both cases $|\mathrm{Rm}^{T}|\leq CR^{T}$ and
$|\mathrm{Rm}_{g}|\leq C(1+R^{T})=C(5+R_{g})$ on $M$.
\end{theorem}

The dichotomy and the lower bound $R^{T}f\geq c$ are proved in
Section~\ref{sec:identities} (Lemma~\ref{lem:nonneg} and
Proposition~\ref{prop:flat}), following B.-L.~Chen \cite{Che} and
Chow--Lu--Yang \cite{CLY} in the Riemannian case, and the curvature bounds
are Theorem~\ref{thm:curvature}. The proof does not use a leaf space: all
tensor computations are made in local transverse K\"{a}hler charts, and the
global maximum principles are carried out on the exhaustion space of
Lemma~\ref{lem:reeb}, which is $M$ itself when the closure of the Reeb flow
is a compact torus and the smooth quotient $M/\mathbb{R}$ when the Reeb flow
is a proper $\mathbb{R}$-action. The same method gives, in the nonflat case,
the derivative estimates $|\nabla ^{T}\mathrm{Rm}^{T}|\leq CR^{T}$ and
$|\nabla ^{T}\log R^{T}|\leq C$ (Proposition~\ref{prop:gradT}).

\begin{remark}
\label{rem:hypothesis}The upper bound for $R^{T}$ in
Theorem~\ref{thm:intro-curvature} is a genuine hypothesis only when no
complete smooth leaf space is available. Indeed, if the structure is
regular, or if the Reeb flow is a proper $\mathbb{R}$-action as in
Lemma~\ref{lem:reeb}, then the leaf space $Z$ is a complete gradient
K\"{a}hler--Ricci shrinker surface, and its scalar curvature is bounded by
the theorem of Li--Wang \cite[Theorem~1.2]{LW} (see also \cite{LZ}). The
same holds, after passing to a finite cover, for a quasi-regular structure
whose leaf space is a developable orbifold; see
\cite[Theorem~\IIref{thm:quotient}{4.1}]{PartII}. In these cases
$\sup_{M}R^{T}<\infty $ is automatic, and on the smooth leaf space the
estimate of Theorem~\ref{thm:intro-curvature} is the theorem of
Munteanu--Wang \cite{MW1}. Thus the boundedness of $R^{T}$ has to be assumed
only for quasi-regular structures whose leaf space is not developable and
for irregular structures. In these cases there is no smooth leaf space to
which \cite{LW} and \cite{MW1} could be applied, and the proof is carried
out on the exhaustion space of Lemma~\ref{lem:reeb}.
\end{remark}

For the analysis at transverse infinity we use the ancient flow associated
with a soliton. It is convenient to keep the normal bundle
$TM/\mathbb{R}\xi $ fixed; the associated \emph{transverse
K\"{a}hler--Ricci flow} is a family of transverse metrics satisfying
\begin{equation}
\partial _{t}h(t)=-2\,\mathrm{Ric}^{T}(h(t)).  \label{eq:marked-flow}
\end{equation}

\begin{remark}
\label{rem:flow-normal}We always take $h_{j\bar{k}}=h(\partial
_{z^{j}},\partial _{\bar{z}^{k}})$
and $\mathrm{Ric}_{j\bar{k}}=\mathrm{Ric}(\partial _{z^{j}},\partial
_{\bar{z}^{k}})$ by complex-bilinear extension of the real tensors. Then
$\partial _{t}h_{j\bar{k}}=-2\mathrm{Ric}_{j\bar{k}}$. The flow of
\cite{SWZ}, $\partial _{s}\widehat{h}=-\mathrm{Ric}^{T}(\widehat{h})$, is
related to \eqref{eq:marked-flow} by $h(t)=\widehat{h}(2t)$. All curvature,
reduced-length and reduced-volume identities in Appendix~\ref{app:reduced}
use \eqref{eq:marked-flow}; in the backward time $\tau =-t$ it reads
$\partial _{\tau }h=2\,\mathrm{Ric}^{T}(h)$, which is Perelman's
normalization \cite{Per1}. The letter $\tau$ also denotes the parameter of
a transverse $D$-homothety (Remark~\ref{rem:Dhom}); the two uses never
occur in the same argument.
\end{remark}

Perelman's reduced geometry involves integrals over a leaf space. For
its global statements we therefore work in the following setting.

\textbf{Setting (S)} \quad 
Let $h(t)$, $t\in I$, be the transverse flow
\eqref{eq:marked-flow} on the fixed normal bundle
$TM/\mathbb{R}\xi$. We assume that the following data have been chosen.

\begin{enumerate}
\item A connected finite covering map $q:\widetilde M\to M$, possibly
the identity, is fixed for the whole interval $I$. On $\widetilde M$ we
use the lifted Reeb field $\widetilde\xi$ and the lifted transverse
metrics $\widetilde h(t)=q^*h(t)$. The lifted Reeb flow defines either
a free circle action or a free proper $\mathbb R$-action. Thus its
leaf space is a smooth Hausdorff manifold
\[
Z=\widetilde M/\mathcal F_{\widetilde\xi},\qquad
\pi:\widetilde M\longrightarrow Z.
\]
The covering $q$ is a covering of the total space. It is not required
to be a normal (Galois) covering; regularity of the lifted Reeb
foliation means freeness of the circle action.

\item The manifold $Z$, its complex structure $J_Z$, and the quotient
map $\pi$ are independent of $t$. For every $t\in I$, the metric
$\widetilde h(t)$ descends to a complete K\"{a}hler metric $h_Z(t)$ on
$Z$, and
\[
\widetilde h(t)=\pi^*h_Z(t),\qquad
\partial_t h_Z(t)=-2\,\mathrm{Ric}_{h_Z(t)}.
\]
The pullback identity is an identity of transverse tensors on
$T\widetilde M/\mathbb R\widetilde\xi$. Completeness refers to each
time slice; bounded curvature and other analytic hypotheses are
stated separately in the results where they are needed.

\item For a pointed quantity based at $(x,t_0)$, choose a lift
$\widetilde x\in q^{-1}(x)$ and use
$z=\pi(\widetilde x)$ as the quotient basepoint. All global reduced
lengths, reduced volumes, and noncollapsing inequalities in setting
(S) are defined on this chosen quotient flow, using
\[
d\mu_t^T=d\mathrm{vol}_{h_Z(t)}.
\]
In particular, the volume of a transverse ball based at the chosen
lift is
\[
\operatorname{Vol}^T_t\bigl(B^T(\widetilde x,r)\bigr)
=\operatorname{Vol}_{h_Z(t)}\bigl(B_{h_Z(t)}(z,r)\bigr).
\]
It is a $2n$-dimensional quotient volume, not the
$(2n+1)$-dimensional volume of its inverse image in the total space.
These quantities and their noncollapsing constants refer to the
specified cover and quotient; transferring them to a different
quotient requires a separate argument.
\end{enumerate}

Once these choices have been made, statements in setting (S) use $M$
and $x$ for the chosen total space and chosen lift, unless the
original space or the covering map $q$ is explicitly mentioned.

For a regular structure, or for a closed free proper
$\mathbb R$-action as in Lemma~\ref{lem:reeb}, one takes
$q=\mathrm{id}_M$. Developable quasi-regular shrinkers enter this
setting after passage to a finite cover whose lifted Reeb foliation
is regular; see \cite[Lemma~\IIref{lem:finite}{4.2} and
Theorem~\IIref{thm:good}{4.4}]{PartII}. For these shrinkers, the lifted
soliton gradient preserves the Reeb foliation and projects to the smooth
quotient; hence the same cover and quotient carry the entire canonical
ancient flow of Proposition~\ref{lem:selfsim}. For a general flow, the
existence of one such cover and smooth quotient valid for every
$t\in I$ is an additional hypothesis. It is not inferred merely
from quasi-regularity at a single time.

The data in setting (S) are assumptions of every global result below
that invokes this setting. References to \cite{PartII} describe
classes of examples and applications; they are not used to construct
the quotient within those proofs. In particular, the regular case
has $q=\mathrm{id}_M$ and is treated directly. Local transverse
statements that do not assume setting (S) retain their stated scope.

Our first result on the flow is the transverse version of Perelman's
no-local-collapsing theorem \cite{Per1} (Theorem~\ref{thm:perelmanNC}).

\begin{theorem}
\label{thm:perelmanNC-intro} Let $(M,g(t))$, $t\in \lbrack 0,T]$, be a
Sasaki--Ricci flow whose transverse metrics $h(t)$, normalized by
\eqref{eq:marked-flow}, are as in setting (S), and let $x\in M$, $0<\bar{\tau}\leq T$
and $V>0$. There is $\kappa =\kappa (n,V)>0$ such that if
$0<r\leq \sqrt{\bar{\tau}}$, $|\mathrm{Rm}^{T}|\leq r^{-2}$ on the
transverse parabolic tube $P^{T}(x,T,r)$, and the reduced transverse volume
satisfies $V_{(x,T)}^{T}(\bar{\tau})\geq V$, then
\begin{equation*}
\operatorname{Vol}^{T}\!\big(B_{r}^{T}(x)\big)\geq \kappa r^{2n},
\end{equation*}
where $\operatorname{Vol}^{T}$ is the quotient measure \eqref{eq:analyticspace}.
\end{theorem}

Applied to the ancient flow of a soliton with bounded transverse curvature,
together with the asymptotic reduced volume (compare Naber
\cite[Section~3]{Nab} and Carrillo--Ni \cite{CN} in the Riemannian case),
this gives a uniform noncollapsing constant (Theorem~\ref{thm:kappa} and
Corollary~\ref{cor:controlled}).

\begin{corollary}
Let $(M,\xi ,\eta ,\Phi ,g,f)$ be a normalized complete gradient shrinking
Sasaki--Ricci soliton with $\lambda=\tfrac12$ and bounded transverse
curvature in setting (S).
Then there is
\begin{equation*}
\kappa =\kappa \big(n,\operatorname{Vol}_{f}^{T}(M)\big)>0
\end{equation*}
such that the associated Sasaki--Ricci flow of Proposition~\ref{lem:selfsim}
is transversely $\kappa $-noncollapsed.
\end{corollary}

We next derive the Sasakian version of the splitting theorem of Naber
\cite{Nab} at infinity. In setting (S) it is a statement about the leaf
space (Theorem~\ref{P54}).

\begin{theorem}
\label{P54-intro}Let $(M^{5},\xi ,\eta ,\Phi ,g,f)$ be a normalized complete gradient
shrinking Sasaki--Ricci soliton with $\lambda=\tfrac12$ in setting (S). Suppose that its smooth
quotient $(Z,h,J,f_{Z})$ is noncompact and $|\mathrm{Rm}(h)|\leq C$, and let
$h(t)=(-t)\Psi _{t}^{\ast }h$, $t<0$, be the quotient flow induced by the
ancient flow of Proposition~\ref{lem:selfsim}, so that $h(-1)=h$. For any
$q_{i}\in Z$ with $f_{Z}(q_{i})\rightarrow \infty $, a subsequence of
$(Z,h(t),q_{i})$ converges smoothly on compact space-time sets, $t<0$, to a
complete noncollapsed ancient K\"{a}hler--Ricci flow
$(Z_{\infty },h_{\infty }(t),q_{\infty })$ with
$|\mathrm{Rm}(h_{\infty }(t))|\leq C/|t|$, whose time $-1$ slice is the
pointed limit of $(Z,h,q_{i})$. Its universal cover splits as
\begin{equation*}
(\widehat{Z}_{\infty },\widehat{h}_{\infty }(t))=(\mathbb{C},h_{\mathrm{E}})
\times (N,k(t)),
\end{equation*}
where $N$ is complete and $k(t)$ is a noncollapsed ancient Ricci flow on a
surface with $|\mathrm{Rm}(k(t))|\leq C/|t|$. In particular, the limit
$(Z_{\infty },h_{\infty }(t))$ itself splits locally in this fashion.
\end{theorem}

For an arbitrary Reeb foliation there is in general no global leaf space,
and we formulate the splitting for developed local limits, that is, limits
of the structures pulled back by the exponential map to tangent balls of a
fixed radius (Definition~\ref{def:developed}). The following result is
Theorem~\ref{thm:splitting}; its Riemannian counterpart is due to Naber
\cite{Nab}.

\begin{theorem}
\label{thm:splitting-intro}Let $p_{i}\in M$ with $f(p_{i})\rightarrow \infty $
on a complete five-dimensional gradient shrinking Sasaki--Ricci soliton $M$
with bounded transverse curvature (for instance, with $R^{T}$ bounded
above, by Theorem~\ref{thm:intro-curvature}). Every
developed local transverse limit at these basepoints has orthonormal parallel fields
$V$, $JV$. Consequently it splits locally as
\begin{equation}
(U,h_{\mathrm{E}})\times (\Sigma ,h_{\Sigma }),\qquad U\subset \mathbb{C},
\label{eq:local-split-intro}
\end{equation}
with its product K\"{a}hler structure. Suppose that a nonflat product limit
satisfies the local soliton equation, that is, on a neighborhood containing
its product chart there is a smooth function $F$ with
$\mathrm{Ric}(h_{\lim })+\operatorname{Hess}F=\lambda h_{\lim }$. Then
\begin{equation*}
F(z,y)=\frac{\lambda }{2}|z|^{2}+\ell (z)+\psi (y),\qquad
\operatorname{Hess}_{h_{\Sigma }}\psi =(\lambda -K)h_{\Sigma },
\end{equation*}
where $\ell $ is affine and $K$ is the Gaussian curvature of $h_{\Sigma }$;
in particular $dK=K\,d\psi $. If the surface factor, together with its
potential, extends to a complete smooth connected two-dimensional gradient
shrinker with the same constant $\lambda $, then $K\equiv \lambda $, and
the limit is locally the product $\mathbb{C}\times \mathbb{CP}^{1}$ with the
round metric of curvature $\lambda $ on the second factor.
\end{theorem}

Moreover, in the nonflat transverse case, the curvature operator is
asymptotically of rank one at transverse infinity. If $R^{T}$ is bounded above, then outside a fixed $f$-sublevel, where $\nabla^{T}f\ne0$,
$\Vert \mathcal{R}^{T}-\frac{1}{2}R^{T}P_{E}\Vert _{\mathrm{op}}\leq
CR^{T}/\sqrt{f}$. Here
$E=\operatorname{span}\{\nabla^{T}f,J\nabla^{T}f\}^{\perp}\subset D$
is a complex line and $P_E$ is the orthogonal projection in $\Lambda^2D$
onto its unit area bivector. Thus $\mathcal{R}^{T}/R^{T}$ is within
$Cf^{-1/2}$ of the rank-one operator $\frac12P_E$. The curvature decays
quadratically if $R^{T}\rightarrow0$ uniformly at transverse infinity
(Theorem~\ref{thm:general}).
Consequently every developed local limit at transverse infinity has
nonnegative transverse curvature operator and its surface factor has
nonnegative Gaussian curvature (Corollary~\ref{prop:conditional}); if in
addition $R^{T}(p_{i})\rightarrow 0$, the limit is transversely flat and
locally isomorphic to the Heisenberg model (Corollary~\ref{cor:flatlimit}).

These local conclusions do not by themselves produce a global complex
surface: developed local limits need not be injective images of
neighborhoods in $M$, different subsequences may have different limits, and
for irregular foliations the surface factor need not be round. The
classification of the global quotients, the structure theory of
quasi-regular and irregular Reeb foliations and the relevant examples are
given in \cite{PartII}.

\textbf{Organization.} Section~\ref{sec:prelim} collects the transverse
K\"{a}hler geometry of Sasakian manifolds, Sasaki--Ricci solitons and the
types of Reeb foliations. Section~\ref{sec:identities} derives the soliton
identities, the exhaustion space, the nonnegativity of $R^{T}$ and the
growth of the potential. Section~\ref{sec:curvature} proves the curvature
estimate and the derivative estimate of Proposition~\ref{prop:gradT}.
Section~\ref{sec:flows} introduces the controlled Sasaki--Ricci
flows associated with solitons and the developed local limits.
Section~\ref{sec:splitting} contains the noncollapsing theorem, the
Naber-type splitting on the leaf space, the local splitting for arbitrary
Reeb type and the asymptotic structure at transverse infinity.
Appendix~\ref{app:reduced} records the transverse version of Perelman's
reduced geometry used in Section~\ref{sec:splitting}.

\section{Preliminaries}

\label{sec:prelim}

\subsection{Sasakian Structures and Transverse K\"{a}hler Geometry}

Sasakian geometry is the odd-dimensional analogue of K\"{a}hler geometry; we
refer to \cite{BG} for background. A Riemannian manifold $(M^{2n+1},g)$ is
called a \emph{Sasakian manifold} if its metric cone
\begin{equation*}
(C(M),\overline{g},J)=(M\times \mathbb{R}^{+},\,r^{2}g+dr^{2},\,J)
\end{equation*}
is K\"{a}hler. We identify $M$ with $\{r=1\}\subset C(M)$. The \emph{Reeb
vector field} and the \emph{contact form} are
\begin{equation*}
\xi =J\Bigl(r\frac{\partial }{\partial r}\Bigr),\qquad \eta (Y)=g(\xi ,Y).
\end{equation*}
Then $\eta (\xi )=1$ and $\iota _{\xi }d\eta =0$, and $\xi $ is a Killing
field of unit length. The following conditions on a Riemannian manifold
$(M^{2n+1},g)$ are equivalent (see \cite{BG}):

\begin{enumerate}
\item there is a Killing vector field $\xi $ of unit length on $M$ such that
the tensor field $\Phi $ of type $(1,1)$ defined by $\Phi (Y)=\nabla _{Y}\xi $
satisfies
\begin{equation*}
(\nabla _{X}\Phi )(Y)=g(\xi ,Y)X-g(X,Y)\xi
\end{equation*}
for all vector fields $X$ and $Y$ on $M$;

\item there is a Killing vector field $\xi $ of unit length on $M$ such that
the Riemann curvature tensor satisfies
\begin{equation*}
R(X,\xi )Y=g(\xi ,Y)X-g(X,Y)\xi
\end{equation*}
for all vector fields $X$ and $Y$ on $M$;

\item the metric cone $(C(M),\overline{g})$ is K\"{a}hler.
\end{enumerate}

Moreover, there is a natural splitting
\begin{equation*}
TC(M)=L_{r\partial _{r}}\oplus L_{\xi }\oplus D,
\end{equation*}
where $L_{r\partial _{r}}$ and $L_{\xi }$ are the line bundles spanned by
$r\partial _{r}$ and $\xi $, and $D=\ker \eta $ is the contact distribution.
The complex structure $J$ of the cone and the tensor $\Phi $ are related by
\begin{equation*}
JY=\Phi (Y)-\eta (Y)\,r\frac{\partial }{\partial r},\qquad \Phi ^{2}=-I+\eta
\otimes \xi ,
\end{equation*}
for vector fields $Y$ on $M$, and
\begin{equation*}
g(\Phi X,\Phi Y)=g(X,Y)-\eta (X)\eta (Y).
\end{equation*}
Thus
\begin{equation*}
g=g^{T}+\eta \otimes \eta ,
\end{equation*}
where $g^{T}(X,Y)=g(\Phi X,\Phi Y)$ is the transverse metric, that is, the
restriction of $g$ to $D$ extended by zero in the direction of $\xi $.
We also write $h=g|_D$ for the positive-definite metric on $D$, or
for the corresponding metric on the normal bundle $TM/\mathbb R\xi$.
Thus $g^T$ on $TM$ is degenerate, whereas $h$ is the transverse
Riemannian metric used for traces and norms. Unless an ambient metric is
explicitly indicated, norms of transverse tensors are taken with $h$;
$R^T=\operatorname{tr}_h\mathrm{Ric}^T$ is the real trace in dimension
$2n$. We call $(M,\xi ,\eta ,\Phi ,g)$ a \emph{Sasakian structure}.

Let $\{U_{\alpha }\}_{\alpha \in A}$ be an open covering of the Sasakian
manifold $(M,\xi ,\eta ,\Phi ,g)$, and let
\begin{equation*}
\pi _{\alpha }:U_{\alpha }\rightarrow V_{\alpha }\subset \mathbb{C}^{n}
\end{equation*}
be foliation-chart submersions whose fibers are the connected plaques of
the Reeb foliation $\mathcal F_\xi$ in $U_\alpha$. After restricting to
sufficiently small components of an overlap, there are biholomorphic
transverse transition maps $\gamma_{\alpha\beta}$ with
$\pi_\alpha=\gamma_{\alpha\beta}\circ\pi_\beta$. These are local
transition maps: no inverse of a submersion $\pi_\beta$ and no global leaf
space are used.
For every $p\in U_{\alpha }$ the differential
\begin{equation*}
d\pi _{\alpha }:D_{p}\rightarrow T_{\pi _{\alpha }(p)}V_{\alpha }
\end{equation*}
is an isomorphism, where $D=\ker \eta \subset TM$. Since $\xi $ generates
isometries, the restriction of the Sasakian metric $g$ to $D$ induces a
well-defined Hermitian metric $g_{\alpha }^{T}$ on $V_{\alpha }$, and this
metric is K\"{a}hler. More precisely, let $z^{1},\ldots ,z^{n}$ be holomorphic
coordinates on $V_{\alpha }$, pulled back to $U_{\alpha }$ and denoted by the
same letters, and let $x$ be a coordinate along the leaves with $\xi =\partial
/\partial x$. Then $\{x,z^{1},\ldots ,z^{n}\}$ are foliated coordinates on
$U_{\alpha }$, and $D\otimes \mathbb{C}$ is spanned by the vector fields
\begin{equation*}
Z_{j}=\frac{\partial }{\partial z^{j}}-\eta \Bigl(\frac{\partial }{\partial
z^{j}}\Bigr)\frac{\partial }{\partial x},\qquad j=1,\ldots ,n,
\end{equation*}
and their conjugates. Since $\iota _{\xi }d\eta =0$,
\begin{equation*}
d\eta (Z_{j},\overline{Z_{k}})=d\eta \Bigl(\frac{\partial }{\partial
z^{j}},\frac{\partial }{\partial \overline{z}^{k}}\Bigr).
\end{equation*}
Hence the K\"{a}hler form $\omega _{\alpha }^{T}$ of $g_{\alpha }^{T}$ on
$V_{\alpha }$, which corresponds to the restriction of $\frac{1}{2}d\eta $ to
a slice $\{x=\mathrm{const}\}$ in $U_{\alpha }$, is closed. The collection of
K\"{a}hler metrics $\{g_{\alpha }^{T}\}$ on $\{V_{\alpha }\}$ is called a
\emph{transverse K\"{a}hler metric}; we write $g^{T}$, $\omega
^{T}=\frac{1}{2}d\eta $ and $J=\Phi |_{D}$ for the transverse metric, its
K\"{a}hler form and the transverse complex structure. For example, on the
Heisenberg group $\mathbb{H}_{n}=\mathbb{C}^{n}\times \mathbb{R}$ with
coordinates $(z,t)$, the contact form
\begin{equation*}
\eta
=dt+\frac{\sqrt{-1}}{2}\sum_{j}(z^{j}d\overline{z}^{j}-\overline{z}^{j}dz^{j})=dt+\sum_{j}(x_{j}dy_{j}-y_{j}dx_{j}),\qquad
z^{j}=x_{j}+\sqrt{-1}y_{j},
\end{equation*}
the Reeb field $\xi =\partial /\partial t$ and the frame $Z_{j}=\partial
/\partial z^{j}+\frac{\sqrt{-1}}{2}\overline{z}^{j}\partial /\partial t$ give
the flat model: a transverse K\"{a}hler potential is $\frac{1}{2}|z|^{2}$, the
transverse metric is Euclidean, and $d\eta =2\omega _{\mathrm{E}}$ with
$\omega _{\mathrm{E}}=\sum_{j}dx_{j}\wedge dy_{j}$. This is the Heisenberg
model \eqref{eq:heisenberg} used below.

Let $\widetilde{X},\widetilde{Y},\widetilde{W},\widetilde{Z}\in \Gamma (D)$
correspond under $d\pi _{\alpha }$ to $X,Y,W,Z\in \Gamma (TV_{\alpha })$. The
Levi-Civita connection $\nabla ^{T}$ of the transverse K\"{a}hler metric
$g^{T}$, its curvature and its Ricci tensor are given by
\begin{equation*}
\begin{array}{rcl}
\nabla _{X}^{T}Y & := & d\pi _{\alpha }(\nabla _{\widetilde{X}}\widetilde{Y}),
\\
\widetilde{\nabla _{X}^{T}Y} & = & \nabla _{\widetilde{X}}\widetilde{Y}+g(\Phi
\widetilde{X},\widetilde{Y})\xi , \\
\mathrm{Rm}^{T}(X,Y,Z,W) & = &
\mathrm{Rm}_{g}(\widetilde{X},\widetilde{Y},\widetilde{Z},\widetilde{W})+g(\Phi
\widetilde{Y},\widetilde{W})g(\Phi \widetilde{X},\widetilde{Z}) \\
& & -g(\Phi \widetilde{X},\widetilde{W})g(\Phi
\widetilde{Y},\widetilde{Z})+2g(\Phi \widetilde{X},\widetilde{Y})g(\Phi
\widetilde{Z},\widetilde{W}), \\
\mathrm{Ric}^{T}(X,Z) & = &
\mathrm{Ric}_{g}(\widetilde{X},\widetilde{Z})+2g(\widetilde{X},\widetilde{Z}).
\end{array}
\end{equation*}
In the last line $\mathrm{Ric}_{g}$ is the Ricci tensor of $g$; the
contraction of $\mathrm{Rm}_{g}$ over $D$ alone is
$\mathrm{Ric}_{g}|_{D}-g^{T}$, by \eqref{eqn-SR3} below. Here and throughout
we use the conventions
\begin{equation}
\begin{array}{rcl}
R(X,Y)Z & := & \nabla _{X}\nabla _{Y}Z-\nabla _{Y}\nabla _{X}Z-\nabla
_{\lbrack X,Y]}Z, \\
R(X,Y,Z,W) & := & \langle R(Z,W)Y,X\rangle , \\
R_{\alpha \beta \gamma \delta } & := & R(e_{\alpha },e_{\beta },e_{\gamma
},e_{\delta })=\langle R(e_{\gamma },e_{\delta })e_{\beta },e_{\alpha }\rangle
, \\
\mathrm{Ric}(X,Y) & := & \sum_{\alpha }R(X,e_{\alpha },Y,e_{\alpha }),
\end{array}
\label{eqn-SR1}
\end{equation}
where $\{e_{\alpha }\}$ is a local orthonormal frame; thus $R_{\alpha \beta
\alpha \beta }$ is a sectional curvature. Let $\{e_{i}\}_{i=1}^{2n}$ be a
local orthonormal frame of $D=\ker \eta $. On a Sasakian manifold
\begin{equation}
R(X,\xi )Y=g(\xi ,Y)X-g(X,Y)\xi ,  \label{eqn-SR2}
\end{equation}
and therefore
\begin{equation*}
\mathrm{Ric}^{T}=\mathrm{Ric}_{g}+2g^{T}\qquad \text{on }D.
\end{equation*}

A $p$-form $\alpha $ on $M$ is \emph{basic} if $\iota _{\xi }\alpha =0$ and
$\iota _{\xi }d\alpha =0$; a function $u$ is basic if $\xi u=0$. We write
$d_{B}=\partial _{B}+\overline{\partial }_{B}$ for the restriction of $d$ to
basic forms, $\Delta _{B}$ for the basic Laplacian, and $c_{1}^{B}(M)$ for
the basic first Chern class.

Equivalently, in terms of the foliation, the \emph{transverse Levi-Civita
connection} $\nabla ^{T}$ on $D$ is
\begin{equation}
\nabla _{X}^{T}Y=\begin{cases}
\nabla _{X}Y+g(\Phi X,Y)\,\xi , & X,Y\in \Gamma (D), \\
\lbrack \xi ,Y], & X=\xi ,\ Y\in \Gamma (D),\end{cases}
\label{eq:pre-nablaT}
\end{equation}
i.e.\ $\nabla _{X}^{T}Y$ is the $D$-component of $\nabla _{X}Y$ for $X,Y\in
\Gamma (D)$.

On each chart $\nabla ^{T}$ is the pullback of the Levi-Civita connection
of $g_{\alpha }^{T}$. We denote by $\mathrm{Rm}^{T}$, $\mathrm{Ric}^{T}$ and
$R^{T}$ the transverse curvature tensor, Ricci tensor and scalar curvature
of $g^{T}$, and by $\mathrm{Rm}_{g}$, $\mathrm{Ric}_{g}$ and $R_{g}=R$ those
of $g$. For a basic function $u$ we set $\nabla ^{T}u=\pi _{D}\nabla
u=\nabla u$, $\operatorname{Hess}^{T}u=(\nabla ^{T})^{2}u$, and
\begin{equation}
L_{f}u=\Delta _{B}u-\langle \nabla ^{T}f,\nabla ^{T}u\rangle =\Delta
_{g}u-\langle \nabla f,\nabla u\rangle .  \label{eq:pre-Lf}
\end{equation}
Our Laplacian has the sign
$\Delta_Bu=\operatorname{tr}_h\operatorname{Hess}^Tu$.
The last equality uses $\Delta_Bu=\Delta_gu$ for basic $u$, because
the unit Reeb field has geodesic, hence minimal, orbits. It does not
assert that $D$ is integrable. For a transverse tensor $U$, the same
notation means
\[
L_fU=\Delta_B^TU-\nabla^T_{\nabla^Tf}U,\qquad
\Delta_B^TU=\sum_{i=1}^{2n}\left(
\nabla^T_{e_i}\nabla^T_{e_i}U-
\nabla^T_{\nabla^T_{e_i}e_i}U\right),
\]
where $\{e_i\}$ is a local $h$-orthonormal frame. Thus the tensor
Laplacian is the rough transverse Laplacian, not the Hodge Laplacian.

The transverse and ambient curvatures are related as follows. Let
$\{e_{i}\}_{i=1}^{2n}$ be a local orthonormal frame of $D$, and write $\Phi
_{ij}=g(\Phi e_{i},e_{j})$. We use the abbreviations
\begin{equation}
K=\mathrm{Rm}_{g}|_{D},\qquad S=\mathrm{Ric}_{g}|_{D},\qquad
T=\mathrm{Rm}^{T},\qquad Q=\mathrm{Ric}^{T},\qquad s=R^{T}.
\label{eq:pre-KSTQ}
\end{equation}
All lower-case indices run over $1,\ldots ,2n$, and semicolons and tensor
Laplacians refer to $\nabla ^{T}$ unless explicitly marked otherwise. With
the convention \eqref{eqn-SR1},
\begin{equation}
T_{ijkl}=K_{ijkl}+\mathring{A}_{ijkl},\qquad \mathring{A}_{ijkl}=\Phi
_{ik}\Phi _{jl}-\Phi _{il}\Phi _{jk}+2\Phi _{ij}\Phi _{kl}.
\label{eq:pre-conversion}
\end{equation}
For an orthonormal pair this reads $T_{ijij}=K_{ijij}+3\Phi _{ij}^{2}$,
which is O'Neill's formula for the local Riemannian submersions $\pi
_{\alpha }$; the tensor $\mathring{A}$ is an algebraic curvature tensor,
that is, it also satisfies the first Bianchi identity.
Contracting \eqref{eq:pre-conversion} and using \eqref{eqn-SR2} gives
\begin{equation}
\begin{gathered} \mathrm{Ric}^{T}=\mathrm{Ric}_{g}|_{D}+2g^{T},\qquad
R^{T}=R_{g}+2n,\qquad
\sum_{i=1}^{2n}\mathrm{Ric}_{g}(e_{i},e_{i})=R_{g}-2n,\\
\text{that
is,}\quad Q=S+2h,\qquad s=R+2n,\qquad \operatorname{tr}_{h}S=R-2n,\\
\sum_{i}K_{ijil}=S_{jl}-h_{jl},\qquad |Q|^{2}=|S|^{2}+4R. \end{gathered}
\label{eq:pre-ric}
\end{equation}
Moreover $\nabla ^{T}\mathring{A}=0$, and since $\xi $ is Killing all
curvature components of $g$ with a Reeb index are fixed algebraically by
\eqref{eqn-SR2}. In particular,
\begin{equation}
|\mathrm{Rm}_{g}|\leq C_{n}\bigl(1+|\mathrm{Rm}^{T}|\bigr).
\label{eq:pre-Rmg}
\end{equation}

\subsection{Sasaki--Ricci Solitons}

\begin{definition}
A Sasakian manifold $(M,\xi ,\eta ,\Phi ,g)$ is \emph{$\eta $-Einstein} if
there is a constant $A$ such that
\begin{equation*}
\mathrm{Ric}_{g}=Ag+(2n-A)\eta \otimes \eta .
\end{equation*}
\end{definition}

Since $g=g^{T}+\eta \otimes \eta $, the $\eta $-Einstein condition is
equivalent to
\begin{equation*}
\mathrm{Ric}^{T}=(A+2)g^{T}\quad \text{on }D=\ker \eta ,\qquad \text{that
is,}\qquad \rho ^{T}=(A+2)\omega ^{T},
\end{equation*}
where $\rho ^{T}=\mathrm{Ric}^{T}(\Phi \,\cdot ,\cdot )$ is the transverse
Ricci form; such a structure is also called \emph{transverse
K\"{a}hler--Einstein}. For $A=2n$ one has
\begin{equation*}
\mathrm{Ric}_{g}=2ng,
\end{equation*}
and the structure is \emph{Sasaki--Einstein}. In summary, $(M,g)$ is
Sasaki--Einstein if and only if the K\"{a}hler cone $(C(M),\overline{g})$ is
Ricci-flat, and if and only if the transverse K\"{a}hler structure of the Reeb
foliation $\mathcal{F}_{\xi }$ is K\"{a}hler--Einstein with
\begin{equation*}
\mathrm{Ric}^{T}=2(n+1)g^{T}.
\end{equation*}
We record the following curvature identities on a Sasakian manifold. By
\eqref{eqn-SR2},
\begin{equation}
R(e_{i},\xi ,e_{j},\xi )=g(e_{i},e_{j})=\delta _{ij}.  \label{eqn-SR3}
\end{equation}
Therefore
\begin{equation}
\mathrm{Ric}_{g}(\xi ,\xi )=\sum_{i=1}^{2n}R(e_{i},\xi ,e_{i},\xi )=2n.
\label{eqn-SR4}
\end{equation}
Consequently
\begin{equation*}
R^{T}=R_{g}+2n
\end{equation*}
and
\begin{equation}
\sum_{i=1}^{2n}\mathrm{Ric}_{g}(e_{i},e_{i})=R_{g}-2n=R^{T}-4n.  \label{eqn-SR6}
\end{equation}

\begin{definition}[{\cite{FOW,CLL}}]
Let $\mathcal X$ be a Reeb-invariant complex vector field whose projection
to the complexified normal bundle of the Reeb foliation is
of type $(1,0)$ and transversely holomorphic. It is \emph{Hamiltonian
holomorphic} with basic, possibly complex-valued, potential $f$ if
\begin{equation}
f=\sqrt{-1}\,\eta(\mathcal X),\qquad
\iota_{\mathcal X}\omega^T=\sqrt{-1}\,\overline\partial_B f.
\label{pr1}
\end{equation}
Put $X=\mathcal X+\overline{\mathcal X}=2\operatorname{Re}\mathcal X$.
The Sasakian structure is a \emph{Sasaki--Ricci soliton} if the real field
$X$ satisfies
\[
\mathrm{Ric}^T+\tfrac12\mathcal L_Xg^T=(A+2)g^T
\]
for a real constant $A$. The soliton is expanding, steady or shrinking
according as $A<-2$, $A=-2$ or $A>-2$, respectively.

In the \emph{gradient} case $f$ is real. With
$(\nabla^Tf)^{1,0}=\tfrac12(\nabla^Tf-\sqrt{-1}J\nabla^Tf)$,
\eqref{pr1} then gives
\[
\mathcal X=(\nabla^Tf)^{1,0}-\sqrt{-1}f\xi,\qquad
X=\nabla^Tf,\qquad
\tfrac12\mathcal L_Xg^T=\operatorname{Hess}^Tf.
\]
Thus the gradient equation is the full real equation
\eqref{eq:intro-soliton} with $\lambda=A+2$. Its mixed complex components
are
\[
R^T_{j\bar k}+f_{j\bar k}=(A+2)g^T_{j\bar k},\qquad
R_{j\bar k}+f_{j\bar k}=Ag^T_{j\bar k},
\]
where $R_{j\bar k}$ are the components of $\mathrm{Ric}_g|_D$;
the pure-type Hessian components vanish as well. In particular, the
field in the real Lie-derivative equation is $X$, not its complex
Hamiltonian lift $\mathcal X$.
\end{definition}

\begin{remark}
\label{rem:Dhom}For a constant $\tau >0$, the \emph{transverse $D$-homothetic
deformation} of $(M,\xi ,\eta ,\Phi ,g)$ is the Sasakian structure (see
\cite{BG})
\begin{equation}
\eta _{\tau }=\tau \eta ,\qquad \xi _{\tau }=\tau ^{-1}\xi ,\qquad \Phi _{\tau
}=\Phi ,\qquad g_{\tau }=\tau g+\tau (\tau -1)\eta \otimes \eta =\tau
g^{T}+\tau ^{2}\eta \otimes \eta .  \label{eq:pre-Dhom}
\end{equation}
Then $g_{\tau }^{T}=\tau g^{T}$, the transverse Levi-Civita connection is
unchanged, and $\mathrm{Ric}_{\tau }^{T}=\mathrm{Ric}^{T}$,
$\operatorname{Hess}_{\tau }^{T}f=\operatorname{Hess}^{T}f$ as
$(0,2)$-tensors, while $R^{T}$ and $|\nabla ^{T}f|^{2}$ are multiplied by
$\tau ^{-1}$. Consequently, if $(g,f)$ satisfies \eqref{eq:intro-soliton} with
constant $\lambda $, then $(g_{\tau },f)$ satisfies it with constant $\lambda
/\tau $. For $\lambda >0$ we may therefore choose either normalization in
\eqref{eq:intro-normal}. If $L$ is a Reeb period, then
$L_\tau=\tau L$, transverse lengths are multiplied by $\sqrt\tau$,
and $\lambda_\tau=\lambda/\tau$; the function $f$ is unchanged.
Consequently the Hamilton normalization \eqref{eq:pre-normal} is preserved,
as is the product $\lambda L$. A fixed $D$-homothety preserves completeness
and the regular,
quasi-regular or irregular type of the Reeb foliation.
\end{remark}

Throughout this paper the potential of a Sasaki--Ricci soliton is denoted by
$f$ and we write $\lambda =A+2$, while $\psi $ denotes test functions
(Appendix~\ref{app:reduced}) and
potentials of local K\"{a}hler quotients, in
particular of two-dimensional factors in \cite{PartII}. For a gradient
soliton, the full real equation \eqref{eq:intro-soliton} contains both
\begin{equation*}
R_{j\overline{k}}^{T}+f_{j\overline{k}}=\lambda g_{j\overline{k}}^{T}
\end{equation*}
and the vanishing of the $(2,0)+(0,2)$-part of $\operatorname{Hess}^{T}f$; the
latter says that $\nabla ^{T,1,0}f$ is transversely holomorphic, and it is
not implied by the $(1,1)$-equation alone. With $d^{c}u=-du\circ J$, the
full real equation gives the form equation
\begin{equation}
\rho ^{T}+\tfrac{1}{2}d_{B}d_{B}^{c}f=\lambda \omega ^{T},\qquad \rho
^{T}=\mathrm{Ric}^{T}(J\cdot ,\cdot ).  \label{eq:pre-form}
\end{equation}
The transverse trace and divergence of \eqref{eq:intro-soliton} give the
identities
\begin{equation}
R^{T}+\Delta _{B}f=2n\lambda ,\quad \nabla ^{T}R^{T}=2\mathrm{Ric}^{T}(\nabla
^{T}f,\cdot ),\quad R^{T}+|\nabla ^{T}f|^{2}-2\lambda f=\mathrm{const},
\label{eq:pre-identities}
\end{equation}
which are proved in Lemma~\ref{lem:identities}. For $\lambda >0$ we
normalize the additive constant in $f$ so that the constant in
\eqref{eq:pre-identities} is zero,
\begin{equation}
R^{T}+|\nabla ^{T}f|^{2}=2\lambda f,  \label{eq:pre-normal}
\end{equation}
and call the soliton \emph{normalized}. Under the normalization $\lambda
=2n+2$ the identities \eqref{eq:pre-identities} read, in terms of the
ambient curvature,
\begin{equation*}
\mathrm{Ric}_{g}|_{D}=2ng^{T}-\operatorname{Hess}^{T}f,\quad R_{g}+\Delta
_{B}f=4n^{2}+2n,\quad R_{g}+|\nabla ^{T}f|^{2}=4(n+1)f-2n.
\end{equation*}
In particular a gradient shrinking Sasaki--Ricci soliton with $\lambda =2n+2$
is Sasaki--Einstein if and only if $f$ is constant.

\begin{example}
\label{ex:compact}Unlike compact gradient K\"{a}hler--Ricci shrinker surfaces,
which are classified (see \cite{Cao} and the references therein), compact
five-dimensional gradient shrinking Sasaki--Ricci solitons are not classified.
Examples include:

\begin{enumerate}
\item the regular Sasaki--Einstein metrics on $\mathbb{S}^{5}$,
$\mathbb{S}^{2}\times \mathbb{S}^{3}$ and $\#k(\mathbb{S}^{2}\times
\mathbb{S}^{3})$, $3\leq k\leq 8$, obtained as circle bundles over the
K\"{a}hler--Einstein del Pezzo surfaces \cite{BG};

\item the toric Sasaki--Einstein metrics on $\#k(\mathbb{S}^{2}\times
\mathbb{S}^{3})$ for every $k\geq 1$, most of which are quasi-regular or
irregular \cite{FOW}; see \cite{BG} and the references therein;

\item the regular lifts to $\mathbb{S}^{2}\times \mathbb{S}^{3}$ and to
$\#2(\mathbb{S}^{2}\times \mathbb{S}^{3})$ of the Koiso--Cao shrinker on
$\mathbb{CP}^{2}\#\overline{\mathbb{CP}}^{2}$ \cite{Koi,Cao96} (see also
\cite{Cao}) and of the Wang--Zhu shrinker on
$\mathbb{CP}^{2}\#2\overline{\mathbb{CP}}^{2}$ \cite{WZ}, respectively.
\end{enumerate}
\end{example}

\subsection{Types of Reeb Foliations and Developments}

A Sasakian structure is \emph{quasi-regular} if all leaves of the Reeb
foliation $\mathcal{F}_{\xi }$ are compact, so that $\xi $ generates a locally
free circle action. It is \emph{regular} if this action is free, that is, the
Reeb flow is a free effective circle action with a common minimal period
$L>0$. We call the structure \emph{irregular} if some Reeb leaf is not closed.
On a noncompact manifold, there is also the possibility of a free proper
$\mathbb{R}$-action which has a smooth global quotient and is treated
separately. In that case every leaf is a closed but noncompact line, so the
structure is neither quasi-regular nor irregular; the phrase ``not
quasi-regular'' always includes this case.

In the quasi-regular case the leaf space $\mathcal{Z}=M/\mathbb{S}^{1}$ is a
K\"{a}hler orbifold and
\begin{equation*}
\pi :M\rightarrow \mathcal{Z}
\end{equation*}
is a Seifert circle orbibundle \cite{BG,Kol1}. In the regular case
$Z=M/\mathbb{S}^{1}$ is a smooth K\"{a}hler manifold and $\pi $ is a principal
circle bundle. In both cases
\begin{equation}
d\eta =2\pi ^{\ast }\omega _{Z},\qquad g=\pi ^{\ast }h_{Z}+\eta \otimes \eta
,\qquad f=\pi ^{\ast }f_{Z},  \label{eq:pre-BW}
\end{equation}
and \eqref{eq:intro-soliton} descends to the gradient K\"{a}hler--Ricci
shrinker equation
\begin{equation*}
\mathrm{Ric}_{Z}+\operatorname{Hess}_{Z}f_{Z}=\lambda h_{Z}
\end{equation*}
on the leaf space. If $L$ is the generic Reeb period, normalize the
circle generator and connection by
$U=(L/2\pi)\xi$ and $\theta=(2\pi/L)\eta$, so that $U$ has period
$2\pi$ and $\theta(U)=1$. Write $d\theta=\pi^*F_\theta$.
We use the curvature convention
\[
e_{\mathbb R}=[F_\theta/(2\pi)]=\frac2L[\omega_Z]
=\frac{4\pi}{\lambda L}c_1(Z)_{\mathbb R}.
\]
The last equality follows from
$\rho_Z+\tfrac12dd^cf_Z=\lambda\omega_Z$ and
$[\rho_Z]=2\pi c_1(Z)_{\mathbb R}$. In the quasi-regular case, these
identities are understood on $\mathcal Z$, in orbifold cohomology, with
$c_1^{\mathrm{orb}}(T\mathcal Z)$ in place of $c_1(Z)$.
Circle bundles over a smooth base are classified by the full integral
Euler class \cite{Ste,Hat}; the displayed curvature class records its
real image and does not determine any torsion component. The orbifold
version is described in \cite{Kol1,Kol2}.

In summary, without any regularity assumption, Lemma~\ref{lem:reeb} shows that
exactly one of the following holds.

\begin{enumerate}
\item[(i)] The closure of the Reeb flow in the isometry group of $(M,g)$ is a
compact torus $H$. This comprises
\begin{enumerate}
\item[(a)] the regular case;

\item[(b)] the quasi-regular but nonregular case;

\item[(c)] the irregular case, in which the Reeb flow is dense in a compact
torus $H$ of dimension greater than one. Here the orbit-closure space $M/H$ is
used only as a metric quotient; it need not be a transverse orbifold.
\end{enumerate}

\item[(ii)] The Reeb field $\xi $ generates a closed, free and proper
$\mathbb{R}$-action. Then $Z=M/\mathbb{R}$ is a smooth complete K\"{a}hler
manifold; this case includes the Heisenberg model.
\end{enumerate}

\begin{definition}
\begin{enumerate}
\item A \emph{local development} is an atlas of transverse charts on which the
transverse metric and the transverse complex structure are realized, together
with their transition maps. Such a development need not be injective into the leaf space. An \emph{unwrapped chart} is allowed to be a local immersion
rather than an embedded transverse neighborhood.

\item A \emph{complete smooth development} of a specified transverse
local model is a connected complete smooth K\"{a}hler manifold extending
that model: its metric and complex structure agree with the given germs
under local biholomorphic isometries. It also carries a single-valued
potential satisfying the same soliton equation and agreeing with the
specified local potential wherever that potential is part of the model.
Existence of such a complete extension is an additional global condition;
it does not assert that the original manifold has a global smooth leaf
space or a quotient map onto this extension.
\end{enumerate}
\end{definition}

\begin{proposition}
\label{prop:contactization}Let $(V,h,J,\psi )$ be a K\"{a}hler gradient
shrinker on a contractible open set,
\begin{equation*}
\mathrm{Ric}(h)+\operatorname{Hess}_{h}\psi =\lambda h,
\end{equation*}
with K\"{a}hler form $\omega (X,Y)=h(JX,Y)$. Choose a one-form $A$ with
$dA=2\omega $. On $(-\varepsilon ,\varepsilon )\times V$ put
\begin{equation*}
\eta =du+A,\qquad g=h+(du+A)^{2},\qquad \xi =\partial _{u},\qquad f=\psi .
\end{equation*}
These define a local Sasakian gradient shrinker with transverse K\"{a}hler
structure $(h,J)$. Two choices of $A$ differ locally by a Reeb-coordinate
change.
\end{proposition}

\begin{proof}
The horizontal lift of $J$ to $D=\ker \eta $ is compatible with $d\eta
=2\omega $. Its cone complex distribution is integrable: for horizontal
lifts $X^{H},Y^{H}$ of two $(1,0)$ fields, the vertical component of the
bracket is $-d\eta (X^{H},Y^{H})\xi =-2\omega (X,Y)\xi =0$, and invariant,
radially constant lifts commute with the radial and Reeb fields. The cone
two-form is $\frac{1}{2}d(r^{2}\eta )$ and is closed. Thus the cone is
K\"{a}hler. Transverse Ricci and Hessian tensors are exactly those of $V$.
Finally, if $A^{\prime }=A+db$, then $du+A^{\prime }=d(u+b)+A$, so the
change $u\mapsto u+b$ of the Reeb coordinate identifies the two contact
forms.
\end{proof}

We use the \emph{Heisenberg model} on $\mathbb{H}_{2}=\mathbb{C}^{2}\times
\mathbb{R}$, the global form of the contactization of Proposition
\ref{prop:contactization} for the Gaussian shrinker on $\mathbb{C}^{2}$ with
$A=\sum_{a}(x_{a}dy_{a}-y_{a}dx_{a})$:
\begin{equation}
\eta _{0}=dt+\sum_{a=1}^{2}(x_{a}dy_{a}-y_{a}dx_{a}),\quad g_{0}=|dz|^{2}+\eta
_{0}^{2},\quad f_{0}=\frac{\lambda
}{2}|z|^{2}.  \label{eq:heisenberg}
\end{equation}
This complete transverse Gaussian has
$\mathrm{Rm}^{T}=\mathrm{Ric}^{T}=R^{T}=0$ and $R_{g}=-4$. It is not the
ordinary Riemannian product
metric.

\section{Identities of Gradient Shrinking Sasaki--Ricci Solitons}

\label{sec:identities}

In this section we collect the fundamental identities of gradient shrinking
Sasaki--Ricci solitons (compare Chang--Li--Lin \cite{CLL}), construct the
exhaustion space on which basic cutoffs are compactly supported, and prove
the growth of the potential. The exhaustion space is constructed before any
cutoff function of the potential is used: the properness of the potential is
a consequence, not an ingredient, of the nonnegativity of the scalar
curvature.

\subsection{Soliton Identities}

The geometric hypothesis in this section and in Section~\ref{sec:curvature}
is
\begin{equation}
\mathrm{Ric}^{T}+\operatorname{Hess}^{T}f=\lambda h,\qquad \lambda
=2n+2,\qquad \xi
f=0,\qquad h=g|_{D},\quad D=\ker \eta .  \label{eq:soliton}
\end{equation}
It is a \emph{full real symmetric-tensor equation}: besides the
$(1,1)$-part of $\operatorname{Hess}^{T}f$ it prescribes the vanishing of
its $(2,0)+(0,2)$-part, which says that $\nabla^{T}f$ is a transversely real
holomorphic vector field. An equation for the $(1,1)$-part alone does not
imply \eqref{eq:soliton}. Recall that
\begin{equation*}
K=\mathrm{Rm}_{g}|_{D},\qquad S=\mathrm{Ric}_{g}|_{D},\qquad
T=\mathrm{Rm}^{T},\qquad Q=\mathrm{Ric}^{T},\qquad s=R^{T}.
\end{equation*}

\begin{remark}
We recall the normalizations of the Sasaki--Ricci flow. As in
Remark~\ref{rem:flow-normal}, complex components are taken by complex-bilinear
extension of the real tensors, in the frame $Z_{j}=\partial /\partial
z^{j}-\eta (\partial /\partial z^{j})\,\xi $ of $D\otimes \mathbb{C}$. With
this convention the flow of \cite{SWZ} reads
\begin{equation}
\frac{\partial }{\partial t}g_{j\overline{k}}^{T}=-R_{j\overline{k}}^{T},
\label{fe1}
\end{equation}
which is the real equation $\partial _{t}g^{T}=-\mathrm{Ric}^{T}$, and
\eqref{eq:marked-flow} is obtained from it by the change of time $t\mapsto
2t$. If the components are computed in a differently normalized unitary frame,
constant factors appear in \eqref{fe1}; all these normalizations differ only
by a linear change of the time variable. The curvature estimates of
Sections~\ref{sec:identities} and~\ref{sec:curvature} do not involve time,
while all statements about flows in Sections~\ref{sec:flows},
\ref{sec:splitting} and Appendix~\ref{app:reduced} refer to the normalization
\eqref{eq:marked-flow}; in particular the constants of the reduced geometry,
such as $n/\bar{\tau}$ in \eqref{eq:e3}, are those of \eqref{eq:marked-flow}.
Equation \eqref{fe1} is understood for the classes of $Z_{j}$ and
$\overline{Z_{k}}$ in the normal bundle $\nu (\mathcal{F}_{\xi })$ of the Reeb
foliation.
\end{remark}

Let $(M^{2n+1},\xi ,\eta ,\Phi ,g,f)$ satisfy \eqref{eq:soliton}, and let
$\{e_{j}\}_{j=1}^{2n}$ be a local orthonormal frame of $D$. The following
identities are proved in \cite{CLL}; we include the short derivation.

\begin{lemma}
\label{lem:identities}With the notation \eqref{eq:pre-KSTQ}, the following
identities hold.
\begin{equation}
\begin{array}{cl}
\left( \mathrm{i}\right) & S_{jk}=2nh_{jk}-f_{;jk}, \\
\left( \mathrm{ii}\right) & R+\Delta _{B}f=4n^{2}+2n, \\
\left( \mathrm{iii}\right) & R_{;i}=2S_{ij}f_{;j}+4f_{;i}, \\
\left( \mathrm{iv}\right) & R+|\nabla ^{T}f|^{2}=4(n+1)f+C_{1}, \\
\left( \mathrm{v}\right) & L_{f}R+2|S|^{2}-4(n-1)R-8n(n+1)=0.\end{array}
\label{n1}
\end{equation}
\end{lemma}

Equivalently, the purely transverse formulas are
\begin{equation}
s+\Delta _{B}f=2n\lambda ,\qquad s_{;i}=2Q_{ij}f_{;j},\qquad L_{f}s=2\lambda
s-2|Q|^{2}.  \label{eq:trans-scalar}
\end{equation}

Here $|S|^{2}=\sum_{i,j=1}^{2n}S_{ij}^{2}$ is the squared norm of the Ricci
curvature of $g$ restricted to $D$, and $C_{1}$ is a constant; under the
normalization \eqref{eq:canonical} of Lemma~\ref{lem:nonneg}, $C_{1}=-2n$.

\begin{proof}
Taking the transverse trace and divergence of \eqref{eq:soliton} gives
\begin{equation*}
s+\Delta _{B}f=2n\lambda ,\qquad \tfrac{1}{2}s_{;i}+(\Delta
_{B}f)_{;i}+Q_{ij}f_{;j}=0.
\end{equation*}
The second equation and the derivative of the first imply
$s_{;i}=2Q_{ij}f_{;j}$. Therefore
\begin{equation*}
(s+|\nabla^{T} f|^{2}-2\lambda f)_{;i}=2Q_{ij}f_{;j}+2(\lambda
h_{ij}-Q_{ij})f_{;j}-2\lambda f_{;i}=0.
\end{equation*}
This basic quantity is constant on a connected manifold. Differentiating
$s_{;i}=2Q_{ij}f_{;j}$ once more and contracting yields
\begin{equation*}
\Delta _{B}s=\left\langle \nabla^{T} s,\nabla^{T} f\right\rangle+2\lambda
s-2|Q|^{2}.
\end{equation*}
Now insert \eqref{eq:pre-ric}. In particular, $2\lambda (R+2n)-2(|S|^{2}+4R)$
is precisely (v).
\end{proof}

For an algebraic curvature tensor $U$, define
\begin{equation}
\mathscr B(U)_{ijkl}=\sum_{p,q}(U_{ipkq}U_{jplq}-U_{iplq}U_{jpkq}),\qquad
\mathscr P(U)_{ijkl}=\sum_{p,q}U_{pqij}U_{pqkl}.  \label{eq:B-P}
\end{equation}

\begin{proposition}
\label{prop:tensor} With convention \eqref{eqn-SR1}, the tensor
identities are
\begin{align}
\sum_i T_{ijkl;i} &=Q_{jl;k}-Q_{jk;l}=\sum_i T_{ijkl}f_{;i},  \label{eq:divT}
\\
-T_{ijkl}f_{;j}f_{;k}+\tfrac12s_{;il} &=Q_{il;p}f_{;p}+\lambda
Q_{il}-Q_{ip}Q_{pl},  \label{eq:second} \\
L_f T&=2\lambda T-2\mathscr B(T)-\mathscr P(T),  \label{eq:LT} \\
(L_f Q)_{il}&=2\lambda Q_{il}-2T_{iplq}Q_{pq}.  \label{eq:LQ}
\end{align}
In ambient-horizontal notation the last two equations become
\begin{align}
L_f K&=2\lambda(K+\mathring{A})-2\mathscr B(K+\mathring{A})-\mathscr
P(K+\mathring{A}),  \label{eq:LK} \\
(L_f S)_{il}&=-2K_{iplq}S_{pq}+(4n-6)S_{il}+8nh_{il}.  \label{eq:LS}
\end{align}
The ambient-horizontal version of \eqref{eq:second} is
\begin{equation}
  \label{eq:second-ambient}
\begin{split}
-(K+\mathring{A})_{ijkl}f_{;j}f_{;k}+\tfrac12R_{;il}
=S_{il;p}f_{;p}+(2n-2)S_{il}+4nh_{il}-S_{ip}S_{pl}.
\end{split}
\end{equation}
\end{proposition}

\begin{proof}
All computations can be made on a local transverse K\"{a}hler chart. For the
convention \eqref{eqn-SR1}, the contracted second Bianchi identity is
$T_{ijkl;i}=Q_{jl;k}-Q_{jk;l}$. Differentiating the soliton equation and
commuting covariant derivatives of $df$ gives
\begin{equation*}
Q_{jl;k}-Q_{jk;l}=T_{pjkl}f_{;p}.
\end{equation*}
This proves \eqref{eq:divT}. Next,
\begin{equation*}
\tfrac12s_{;il}=Q_{lp;i}f_{;p}+Q_{lp}f_{;pi}
=Q_{il;p}f_{;p}+T_{alip}f_{;a}f_{;p} +\lambda Q_{il}-Q_{lp}Q_{pi}.
\end{equation*}
Pair symmetry gives $T_{alip}f_{;a}f_{;p}=T_{ijkl}f_{;j}f_{;k}$, proving
\eqref{eq:second}.

For completeness, the standard elliptic curvature-operator computation on
this transverse chart gives
\begin{equation*}
L_{f}\mathbf{T}=2\lambda\mathbf{T}-2(\mathbf{T}^2+\mathbf{T}^\#).
\end{equation*}
Here $e_i\wedge e_j$, $i<j$, are orthonormal, and the components of
$\mathbf{T}^2$ and $\mathbf{T}^\#$ are respectively $\frac12\mathscr P(T)$ and
$\mathscr B(T)$. This is the curvature-operator identity of Petersen--Wylie
\cite[Lemma~2.1]{PW}, with its sign convention converted to \eqref{eqn-SR1};
it immediately gives \eqref{eq:LT}. Equivalently one
differentiates
\begin{equation*}
T_{ijkl;p}=T_{klpj;i}-T_{klpi;j},
\end{equation*}
uses \eqref{eq:divT}, substitutes $f_{;ij}=\lambda h_{ij}-Q_{ij}$, and
commutes the remaining derivatives. The Ricci-linear commutators cancel the
four terms from the Hessian of $f$, leaving the two negative quadratic
contractions in \eqref{eq:LT}. Taking the $i,k$ trace yields \eqref{eq:LQ};
in components the quadratic trace is
\begin{equation*}
\sum_i\{2\mathscr B(T)_{ijil}+\mathscr P(T)_{ijil}\}
=2\sum_{p,q}Q_{pq}T_{jplq}.
\end{equation*}
Since $\nabla^{T}\mathring{A}=0$, \eqref{eq:LK} follows. The tensor $S$ is
$\Phi$-invariant, so $\mathring{A}_{iplq}S_{pq}=3S_{il}$. Inserting $Q=S+2h$,
$T=K+\mathring{A}$ into \eqref{eq:LQ} then gives
\begin{equation*}
L_{f} S=2\lambda S+4\lambda h-2K(S)-6S-4Q
=-2K(S)+(2\lambda-10)S+(4\lambda-8)h,
\end{equation*}
which is \eqref{eq:LS}.
\end{proof}

\begin{lemma}
\label{lem:reeb}Let $(M,\xi ,\eta ,\Phi ,g)$ be a connected complete Sasakian
manifold.

\begin{enumerate}
\item Either the closure of the Reeb one-parameter group in
$\operatorname{Isom}(M,g) $ is a compact torus $H$,

\item or the Reeb group is a closed, free and proper $\mathbb{R}$-action.
\end{enumerate}

In the first case put $\mathcal{X}=M$ and choose a compact orbit $P=Hp$. In
the second case put $\mathcal{X}=Z=M/\mathbb{R}$, which is a smooth complete
K\"{a}hler manifold, and let $P$ be a point of $Z$.

In both cases the distance $r=d(P,\cdot)$ is proper on $\mathcal X$.
If, in addition, the Sasakian structure carries a basic potential $f$
satisfying the shrinking soliton equation \eqref{eq:soliton} with
$\lambda=2n+2$, then
\begin{equation}
L_{f}r\leq C-ar  \label{eq:radial-drift}
\end{equation}
outside a fixed neighborhood of $P$ with
\begin{equation*}
a=\begin{cases}
2n, & \mathcal{X}=M, \\
\lambda , & \mathcal{X}=Z.\end{cases}
\end{equation*}
For $\mathcal{X}=M$, the distance $r$ is basic. The distance inequality is
understood in the usual upper-barrier sense at its cut locus.
\end{lemma}

\begin{proof}
The closure of a one-parameter subgroup is a connected abelian Lie group. If
its projection onto a Euclidean factor is nonzero, convergence in the group
forces boundedness of the parameter; the subgroup is therefore closed and
isomorphic to $\mathbb{R}$. Otherwise its closure is a torus. A closed
subgroup of the isometry group acts properly on a complete Riemannian
manifold. A stabilizer for an $\mathbb{R}$-action is then a compact subgroup
of $\mathbb{R}$, hence trivial. The quotient metric is complete, and its
local K\"{a}hler structures glue.

In the compact-torus case all basic tensors under consideration are
$H$-invariant, so $r=d(P,\cdot)$ is basic and proper. A geodesic minimizing
from $P$ starts normal to $P$, hence normal to the Reeb field, which is
tangent to the orbit $P=Hp$; the Killing identity implies
$g(\dot\gamma,\xi)=0$ along it. At regular points of $r$, the ordinary
Bochner or Riccati formula gives
\begin{equation*}
\nabla r(L_{f} r)=-|\operatorname{Hess}_g r|^2
-(\mathrm{Ric}_g+\operatorname{Hess}_g
f)(\nabla r,\nabla r)\le-2n.
\end{equation*}
Here the last evaluation uses only the horizontal soliton equation. For
basic $f$ and horizontal $X$ one has $\operatorname{Hess}_g
f(X,X)=\operatorname{Hess}^{T}f(X,X)$, so \eqref{eq:soliton} and
\eqref{eq:pre-ric} give
\begin{equation*}
(\mathrm{Ric}_g+\operatorname{Hess}_g f)(X,X)=\lambda -2=2n
\end{equation*}
for a unit horizontal vector $X$. On a small fixed tubular boundary of the
compact submanifold $P$, $L_{f} r$ has a finite upper bound. Integration
along normal minimizing rays proves \eqref{eq:radial-drift}. At a cut point $x$
we use, as usual, the upper barrier $\varepsilon+d(\gamma(\varepsilon),\cdot)$,
where $\gamma$ is a minimizing geodesic from $P$ to $x$: it lies above $r$,
touches it at $x$ and is smooth near $x$, and the same second-variation
bound gives the upper-barrier inequality with constants independent of
$\varepsilon$. In the quotient case the argument is the
ordinary one, with $\mathrm{Ric}_Z+\operatorname{Hess}_Z f=\lambda h_Z$ and
$a=\lambda$.
\end{proof}

\subsection{Nonnegativity of the Transverse Scalar Curvature}

\begin{lemma}
\label{lem:nonneg}The transverse scalar curvature of any complete noncompact
gradient shrinking Sasaki--Ricci soliton $(M^{2n+1},\xi ,\eta ,\Phi ,g,f)$ is
nonnegative.

\begin{enumerate}
\item Either $R^{T}>0$ everywhere,

\item or $R^{T}\equiv 0$.
\end{enumerate}

Normalize the additive constant in $f$ by
\begin{equation}
R^{T}+|\nabla ^{T}f|^{2}=2\lambda f.  \label{eq:canonical}
\end{equation}
Then $f\geq 0$, and
\begin{equation}
L_{f}f=2n\lambda -2\lambda f,\qquad |\nabla ^{T}\sqrt{f}|\leq \sqrt{\lambda
/2}.  \label{eq:Lf-f}
\end{equation}
The second inequality holds in the Lipschitz sense also at zeros of $f$.
\end{lemma}

\begin{proof}
The proof uses only the local identities of Lemma~\ref{lem:identities} and
the radial barrier of Lemma~\ref{lem:reeb}; in particular, neither $f\geq 0$
nor the properness of $f$ is assumed. (For Riemannian shrinkers, the
nonnegativity of the scalar curvature is due to B.-L. Chen \cite{Che}.) A
fixed $D$-homothety reduces the proof to $\lambda =2n+2$. Put $s=R^{T}$ and
work on the space $\mathcal{X}$ of Lemma~\ref{lem:reeb}. In its $M$ case all
functions below are basic and $L_{f}$ on them is the ambient drift
Laplacian; in its $Z$ case it is the usual drift Laplacian of the complete
smooth quotient. The scalar identity and the transverse trace inequality
give
\begin{equation}
L_{f}s=2\lambda s-2|\mathrm{Ric}^{T}|^{2}\leq 2\lambda s-\frac{s^{2}}{n}.
\label{eq:scalar-ineq}
\end{equation}

Choose a smooth nonincreasing $\chi :[0,\infty )\rightarrow \lbrack 0,1]$
equal to one on $[0,1]$, positive on $[0,2)$ and zero on $[2,\infty )$, with
$|\chi ^{\prime \prime }|\leq C$ and $|\chi ^{\prime }|^{2}/\chi \leq C$
wherever $\chi >0$. Such a cutoff is obtained by squaring a smooth
nonincreasing cutoff. For large $L$ set $\zeta =\chi (r/L)$. Its support is
compact on $\mathcal{X}$ by properness of $r$, independently of the
potential. On the transition annulus, $L_{f}r\leq C_{0}-ar\leq 0$ for
sufficiently large $L$. Since $\chi ^{\prime }\leq 0$, the chain rule gives,
in the lower-barrier sense,
\begin{equation*}
L_{f}\zeta =\frac{\chi ^{\prime \prime }}{L^{2}}|\nabla r|^{2}+\frac{\chi
^{\prime }}{L}L_{f}r\geq -\frac{C}{L^{2}},\qquad \frac{|\nabla \zeta
|^{2}}{\zeta }\leq \frac{C}{L^{2}}.
\end{equation*}
These inequalities also hold where $\zeta $ is constant.

If $v=\zeta s$ is negative somewhere, it has a negative minimum at an
interior point where $\zeta>0$ and $s<0$. At a smooth minimum, $\nabla
s=-s\nabla\zeta/\zeta$ and $L_fv\ge0$, whence
\begin{equation*}
\begin{split}
0\le\zeta L_fv &\le 2\lambda\zeta v-\frac{v^2}{n}
+\left(L_f\zeta-2\frac{|\nabla\zeta|^2}{\zeta}\right)v \\
&\le-\frac{v^2}{n}+\frac{C|v|}{L^2}.
\end{split}
\end{equation*}
The last inequality uses $v<0$ and the lower bound for the quantity in
parentheses. Thus $v\ge-Cn/L^2$. At a cut point use the upper supporting
distances of Lemma~\ref{lem:reeb}, with their arbitrarily small barrier
error. Since $s<0$ near the minimum and $\chi$ decreases, the supported
product lies above $v$, touches it at the minimum, and has a minimum there
too. The same computation followed by letting the barrier error tend to zero
proves the identical bound. For any fixed point $x$, take $L>r(x)$ and let
$L\to\infty$ to get $s(x)\ge0$. If $\mathcal{X}$ is compact, a negative
minimum of $s$ directly contradicts \eqref{eq:scalar-ineq}.

For the zero alternative apply the strong maximum principle to $u=-s\leq 0$, which satisfies
\begin{equation*}
(L_{f}-2\lambda )u=2|\mathrm{Ric}^{T}|^{2}\geq 0,
\end{equation*}
with nonpositive zeroth-order coefficient. An interior maximum $u=0$ forces
$u$ to be identically zero on the connected space $\mathcal{X}$. Otherwise
$s>0$.

The Hamilton identity now permits the additive normalization
\eqref{eq:canonical}, which yields $f\ge0$ and $|\nabla^Tf|^2\le2\lambda f$.
Taking the trace of the soliton equation and subtracting $|\nabla^Tf|^2$
gives the first identity in \eqref{eq:Lf-f}. For $\varepsilon>0$,
\begin{equation*}
|\nabla^T\sqrt{f+\varepsilon}|^2
=\frac{|\nabla^Tf|^2}{4(f+\varepsilon)}\le\frac\lambda2.
\end{equation*}
Local uniform convergence as $\varepsilon\downarrow0$ proves the Lipschitz
assertion at zeros as well. Undoing the fixed $D$-homothety gives the
statements at every $\lambda>0$.
\end{proof}

\begin{lemma}
\label{lem:growth}At the normalization $\lambda =2n+2$, with $a$
and $r$ as in Lemma~\ref{lem:reeb}, there is a constant $C$ such that
\begin{equation}
\frac{(ar-C)_{+}^{2}}{2\lambda }\leq f\leq \frac{(\lambda r+C)^{2}}{2\lambda
}\qquad \text{on }\mathcal{X}.  \label{eq:potential-growth}
\end{equation}
In particular, the sublevel sets of $f$ are compact on $\mathcal X$.
In the compact-torus case this is properness on $M$; in the proper-line
case it is properness of the descended potential on $Z=M/\mathbb R$.
We use \emph{transverse infinity} to mean $f\to\infty$, equivalently
escape from all compact subsets of $\mathcal X$. In the proper-line case,
this concerns the projected points in $Z$, not escape along a Reeb fiber.
\end{lemma}

For Riemannian shrinkers the two-sided quadratic estimate of the potential
is due to Cao--Zhou \cite{CaZ}; the proof below adapts the second-variation
argument used there (compare also \cite[Lemma~5.1]{Nab}).

\begin{proof}
We give an argument for both branches which uses the already proved
nonnegativity of $R^T$, not nonnegativity of the ambient scalar curvature.
Let $d$ be the Riemannian dimension of $\mathcal X$: $d=2n+1$ on $M$
and $d=2n$ on $Z$. Write $\widehat g=g$ or $h_Z$, respectively. The
Hamilton identity gives
\[
 f\ge0,\qquad |\nabla f|\le\sqrt{2\lambda f},\qquad
 |\nabla\sqrt f|\le\sqrt{\lambda/2}
\]
on $\mathcal X$, the last inequality in the Lipschitz sense. The value
of $f$ is constant on the compact orbit $P$ because $f$ is basic and
$H$ is the closure of the Reeb flow. Thus, along a minimizing geodesic
from $P$, the last inequality immediately gives
$\sqrt{2\lambda f(x)}\le\lambda r(x)+C$.

For the lower bound take $r=r(x)>2$ and a unit-speed minimizing
geodesic $\gamma:[0,r]\to\mathcal X$ from $P$ to $x$. In the $M$ branch
it is horizontal by Lemma~\ref{lem:reeb}. Put $F(t)=f(\gamma(t))$ and
$\rho(t)=\mathrm{Ric}_{\widehat g}(\dot\gamma,\dot\gamma)$.
The full horizontal soliton equation gives, in the two branches,
\[
 \rho(t)=a-F''(t),\qquad
 a=\lambda-2=2n\quad(M),\qquad a=\lambda\quad(Z).
\]
Let $\phi(t)=t$ on $[0,1]$, $\phi=1$ on $[1,r-1]$, and
$\phi(t)=r-t$ on $[r-1,r]$. Applying the fixed-endpoint index form to
$\phi E_j$ for a parallel orthonormal frame normal to $\dot\gamma$ gives
\[
 \int_0^r\phi^2\rho\,dt\le(d-1)\int_0^r|\phi'|^2dt=2(d-1).
\]
The first-end defect $\int_0^1(1-t^2)\rho(t)dt$ is bounded above
uniformly: the radius-one neighborhood of the compact set $P$ is compact.
The remaining defect has the exact expression
\[
 \begin{split}
 \int_{r-1}^r\bigl(1-(r-t)^2\bigr)\rho(t)dt
 &=\frac{2a}{3}-F'(r)
       +2\int_{r-1}^r(r-t)F'(t)dt.
 \end{split}
\]
Here integration by parts has no contribution at $r-1$.
Since $\int_0^r\rho=ar-F'(r)+F'(0)$ and $F'(0)$ is uniformly bounded
on $P$, the terms $F'(r)$ cancel, leaving
\[
 ar-C\le 2\int_{r-1}^r(r-t)F'(t)dt.
\]
The gradient and Lipschitz estimates imply
\[
 F'(t)\le\sqrt{2\lambda F(t)}
 \le\sqrt{2\lambda F(r)}+\lambda(r-t).
\]
Consequently $ar-C\le\sqrt{2\lambda F(r)}+2\lambda/3$.
Absorb the last constant into $C$ and square the positive part.
Enlarging $C$ also covers $r\le2$ and proves both inequalities in
\eqref{eq:potential-growth}. Since $a>0$ and $r$ is proper on
$\mathcal X$, so is $f$. In the proper-line branch this assertion is
on $Z$, not along the noncompact Reeb fibers in $M$.
\end{proof}

\begin{proposition}
\label{prop:flat}If $R^{T}\equiv 0$, then $\mathrm{Ric}^{T}\equiv 0$ and
$\mathrm{Rm}^{T}\equiv 0$. If $R^{T}>0$, then $f>0$ and, for some $c>0$,
\begin{equation}
R^{T}f\geq c>0.  \label{eq:sf}
\end{equation}
\end{proposition}

For Riemannian shrinkers the analogous lower bound $R\ge c/f$ is due to
Chow--Lu--Yang \cite{CLY}.

\begin{proof}
If $R^{T}=0$, its drift equation implies $Q=0$ and hence $(\nabla
^{T})^{2}f=\lambda h$. The vector field $\nabla ^{T}f$ is complete by the
upper growth estimate, and its flow $F_{t}$ satisfies $F_{t}^{\ast
}h=e^{2\lambda t}h$ transversely and $f(F_{t}x)=e^{2\lambda t}f(x)$. Its
backward trajectories, considered on $\mathcal{X}$ (on $Z$ in the
proper-line case), remain in a compact sublevel. Curvature is bounded on
that sublevel, whereas homothety gives $|T|(F_{t}x)=e^{-2\lambda t}|T|(x)$.
Letting $t\rightarrow -\infty $ proves $T(x)=0$.

Now assume $s=R^{T}>0$ and put $m=2n$. The inverse powers of $f$ used below
are smooth: the normalized Hamilton identity \eqref{eq:canonical} gives
\begin{equation*}
2\lambda f=s+|\nabla^Tf|^2>0
\end{equation*}
at every point of $\mathcal{X}$. Moreover $f$ is proper on $\mathcal{X}$ by
Lemma~\ref{lem:growth}. Consequently $f$ has a positive minimum $f_*>0$ on
$\mathcal{X}$, and every compact $f$-sublevel has a positive minimum of $s$.

For $p>0$, the chain rule, $L_ff=m\lambda-2\lambda f$, and
$|\nabla^Tf|^2=2\lambda f-s$ give the \emph{exact} identity
\begin{equation}
L_f(f^{-p})=2p\lambda f^{-p} +p\lambda(2p+2-m)f^{-p-1} -p(p+1)s f^{-p-2}.
\label{eq:inversepowers-general}
\end{equation}
In particular,
\begin{align*}
L_f(f^{-1})&=2\lambda f^{-1}+(4-m)\lambda f^{-2}-2sf^{-3}, \\
L_f(f^{-2})&=4\lambda f^{-2}+(12-2m)\lambda f^{-3}-6sf^{-4}.
\end{align*}
Let $F=4m$ and choose
\begin{equation*}
0<c<\min_{\{f\le F\}} \frac{s}{f^{-1}+mf^{-2}}.
\end{equation*}
The sublevel is nonempty: at a minimum of $f$, $\nabla f=0$ and
$\Delta_Bf\ge0$, so $L_ff=m\lambda-2\lambda f\ge0$ and $\min f\le m/2<F$. The
minimum in the displayed choice is strictly positive, since the sublevel is
compact and $s,f>0$ there.

Set $w=s-cf^{-1}-mcf^{-2}$. Using $(L_{f}-2\lambda )s=-2|Q|^{2}\leq 0$ and
$0<s\leq 2\lambda f$, the exact inverse-power formulas yield, for $f\geq F$,
\begin{align*}
(L_{f}-2\lambda )w& \leq -(4-m)\lambda cf^{-2}+2csf^{-3}-2m\lambda
cf^{-2}-m(12-2m)\lambda cf^{-3}+6mcsf^{-4} \\
& \leq -m\lambda cf^{-2}+2m^{2}\lambda cf^{-3}\leq -\frac{m\lambda
c}{2}f^{-2}<0.
\end{align*}
By construction $w\geq 0$ on $\{f\leq F\}$. Also $w\geq -c/f-mc/f^{2}$, so
\begin{equation*}
\liminf_{f\rightarrow \infty }w\geq 0.
\end{equation*}
If $w$ were negative anywhere on $\{f>F\}$, a negative minimum would be
attained in a compact subset of that open region: choose a negative value
$-\delta $ and use $f\rightarrow \infty $ to exclude the end. At the
minimizing point $L_{f}w\geq 0$, whereas $(L_{f}-2\lambda )w<0$ implies
$L_{f}w<2\lambda w<0$, a contradiction. Thus $w\geq 0$ everywhere, and
\begin{equation*}
sf\geq c\left( 1+\frac{m}{f}\right) \geq c>0.
\end{equation*}
If $\mathcal{X}$ is compact, the same conclusion follows directly from the
positive minimum of $sf$.
\end{proof}

\section{Boundedness of the Curvature Tensor}

\label{sec:curvature}

In this section $n=2$, so the transverse real dimension is four and $\lambda
=6$. Assume
\begin{equation*}
R^{T}\leq B<\infty .
\end{equation*}
Since $R^{T}=R+4$ by \eqref{eq:pre-ric}, this is equivalent to an upper
bound for the ambient scalar curvature $R$. The flat alternative has already
been settled, so all ratios below are used only when $R^{T}>0$. All
constants in this section may depend on $B$, on the constant $c$ in
\eqref{eq:sf}, and on the transverse geometry and the normalized potential
on a fixed compact sublevel in $\mathcal{X}$; they are independent of the
cutoff parameter $r$ used below.

\begin{theorem}
\label{thm:curvature} Let $(M^{5},\xi ,\eta ,\Phi ,g,f)$ be connected,
complete, and satisfy \eqref{eq:soliton}. If $R^{T}$ is bounded above, then
the transverse curvature tensor $\mathrm{Rm}^{T}$ and the curvature tensor
$\mathrm{Rm}_{g}$ of $g$ are bounded. In the nonflat transverse case there is
a constant $C>0$ such that
\begin{equation}
|\mathrm{Rm}^{T}|\leq CR^{T}\quad \text{and}\quad |\mathrm{Rm}_{g}|\leq
C(1+R^{T})=C(5+R).  \label{eq:final-transverse}
\end{equation}
The constant $C$ depends on $B$ and the transverse geometry and normalized
potential on a fixed compact sublevel in $\mathcal{X}$.
\end{theorem}

The tensor identities below use the full real transverse soliton equation.
The scalar barrier \eqref{eq:sf} holds on the whole space $\mathcal{X}$. The
local curvature computations follow the four-dimensional strategy of
Munteanu--Wang \cite{MW1}; the points at which the argument differs from
\cite{MW1} are collected in Remark~\ref{rem:MW-comparison}. We keep $\lambda $
in the calculations.

\textbf{Standing conventions.} If $R^{T}\equiv 0$, then
$\mathrm{Rm}^{T}\equiv 0$ by Proposition~\ref{prop:flat} and
$\mathrm{Rm}_{g}$ is bounded by \eqref{eq:pre-Rmg}. From now on we therefore
assume that $s=R^{T}>0$ on $\mathcal{X}$ and that $f$ is normalized by
\eqref{eq:canonical}. Then $f>0$, the sublevel sets of $f$ are compact in
$\mathcal{X}$ (Lemma~\ref{lem:growth}), and $sf\geq c>0$
(Proposition~\ref{prop:flat}). The functions $s$, $f$ and the norms of the
transverse tensors $T$, $Q$, $\nabla ^{T}T,\ldots $ are basic. On
$\mathcal{X}=M$ the operator $L_{f}$ acts on basic functions as the ambient
drift Laplacian $\Delta _{g}-\langle \nabla f,\nabla \,\cdot \,\rangle $ by
\eqref{eq:pre-Lf}, and on $\mathcal{X}=Z$ it is the drift Laplacian of the
quotient. In both cases it is an elliptic operator on the manifold
$\mathcal{X}$, and its value on such a function is given by the transverse
formulas of Section~\ref{sec:identities}. This is what allows us to apply
the maximum principle at interior extrema on $\mathcal{X}$ and, later, to
integrate by parts on $\mathcal{X}$.

\subsection{The Four-Dimensional Algebraic Estimate}

Write $\Omega (t)=\{f\leq t\}\subset \mathcal{X}$. By \eqref{eq:canonical},
$|\nabla ^{T}f|^{2}=2\lambda f-s\geq 2\lambda f-B$. Hence outside $\Omega
(t_{0})$, with $t_{0}=t_{0}(B)\geq B/\lambda $ sufficiently large,
\eqref{eq:canonical} implies $|\nabla ^{T}f|\geq c\sqrt{f}$; the constant
$t_{0}$ will be enlarged a finite number of times below. On
$\mathcal{X}\setminus \Omega (t_{0})$ choose $e_{4}=\nabla ^{T}f/|\nabla
^{T}f|$ and complete it to an orthonormal frame $\{e_{1},\ldots ,e_{4}\}$ of
$D$. By \eqref{eq:divT}, $|\nabla ^{T}f|\,T_{4jkl}=Q_{jl;k}-Q_{jk;l}$, so
from the curvature symmetries every component with at least one index equal
to four satisfies
\begin{equation*}
|T_{ijk4}|\leq C\frac{|\nabla ^{T}Q|}{\sqrt{f}}.
\end{equation*}
An algebraic curvature tensor on a four-dimensional vector space is
controlled by its Ricci contraction and the components with at least one
index equal to four. Indeed, restrict to $e_{4}^{\perp }$: the restriction
is an algebraic curvature tensor on a three-dimensional space, its
three-dimensional Ricci tensor is $P_{ab}=Q_{ab}-T_{a4b4}$, and a
three-dimensional algebraic curvature tensor is determined by its Ricci
tensor, namely, for $a,b,c,d\in \{1,2,3\}$,
\begin{equation*}
T_{abcd}=P_{ac}\delta _{bd}+P_{bd}\delta _{ac}-P_{ad}\delta _{bc}-P_{bc}\delta
_{ad}-\tfrac{1}{2}(\operatorname{tr}P)(\delta _{ac}\delta _{bd}-\delta
_{ad}\delta _{bc}).
\end{equation*}
Thus
\begin{equation}
|T|\leq C\left( |Q|+\frac{|\nabla ^{T}Q|}{\sqrt{f}}\right) \quad \text{on
}\mathcal{X}\setminus \Omega (t_{0}).  \label{eq:4-algebra}
\end{equation}
This is a pointwise algebraic argument. In contrast with
\cite[Proposition~1.1]{MW1} it does not use the Gauss equation of the level
sets of $f$, so that no term coming from their second fundamental form
appears. For $K$, the valid estimate has an additional constant:
\begin{equation*}
|K|\leq C(1+|Q|+|\nabla ^{T}Q|/\sqrt{f}).
\end{equation*}

\subsection{Bounded Transverse Ricci Curvature}

We will repeatedly use $\chi =\zeta (f/r)$, where $r\geq 2t_{0}$ is a large
parameter (not to be confused with the distance function of
Lemma~\ref{lem:reeb}), $0\leq \zeta \leq 1$ is smooth, $\zeta =1$ on $[1,2]$,
and $\zeta $ is supported in $(1/2,3)$. Thus $\chi =1$ on $\Omega
(2r)\setminus \Omega (r)$, and the support of $\chi $ is contained in
$\Omega (3r)\setminus \Omega (r/2)\subset \mathcal{X}\setminus \Omega
(t_{0})$. Since $\nabla ^{T}\chi =r^{-1}\zeta ^{\prime }\nabla ^{T}f$ and
$L_{f}\chi =r^{-1}\zeta ^{\prime }L_{f}f+r^{-2}\zeta ^{\prime \prime }|\nabla
^{T}f|^{2}$, equations \eqref{eq:canonical}-\eqref{eq:Lf-f} imply
\begin{equation}
|\nabla ^{T}\chi |^{2}\leq C/r,\qquad |L_{f}\chi |\leq C,  \label{eq:cutoff}
\end{equation}
with $C$ independent of $r$. Their supports are compact on $\mathcal{X}$ by
Lemma~\ref{lem:growth}. If $u\geq 0$ and $G=\chi ^{2}u$ has a positive
maximum, then at that maximum $\chi \nabla ^{T}u=-2u\nabla ^{T}\chi $ and
$L_{f}G\leq 0$, so that expanding $L_{f}(\chi ^{2}u)$ gives
\begin{equation}
0\geq \chi ^{2}L_{f}G=\chi ^{4}L_{f}u+2\chi GL_{f}\chi -6G|\nabla ^{T}\chi
|^{2}.  \label{eq:cutoff-max}
\end{equation}
Consequently an inequality
\begin{equation*}
L_{f}u\geq a_{1}u^{2}-a_{2}u^{3/2}-a_{3}u-a_{4},
\end{equation*}
with $a_{1}>0$ and fixed constants $a_{2},a_{3},a_{4}\geq 0$ outside $\Omega
(t_{0})$, bounds $u$ uniformly on the annuli $\Omega (2r)\setminus \Omega
(r) $, $r\geq 2t_{0}$. Indeed, by \eqref{eq:cutoff}, \eqref{eq:cutoff-max}
and $0\leq \chi \leq 1$, at a positive maximum of $G$
\begin{equation*}
a_{1}G^{2}\leq a_{2}G^{3/2}+(a_{3}+C)G+a_{4},
\end{equation*}
which bounds $\sup G$ by a constant independent of $r$, and $u=G$ on $\Omega
(2r)\setminus \Omega (r)$. The compact remaining region is handled
separately. Norm inequalities at tensor zeros are understood in the usual
weak Kato sense, or by replacing a norm by $(|U|^{2}+\varepsilon )^{1/2}$
and passing to the limit.

Let $Y=|Q|^{2}$. Equation \eqref{eq:LQ} gives the exact identity
\begin{equation}
L_{f}Y=2|\nabla ^{T}Q|^{2}+4\lambda Y-4T_{ipjq}Q_{ij}Q_{pq}.  \label{eq:LY}
\end{equation}
Let $0<a<1$ and set $u=Ys^{-a}$, a smooth function because $s>0$. By the
product rule and $L_{f}s=2\lambda s-2Y$ (see \eqref{eq:trans-scalar}),
\begin{equation*}
\begin{split}
L_{f}u={}& s^{-a}L_{f}Y-2a\lambda u+2aY^{2}s^{-a-1} \\
& +a(a+1)s^{-a-2}Y|\nabla ^{T}s|^{2}-4as^{-a-1}\left\langle Q,\nabla
_{\nabla ^{T}s}^{T}Q\right\rangle .
\end{split}
\end{equation*}
The gradient terms are absorbed as follows:
\begin{align*}
& 2s^{-a}|\nabla ^{T}Q|^{2}-4as^{-a-1}\left\langle Q,\nabla _{\nabla
^{T}s}^{T}Q\right\rangle +a(a+1)s^{-a-2}Y|\nabla ^{T}s|^{2} \\
& \hspace{15mm}\geq \frac{2(1-a)}{1+a}s^{-a}|\nabla ^{T}Q|^{2}, \\
& \frac{C}{\sqrt{f}}s^{-a}Y|\nabla ^{T}Q|\leq \frac{1-a}{1+a}s^{-a}|\nabla
^{T}Q|^{2}+\frac{C_{a}}{f}s^{-a}Y^{2}.
\end{align*}
For the first inequality put $x=|\nabla ^{T}Q|$ and $y=|Q|\,|\nabla
^{T}s|/s$. By the Cauchy--Schwarz inequality the left side is at least
$s^{-a}\{2x^{2}-4axy+a(a+1)y^{2}\}$, and
\begin{equation*}
2x^{2}-4axy+a(a+1)y^{2}-\frac{2(1-a)}{1+a}x^{2}=\frac{a}{1+a}\bigl(2x-(1+a)y\bigr)^{2}\geq
0.
\end{equation*}
The second inequality is Young's inequality; it spends only half of the
positive gradient term from the first inequality, leaving a nonnegative
remainder. Constants $C_{a}$ may depend on the fixed exponent $a\in (0,1)$;
here $C_{a}=C^{2}(1+a)/(4(1-a))$. The remaining reaction term from
$-aYs^{-a-1}L_{f}s$ is $2aY^{2}s^{-a-1}-2a\lambda u$, and by
\eqref{eq:4-algebra} the curvature term in \eqref{eq:LY} satisfies
\begin{equation*}
-4s^{-a}T_{ipjq}Q_{ij}Q_{pq}\geq
-Cs^{-a}Y^{3/2}-\frac{C}{\sqrt{f}}s^{-a}Y|\nabla ^{T}Q|\qquad \text{on
}\mathcal{X}\setminus \Omega (t_{0}).
\end{equation*}
Since $Y^{2}s^{-a-1}=u^{2}s^{a-1}$, $Y^{2}s^{-a}=u^{2}s^{a}$ and
$Y^{3/2}s^{-a}=u^{3/2}s^{a/2}$, the product rule, \eqref{eq:4-algebra}, and
Young's inequality therefore give, on $\mathcal{X}\setminus \Omega (t_{0})$,
\begin{equation}
L_{f}u\geq \left( 2a-\frac{Cs}{(1-a)f}\right)
u^{2}s^{a-1}-Cu^{3/2}s^{a/2}+(4-2a)\lambda u.  \label{eq:quotient-a}
\end{equation}
Take $a=1/2$ and enlarge $t_{0}(B)$ so the leading bracket is at least $1/2$.
Since $s\leq B$, we have $s^{a-1}\geq B^{-1/2}$ and $s^{a/2}\leq B^{1/4}$,
and the inequality has the form needed for \eqref{eq:cutoff-max}. Therefore
$Y/\sqrt{s}\leq C$ outside a compact sublevel, and in particular, since
$s\leq B$ and $Q$ is bounded on compact sets,
\begin{equation}
|Q|\leq C\quad \text{on }\mathcal{X}.  \label{eq:Q-bounded}
\end{equation}

\subsection{Bounded Transverse Curvature}

Put $W=|T|$. Since $L_{f}|T|^{2}=2|\nabla ^{T}T|^{2}+2\langle
T,L_{f}T\rangle $ and $|\nabla ^{T}|T||\leq |\nabla ^{T}T|$, equation
\eqref{eq:LT} and Kato's inequality give
\begin{equation}
L_{f}W\geq 2\lambda W-CW^{2}.  \label{eq:LW}
\end{equation}
Equations \eqref{eq:4-algebra}, \eqref{eq:LY} and \eqref{eq:Q-bounded}
imply, for large $f$,
\begin{equation*}
W^{2}\leq C+C\frac{|\nabla ^{T}Q|^{2}}{f},\qquad L_{f}Y\geq |\nabla
^{T}Q|^{2}-C.
\end{equation*}
For the second bound note that $L_{f}Y\geq 2|\nabla ^{T}Q|^{2}-C|T|\geq
2|\nabla ^{T}Q|^{2}-C(1+|\nabla ^{T}Q|/\sqrt{f})$. Writing
$-CW^{2}=W^{2}-(C+1)W^{2}$ in \eqref{eq:LW} and using the first of these
bounds, we obtain
\begin{equation*}
L_{f}W\geq 2\lambda W+W^{2}-C-C|\nabla ^{T}Q|^{2}/f\geq W^{2}-C-C|\nabla
^{T}Q|^{2}/f.
\end{equation*}
Adding the second bound and increasing $t_{0}$ shows, for $\mathcal{H}=W+Y$,
\begin{equation*}
L_{f}\mathcal{H}\geq W^{2}-C\geq \tfrac{1}{2}\mathcal{H}^{2}-C,
\end{equation*}
where the last step uses that $Y$ is bounded. The annular maximum principle
bounds $T$ globally.

The local differentiated soliton-curvature identity is
\begin{equation}
  \label{eq:deriv}
L_f(\nabla^{T} T)=3\lambda\nabla^{T} T+T*\nabla^{T} T,
\end{equation}
where $*$ is a fixed finite sum of universal tensor contractions. In
particular,
\begin{equation}
  \label{eq:deriv-norm}
L_f|\nabla^{T} T|^2\ge2|(\nabla^{T})^2T|^2+6\lambda|\nabla^{T} T|^2
-C|T|\,|\nabla^{T} T|^2.
\end{equation}
To see the coefficient, commute $\nabla^T_q$ past the rough Laplacian and
drift:
\begin{equation*}
L_f(\nabla^T_qT)=\nabla^T_q(L_f T)+\lambda\nabla^T_qT+T*\nabla^{T} T.
\end{equation*}
The terms involving a derivative of Ricci or a curvature contraction with
$\nabla^{T} f$ combine by \eqref{eq:divT}; the Hessian of $f$ contributes
$\lambda h-Q$, cancelling the Ricci action on the derivative index.
Differentiating \eqref{eq:LT} supplies the remaining $2\lambda\nabla^T_qT$.
This is also the calculation in \cite[Proposition~2.2]{MW1}, whose
normalization $\lambda=1/2$ gives the coefficient $3$ in the squared-norm
equation.

Now $\Xi=|\nabla^{T} T|$ satisfies $L_f \Xi\ge-C\Xi$ and $L_f W^2\ge2\Xi^2-C$,
because $T$ is bounded. Hence $L_{f}(\Xi +W^{2})\geq 2\Xi ^{2}-C\Xi -C\geq
\frac{1}{2}(\Xi +W^{2})^{2}-C$, and applying \eqref{eq:cutoff-max} to
$\Xi+W^2$ gives
\begin{equation}
  \label{eq:T-first-bounded}
|T|+|\nabla^{T} T|\le C\quad\text{on }\mathcal{X}.
\end{equation}
All tensor calculations here are transverse, so no ambient Reeb commutators
are discarded.

Once $T$ is bounded, a fixed exponent suffices. Direct calculation for $u=Y/s
$ (the case $a=1$ of the product rule above) gives
\begin{equation}
L_{f}u=\frac{2}{s}\left\vert \nabla ^{T}Q-Q\otimes \nabla ^{T}\log
s\right\vert ^{2}+2\lambda u-\frac{4}{s}T_{ipjq}Q_{ij}Q_{pq}+2u^{2}\geq
2u^{2}-Cu.  \label{eq:Y-over-s}
\end{equation}
Thus \eqref{eq:cutoff-max} implies
\begin{equation}
|Q|^{2}\leq Cs.  \label{eq:Qsqrt}
\end{equation}
Since $sf\geq c$, equations \eqref{eq:4-algebra} and
\eqref{eq:T-first-bounded} now give
\begin{equation}
|T|\leq C\sqrt{s}.  \label{eq:Tsqrt}
\end{equation}
Indeed, $|\nabla ^{T}Q|/\sqrt{f}\leq C/\sqrt{f}\leq C\sqrt{s}$ outside
$\Omega (t_{0})$. The compact sublevel causes no problem because $s$ has a
positive minimum there.

We next record the measure needed for a global maximum principle. Let
\begin{equation}
d\mu =\begin{cases}
d\operatorname{vol}_{g}, & \mathcal{X}=M, \\
d\operatorname{vol}_{h_{Z}}, & \mathcal{X}=Z,\end{cases}\qquad d\nu
_{f}=e^{-f}d\mu .  \label{eq:mu}
\end{equation}
This is the Riemannian measure on the exhaustion space $\mathcal X$.
When $\mathcal X=M$ it has dimension $2n+1$ and is distinct from the
$2n$-dimensional quotient measure used in setting (S). For compactly
supported basic functions on $M$, or arbitrary compactly supported
functions on $Z$,
\begin{equation}
\int vL_{f}u\,d\nu _{f}=-\int \left\langle \nabla ^{T}v,\nabla
^{T}u\right\rangle \,d\nu _{f}.  \label{eq:ibp}
\end{equation}
Here $u$ is smooth and it suffices that $v$ be Lipschitz. The measure $\nu
_{f}$ has finite mass and all polynomial $f$-moments. Indeed, bounded $T$
bounds the full ambient Ricci tensor on $M$, or the quotient Ricci tensor on
$Z$. Ordinary volume comparison gives at most exponential growth of metric
balls, $\mu (B(P,k))\leq Ce^{Ck}$, while \eqref{eq:potential-growth} gives a
quadratic lower bound for $f$, so that for every $p\geq 0$
\begin{equation*}
\int f^{p}\,d\nu _{f}\leq C_{p}\sum_{k\geq 0}\mu
\bigl(B(P,k+1)\bigr)\,(1+k)^{2p}\exp \Bigl(-\frac{(ak-C)_{+}^{2}}{2\lambda
}\Bigr)<\infty .
\end{equation*}
This proves Gaussian integrability. No volume comparison on a transverse
orbifold is used.

\begin{lemma}
\label{lem:L2} Let $v$ be a smooth function on the exhaustion space
$\mathcal X$, basic when $\mathcal X=M$. Suppose
$v\in L^2(\mathcal X,d\nu_f)$ and
\begin{equation*}
L_fv\geq av-b
\end{equation*}
for constants $a>0$ and $b\in\mathbb R$. Then
\begin{equation*}
v\leq b/a.
\end{equation*}
\end{lemma}

\begin{proof}
Let $w=(v-b/a)_{+}$ and use the test function $\theta _{r}^{2}w$, where
$\theta _{r}=\vartheta (f/r)$ for a fixed smooth nonincreasing $\vartheta $
with $\vartheta =1$ on $[0,1]$ and $\vartheta =0$ on $[2,\infty )$; thus
$\theta _{r}=1$ on $\Omega (r)$, $\theta _{r}$ vanishes outside $\Omega
(2r)$, and $|\nabla ^{T}\theta _{r}|^{2}\leq C/r$. On $\{w>0\}$ we have
$L_{f}v\geq a(v-b/a)=aw$ and $\nabla ^{T}v=\nabla ^{T}w$. Integration by
parts \eqref{eq:ibp}, justified first by a smooth convex approximation to
the positive part, gives
\begin{equation*}
a\int \theta _{r}^{2}w^{2}\,d\nu _{f}\leq -\int \theta _{r}^{2}|\nabla
^{T}w|^{2}\,d\nu _{f}-2\int \theta _{r}w\left\langle \nabla ^{T}\theta
_{r},\nabla ^{T}w\right\rangle \,d\nu _{f}\leq \int |\nabla ^{T}\theta
_{r}|^{2}w^{2}\,d\nu _{f}.
\end{equation*}
The right side is at most $(C/r)\int w^{2}d\nu _{f}$ and tends to zero as
$r\rightarrow \infty $, since $w\leq |v|+|b|/a$ and $\nu _{f}$ is finite, so
that $w\in L^{2}(d\nu _{f})$. Therefore $w=0$.
\end{proof}

Set $v=|\nabla ^{T}T|^{2}/s$. The gradient terms in its product rule form a
nonnegative square after using \eqref{eq:deriv-norm}: with
$L_{f}(s^{-1})=-2\lambda s^{-1}+2|Q|^{2}s^{-2}+2s^{-3}|\nabla ^{T}s|^{2}$ and
$|\nabla ^{T}|\nabla ^{T}T|^{2}|\leq 2|\nabla ^{T}T|\,|(\nabla ^{T})^{2}T|$,
they are bounded below by $2s^{-1}\bigl(|(\nabla ^{T})^{2}T|-|\nabla
^{T}T|\,|\nabla ^{T}\log s|\bigr)^{2}$. Hence
\begin{equation*}
L_{f}v\geq (4\lambda -C|T|)v+2\frac{|Q|^{2}}{s}v.
\end{equation*}
By \eqref{eq:Tsqrt}, \eqref{eq:T-first-bounded} and Young's inequality,
\begin{equation*}
C|T|v\leq C\frac{|\nabla ^{T}T|^{2}}{\sqrt{s}}\leq C\sqrt{v}\leq 2\lambda
v+C.
\end{equation*}
Thus
\begin{equation*}
L_{f}v\geq 2\lambda v-C.
\end{equation*}
Also $v\leq Cf$ by \eqref{eq:sf} and \eqref{eq:T-first-bounded}, so $v\in
L^{2}(d\nu _{f})$. Lemma~\ref{lem:L2} proves
\begin{equation}
|\nabla ^{T}T|^{2}\leq Cs.  \label{eq:gradTsqrt}
\end{equation}
Combining this with \eqref{eq:4-algebra} and $sf\geq c$ improves the
algebraic estimate to
\begin{equation}
|T|\leq C(|Q|+s).  \label{eq:T-linear-Q}
\end{equation}
Indeed, $|\nabla ^{T}Q|/\sqrt{f}\leq C\sqrt{s}/\sqrt{f}\leq Cs$ outside
$\Omega (t_{0})$. On the compact set $\Omega (t_{0})$ the tensor $T$ is
bounded and $s$ has a positive minimum, so, after enlarging $C$,
\eqref{eq:T-linear-Q} holds on all of $\mathcal{X}$.

Define
\begin{equation}
\widehat{L}=L_{f}+2\left\langle \nabla ^{T}\log s,\nabla ^{T}\,\cdot
\,\right\rangle =L_{f-2\log s},\qquad d\widehat{\nu }=s^{2}e^{-f}d\mu .
\label{eq:hatL}
\end{equation}
It is symmetric with respect to $d\widehat{\nu }$ on compactly supported
basic functions on $M$ (on compactly supported functions on $Z$ in the
proper $\mathbb{R}$ case); this is \eqref{eq:ibp} with $f$ replaced by
$f-2\log s$, which is smooth because $s>0$. Put
\begin{equation*}
V=Q/s\quad \text{and}\quad U=|V|^{2}=|Q|^{2}/s^{2}.
\end{equation*}
The tensor quotient identity is
\begin{equation*}
\widehat{L}V=-2T(V)+2\frac{|Q|^{2}}{s}V,\qquad T(V)_{il}=T_{iplq}V_{pq}.
\end{equation*}
Consequently
\begin{align}
\widehat{L}U& =2|\nabla ^{T}V|^{2}-4T_{ipjq}V_{ij}V_{pq}+4sU^{2}  \notag \\
& \geq s\bigl(4U^{2}-CU^{3/2}-CU\bigr)\geq s(3U^{2}-C_{1}U),
\label{eq:final-ratio}
\end{align}
where \eqref{eq:T-linear-Q} was used in the second line, in the form
$|T_{ipjq}V_{ij}V_{pq}|\leq |T|\,U\leq Cs(U^{1/2}+1)U$. This is the key
ratio inequality; it holds on all of $\mathcal{X}$.

The remaining maximum principle does \emph{not} require a prior estimate of
$|\nabla ^{T}\log s|$. Indeed, by \eqref{eq:Qsqrt},
\begin{equation}
\int U^{2}\,d\widehat{\nu }=\int \left( \frac{|Q|^{2}}{s}\right)
^{2}e^{-f}d\mu <\infty .  \label{eq:U-L2}
\end{equation}
Choose $U_{0}>C_{1}$ and set $w=(U-U_{0})_{+}$. The function $U$ is smooth
and $\theta _{r}^{2}w$ is Lipschitz with compact support, so the symmetry of
$\widehat{L}$ gives, exactly as in the proof of Lemma~\ref{lem:L2},
\begin{equation*}
\int \theta _{r}^{2}w\,\widehat{L}U\,d\widehat{\nu }=-\int \theta
_{r}^{2}|\nabla ^{T}w|^{2}\,d\widehat{\nu }-2\int \theta _{r}w\left\langle
\nabla ^{T}\theta _{r},\nabla ^{T}w\right\rangle \,d\widehat{\nu }\leq \int
|\nabla ^{T}\theta _{r}|^{2}w^{2}\,d\widehat{\nu }.
\end{equation*}
Testing \eqref{eq:final-ratio} with $\theta _{r}^{2}w$ and integrating
against $d\widehat{\nu }$ therefore gives
\begin{equation*}
\int \theta _{r}^{2}w\,s(3U^{2}-C_{1}U)\,d\widehat{\nu }\leq \int |\nabla
^{T}\theta _{r}|^{2}w^{2}\,d\widehat{\nu }\leq \frac{C}{r}\int
U^{2}\,d\widehat{\nu }\longrightarrow 0.
\end{equation*}
The left integrand is nonnegative and is strictly positive wherever $w>0$.
Fatou's lemma therefore gives $w=0$. This proves
\begin{equation*}
|Q|\leq Cs,
\end{equation*}
and \eqref{eq:T-linear-Q} yields
\begin{equation*}
|T|\leq Cs.
\end{equation*}
Finally \eqref{eq:pre-Rmg} gives $|\mathrm{Rm}_{g}|\leq C(1+|T|)\leq
C(1+s)$.

This completes the proof of \eqref{eq:final-transverse} and of
Theorem~\ref{thm:curvature}.

\subsection{The Derivative Estimate and Comparison with Munteanu--Wang}

The weighted integral argument used for \eqref{eq:final-ratio} also gives
the Sasakian analogue of \cite[Theorem~2.6]{MW1}. It is not needed for
Theorem~\ref{thm:curvature}, but it sharpens the asymptotic estimates of
Section~\ref{sec:end}: it gives the rank-one estimate of
Theorem~\ref{thm:general} (see Remark~\ref{rem:rank-improved}).

\begin{proposition}
\label{prop:gradT}Under the hypotheses of Theorem~\ref{thm:curvature}, in
the nonflat transverse case there is a constant $C>0$ such that
\begin{equation}
|\nabla ^{T}\mathrm{Rm}^{T}|\leq CR^{T}\quad \text{and}\quad |\nabla
^{T}\log R^{T}|\leq C\qquad \text{on }M.  \label{eq:gradT-linear}
\end{equation}
\end{proposition}

\begin{proof}
Put $A=|\nabla ^{T}T|^{2}$ and $w=A/s^{2}$. By \eqref{eq:trans-scalar},
$L_{f}(s^{-2})=-4\lambda s^{-2}+4|Q|^{2}s^{-3}+6s^{-4}|\nabla ^{T}s|^{2}$.
Together with \eqref{eq:deriv-norm} and the product rule this gives
\begin{equation*}
\widehat{L}w\geq \frac{2|(\nabla ^{T})^{2}T|^{2}}{s^{2}}+(2\lambda
-C|T|)w+4\frac{|Q|^{2}}{s}w+\frac{2|\nabla ^{T}s|^{2}A}{s^{4}}-\frac{2\langle
\nabla ^{T}s,\nabla ^{T}A\rangle }{s^{3}}.
\end{equation*}
Since $|\nabla ^{T}A|\leq 2\sqrt{A}\,|(\nabla ^{T})^{2}T|$, the last term is
at least $-2s^{-2}|(\nabla ^{T})^{2}T|^{2}-2s^{-4}|\nabla ^{T}s|^{2}A$.
Hence, by $|T|\leq Cs$,
\begin{equation*}
\widehat{L}w\geq (2\lambda -C_{2}s)w.
\end{equation*}
Let $s_{1}=\lambda /C_{2}$ and $W_{0}=s_{1}^{-2}\sup_{\mathcal{X}}|\nabla
^{T}T|^{2}$, which is finite by \eqref{eq:T-first-bounded}. Then $w\leq
W_{0}$ on $\{s\geq s_{1}\}$, so that $\widehat{w}=(w-W_{0})_{+}$ is
supported in $\{s<s_{1}\}$, where $\widehat{L}w\geq \lambda w\geq \lambda
\widehat{w}$. Testing with $\theta _{r}^{2}\widehat{w}$ and integrating against
$d\widehat{\nu }$, exactly as in the last step of the proof of
Theorem~\ref{thm:curvature}, we obtain
\begin{equation*}
\lambda \int \theta _{r}^{2}\widehat{w}^{\,2}\,d\widehat{\nu }\leq \int
|\nabla ^{T}\theta _{r}|^{2}\widehat{w}^{\,2}\,d\widehat{\nu }\leq
\frac{C}{r}\int \Bigl(\frac{|\nabla ^{T}T|^{2}}{s}\Bigr)^{2}e^{-f}d\mu
\longrightarrow
0,
\end{equation*}
by \eqref{eq:gradTsqrt} and the finiteness of $\nu _{f}$. Thus
$\widehat{w}=0$, that is, $|\nabla ^{T}T|^{2}\leq W_{0}s^{2}$. Finally $s$ is
a
contraction of $T$, so $|\nabla ^{T}s|\leq C|\nabla ^{T}T|\leq Cs$.
\end{proof}

\begin{remark}
\label{rem:MW-comparison}The proof of Theorem~\ref{thm:curvature} follows
\cite[Sections~1 and~2]{MW1}, with the following differences.

\begin{enumerate}
\item The estimate \eqref{eq:4-algebra} is obtained by restricting the
algebraic curvature tensor $T$ to the orthogonal complement of $\nabla
^{T}f$. The Gauss equation of the level sets of $f$ used in
\cite[Proposition~1.1]{MW1} is not needed, and the additional term
$(|\mathrm{Ric}|^{2}+1)/f$ of that proposition does not appear.

\item The Ricci bound is proved in two steps: first with the fixed exponent
$a=\frac{1}{2}$ in \eqref{eq:quotient-a}, and then, once $T$ is bounded, with
$a=1$ in \eqref{eq:Y-over-s}. In \cite[Proposition~1.3]{MW1} the exponent
$a=1-C/r$ depends on the cutoff.

\item The estimates \eqref{eq:gradTsqrt}, $|Q|\leq Cs$ and
\eqref{eq:gradT-linear} are derived from integral identities for the
weighted measures $e^{-f}d\mu $ and $s^{2}e^{-f}d\mu $, in place of the
maximum principle with the linear cutoff $(r-f)/r$ of
\cite[Proposition~2.2, Theorems~2.5 and~2.6]{MW1}. In particular the logarithmic gradient estimate
$|\nabla \log R|^{2}\leq C\log (f+2)$ of \cite[Lemma~2.3]{MW1} for the
scalar curvature $R$ of the Riemannian shrinker, used there to control the
drift Laplacian $\Delta _{f-2\log R}$ of the cutoff, is not needed here; the bound $|\nabla
^{T}\log s|\leq C$ is a consequence, not an ingredient.

\item All maximum principles and integrations by parts are carried out on
the exhaustion space $\mathcal{X}$ of Lemma~\ref{lem:reeb}, on which $f$ is
proper. No leaf space is used unless the Reeb flow is a proper
$\mathbb{R}$-action.
\end{enumerate}
\end{remark}

\section{Convergence of Controlled Sasaki--Ricci Flows}

\label{sec:flows}

In this section we introduce the controlled Sasaki--Ricci flows, in
particular the ancient flow associated with a soliton, and the developed
local limits used in Section~\ref{sec:splitting}. Flows are normalized by
\eqref{eq:marked-flow} and, for solitons, $\lambda =\frac{1}{2}$ (see
\eqref{eq:intro-normal}).

\subsection{Controlled Sasaki--Ricci Flows}

We first define controlled Sasaki--Ricci flows.

\begin{definition}
\label{def:controlled}Let $(M,g(t))=(M,\xi ,\eta _{t},\Phi _{t},g(t))$, $t\in
(-\infty ,0)$, be a Sasaki--Ricci flow of complete Sasakian
$(2n+1)$-manifolds. It is called \emph{transversely $C$-controlled} if
\begin{equation*}
\bigl|\mathrm{Rm}^{T}(g(t))\bigr|_{g^{T}(t)}\leq \frac{C}{|t|},
\end{equation*}
and, in setting (S), it is called \emph{transversely $(C,\kappa)$-controlled}
if it is, in addition, \emph{transversely $\kappa$-noncollapsed} with respect
to the specified quotient measure: for every $(x,t,r)\in M\times
(-\infty ,0)\times (0,\infty )$ such that
\begin{equation*}
\bigl|\mathrm{Rm}^{T}\bigr|\leq \frac{1}{r^{2}}
\end{equation*}
on the transverse parabolic tube $P^{T}(x,t,r)=B_{g(t)}^{T}(x,r)\times
(t-r^{2},t]$, one has
\begin{equation*}
\frac{\operatorname{Vol}^{T}_{t}\bigl(B_{g(t)}^{T}(x,r)\bigr)}{r^{2n}}\geq \kappa .
\end{equation*}
Here the \emph{transverse ball} (or transverse tube) $B_{g(t)}^{T}(x,r)$ is
the set
\begin{equation*}
\bigl\{y\in M:\ d_{g(t)}^{T}(y,x)<r\bigr\},
\end{equation*}
where the transverse distance $d_{g(t)}^{T}(y,x)$ is the infimum of the
$g(t)$-distances between points of the Reeb orbits through $x$ and $y$. The
transverse distance is a pseudometric on $M$: it vanishes between points
on the same Reeb orbit. In the quasi-regular case it induces the orbifold
distance on the leaf space; in setting (S), on the chosen total space,
it is the pullback of the distance of $(Z,h_Z(t))$. In particular,
\[
B_{g(t)}^{T}(x,r)=\pi^{-1}\bigl(B_{h_Z(t)}(\pi(x),r)\bigr),\qquad
\operatorname{Vol}^{T}_{t}\bigl(B_{g(t)}^{T}(x,r)\bigr)
=\operatorname{Vol}_{h_Z(t)}B_{h_Z(t)}(\pi(x),r).
\]
The spatial set in $P^{T}(x,t,r)$ is fixed using the final-time metric;
the curvature at time $t'\in(t-r^2,t]$ is measured using
$g^T(t')$.
This quotient volume is the measure in \eqref{eq:analyticspace}.
\end{definition}

The local Bernstein--Bando--Shi estimates \cite{Shi} hold for the
transverse curvature along the Sasaki--Ricci flow (see \cite{Co1} and
\cite{SWZ}); by parabolic rescaling they give the following fact.

\begin{lemma}
If $(M^{2n+1},g(t))$ is a transversely $C$-controlled Sasaki--Ricci flow, then
for all nonnegative integers $k,l$ there is a constant $C_{k,l}>0$, depending
only on $n$ and $C$, such that
\begin{equation*}
\Bigl|\frac{\partial ^{k}}{\partial t^{k}}(\nabla
^{T})^{l}\mathrm{Rm}^{T}(g(t))\Bigr|\leq \frac{C_{k,l}}{|t|^{1+k+l/2}}.
\end{equation*}
\end{lemma}

We omit the routine scaling argument (compare \cite[Lemma~2.2]{Nab}). The
lemma is used in the estimates
of the reduced transverse length in Appendix~\ref{app:reduced}.

\begin{proposition}
\label{lem:selfsim}Let $(M,\xi ,\eta ,\Phi ,g,f)$ be a complete gradient
shrinking Sasaki--Ricci soliton with $\lambda =\frac{1}{2}$ (see
\eqref{eq:intro-normal} and Remark~\ref{rem:Dhom}) and bounded transverse
curvature. Write $h=g^T$ for the initial transverse metric. The gradient
$\nabla^T f$ in the following ordinary differential equation is taken
with respect to this initial metric. For $\tau>0$ let $\Upsilon_\tau$ solve
\begin{equation*}
\frac{d}{d\tau }\Upsilon _{\tau }=-\tau ^{-1}\nabla ^{T}f\circ \Upsilon
_{\tau },\qquad \Upsilon _{1}=\mathrm{Id}.
\end{equation*}
Then $h(\tau )=\tau \Upsilon _{\tau }^{\ast }h$ satisfies
\begin{equation*}
\partial _{\tau }h=2\,\mathrm{Ric}^{T}(h).
\end{equation*}
Equivalently, with $t=-\tau $ and $\Psi _{t}=\Upsilon _{-t}$, the family
$g^{T}(t)=(-t)\,\Psi _{t}^{\ast }g^{T}$ solves \eqref{eq:marked-flow} on
$(-\infty ,0)$ with $g^{T}(-1)=g^{T}$, and $|\mathrm{Rm}^{T}(g^{T}(t))|\leq
C/|t|$; in particular the flow is transversely $C$-controlled. The
diffeomorphisms $\Upsilon _{\tau }$ preserve the Reeb foliation and the
transverse holomorphic structure; they need not preserve the original metric
or contact form.
\end{proposition}

\begin{proof}
Boundedness of $\mathrm{Rm}^{T}$ gives boundedness of
$\operatorname{Hess}^{T}f=\lambda g^{T}-\mathrm{Ric}^{T}$, so the gradient has
at most linear
growth and the flow exists for all finite intervals of $\log \tau $. It
commutes with $\xi $. Because $g^{T}$ and $\mathrm{Ric}^{T}$ are
$J$-invariant, so is $\operatorname{Hess}^{T}f$; hence $\nabla ^{T1,0}f$ is a
transverse holomorphic vector field with Hamiltonian potential $f$, and
$\Upsilon _{\tau }$ preserves the transverse holomorphic structure.
Differentiating gives
\begin{equation*}
\partial _{\tau }(\tau \Upsilon _{\tau }^{\ast }h)=\Upsilon _{\tau }^{\ast
}\bigl(h-\mathcal{L}_{\nabla ^{T}f}h\bigr)=\Upsilon _{\tau }^{\ast
}\bigl(h-2\,\operatorname{Hess}^{T}f\bigr)=2\Upsilon _{\tau }^{\ast
}\mathrm{Ric}^{T}.
\end{equation*}
Constant scaling leaves the covariant Ricci tensor unchanged, so the right
side is $2\,\mathrm{Ric}^{T}(h(\tau ))$. The curvature bound follows by
scaling.
\end{proof}

\begin{proposition}
\label{prop:derivatives}Let $(M,\xi ,\eta ,\Phi ,g,f)$ be a complete
gradient shrinking Sasaki--Ricci soliton with bounded transverse curvature;
for instance, assume the hypotheses of Theorem~\ref{thm:curvature}. Then,
for every $k\geq 0$,
\begin{equation}
\sup_{M}|(\nabla ^{T})^{k}T|<\infty .  \label{eq:allderivatives}
\end{equation}
The associated ancient flow $g^{T}(t)$ of Proposition~\ref{lem:selfsim} satisfies
\begin{equation*}
|(\nabla ^{T})^{k}\mathrm{Rm}^{T}(g^{T}(t))|\leq C_{k}|t|^{-1-k/2}.
\end{equation*}
\end{proposition}

\begin{proof}
Commuting transverse derivatives with the curvature equation \eqref{eq:LT},
and using the soliton equation and the Bianchi identity, gives, for $k\geq 1$,
\begin{equation}
L_{f}\bigl((\nabla ^{T})^{k}T\bigr)=(k+2)\lambda (\nabla
^{T})^{k}T+\sum_{j=0}^{k}(\nabla ^{T})^{j}T\ast (\nabla ^{T})^{k-j}T.
\label{eq:deriv-k}
\end{equation}
Here the contractions depend on $k$, and derivatives of $\mathrm{Ric}^{T}$
are contractions of derivatives of $T$. The apparent curvature-$\nabla ^{T}f$
terms produced by commutation are replaced using $\mathrm{div}^{T}T=T(\nabla
^{T}f,\cdot ,\cdot ,\cdot )$ from \eqref{eq:divT}; the
$\operatorname{Hess}^{T}f$
term contributes the extra $k\lambda $. For $k=1$ this is \eqref{eq:deriv}.

Suppose derivatives below order $k$ are bounded. For $\Xi _{k}=|(\nabla
^{T})^{k}T|$, Kato's inequality and \eqref{eq:deriv-k} imply
\begin{equation*}
L_{f}\Xi _{k}\geq -C\Xi _{k}-C,\qquad L_{f}\Xi _{k-1}^{2}\geq 2\Xi
_{k}^{2}-C.
\end{equation*}
The second inequality uses the induction hypothesis. Thus, after changing
constants,
\begin{equation*}
L_{f}(\Xi _{k}+\Xi _{k-1}^{2})\geq c(\Xi _{k}+\Xi _{k-1}^{2})^{2}-C
\end{equation*}
outside a compact sublevel, and the annular maximum calculation
\eqref{eq:cutoff-max} proves the induction step. The flow statement follows
from Proposition~\ref{lem:selfsim} by scaling.
\end{proof}

\subsection{Local Compactness}

\begin{definition}
\label{def:developed}Let $(M,\xi ,\eta ,\Phi ,g)$ be a complete Sasakian
manifold with
\begin{equation*}
C_{0}=\sup_{M}|\mathrm{Rm}_{g}|<\infty \quad \text{and}\quad \sup_{M}|\nabla
^{k}\mathrm{Rm}_{g}|<\infty
\end{equation*}
for every $k\geq 1$, and put $\rho _{0}=\pi /(2\sqrt{C_{0}})$. Since $|\sec
_{g}|\leq C_{0}$, the exponential map $\exp _{p}$ is a local diffeomorphism
on the ball $B_{\rho _{0}}(0)\subset T_{p}M$ for every $p\in M$. Fix $0<\rho
\leq \rho _{0}$.

\begin{enumerate}
\item \emph{Developed Sasakian chart and transverse chart.} The developed
Sasakian chart of radius $\rho$ at $p$ is the $(2n+1)$-dimensional ball
$B_\rho(0)\subset T_pM$ with the pulled-back structure
$\exp_p^*(\xi,\eta,\Phi,g)$. It is Sasakian because $\exp_p$ is a local
isometry for the pulled-back metric. In a sufficiently small flow box
about $0$, quotienting by the pulled-back Reeb flow gives a
$2n$-dimensional K\"{a}hler metric on a transversal through $0$.
This local quotient, well defined up to restriction and holomorphic
isometry, is the \emph{transverse chart} at $p$.

\item \emph{Developed local transverse limit.} Let $p_{i}\in M$, and
identify $T_{p_{i}}M$ with $\mathbb{R}^{2n+1}$ by orthonormal frames whose
last vector is $\xi (p_{i})$. A \emph{developed local limit} at $(p_{i})$ is
a Sasakian structure $(\xi _{\infty },\eta _{\infty },\Phi _{\infty
},g_{\infty })$ on $B_{\rho /2}(0)\subset \mathbb{R}^{2n+1}$ such that,
along a subsequence, the charts of radius $\rho $ at $p_{i}$ converge to it
in $C^{\infty }(B_{\rho /2}(0))$. The K\"{a}hler structure that it induces
on a transversal of $\xi _{\infty }$ through $0$ is the corresponding
\emph{developed local transverse limit}.

\item \emph{Developed local limits are local objects.} The charts need not be
injective images of neighborhoods in $M$, the limits do not determine a
complete limiting quotient, and different subsequences may have different
limits.
\end{enumerate}
\end{definition}

\begin{proposition}
\label{prop:unwrap}Let $(M,\xi ,\eta ,\Phi ,g,f)$ be a complete gradient
shrinking Sasaki--Ricci soliton with bounded transverse curvature. Then
\begin{equation*}
\sup_{M}|\nabla ^{k}\mathrm{Rm}_{g}|<\infty
\end{equation*}
for every $k\geq 0$, and for every $0<\rho \leq \rho _{0},$ every sequence
of basepoints admits developed local limits in the sense of
Definition~\ref{def:developed}.
\end{proposition}

\begin{proof}
Proposition~\ref{prop:derivatives} bounds all transverse derivatives of the
transverse curvature. By \eqref{eq:pre-conversion} and \eqref{eqn-SR2}, the
curvature tensor $\mathrm{Rm}_{g}$ is the sum of $T$, extended by zero in
the Reeb direction, and of a tensor which is algebraic in $g$, $\eta $ and
$\Phi $. The covariant derivatives of the structure tensors are algebraic in
the same tensors: $\nabla \xi =\Phi $, $\nabla \eta =g(\Phi \,\cdot ,\cdot
)$ and $(\nabla _{X}\Phi )Y=\eta (Y)X-g(X,Y)\xi $. By \eqref{eq:pre-nablaT},
on basic horizontal tensors the Levi-Civita connection of $g$ differs from
$\nabla ^{T}$ by terms which are algebraic in $\Phi $, $\eta $ and $\xi $;
in the Reeb direction the derivative of a basic tensor is itself algebraic,
because its Lie derivative along $\xi $ vanishes. By induction on $k$,
$\nabla ^{k}\mathrm{Rm}_{g}$ is a universal polynomial in $(\nabla
^{T})^{j}T$, $j\leq k$, and in $g$, $\eta $, $\xi $, $\Phi $. Thus all full
curvature derivatives are bounded.

Let $0<\rho \leq \rho _{0}$. Since $\rho _{0}$ is smaller than the conjugate
radius, the exponential map at each basepoint is a local diffeomorphism on
its tangent ball of radius $\rho $, whether or not it is injective. Pull
back the metric and the Sasakian tensors to this tangent ball. Jacobi-field
estimates and their differentiated equations give uniform nondegeneracy and
derivative bounds on $B_{\rho /2}(0)$. After choosing frames and passing to
a subsequence, these tensors converge smoothly there. The limiting unit Reeb
field admits a flow, and the transverse tensors converge on its transversal.
This produces local immersions and local metrics. It does not produce an
embedded exhaustion of a complete limiting quotient.
\end{proof}

\section{Splitting Theorems at Infinity}

\label{sec:splitting}

In setting (S) we apply Perelman's reduced geometry on the complete
smooth quotient, as developed in Appendix~\ref{app:reduced}, to obtain
the corresponding global splitting results of Naber \cite{Nab}. In
substance these global results are Naber's theorems applied to the
complete K\"{a}hler--Ricci flow on the quotient $Z$; we state them in
Sasakian form because this is the form in which they are used in
\cite{PartII}.
For arbitrary Reeb type, Theorem~\ref{thm:splitting} instead uses normalized
potentials on developed local charts and requires no global leaf space
or global reduced-volume construction.

Note that $\nabla ^{T}f$ is a complete transverse vector field for a
gradient shrinking Sasaki--Ricci soliton with a basic potential $f$ and
$|\mathrm{Rm}^{T}|\leq C$, since its transverse Hessian is bounded.

\begin{definition}
\label{def:soliton}Recall that $(M,\xi ,\eta ,\Phi ,g,f)$ is a gradient
shrinking Sasaki--Ricci soliton if $f\in C_{B}^{\infty }(M)$ and
\begin{equation*}
\mathrm{Ric}^{T}+\operatorname{Hess}^{T}f=\lambda g^{T}
\end{equation*}
with $\lambda >0$. Taking the trace and the divergence of this
equation gives
\begin{equation*}
R^{T}+\Delta _{B}f=2n\lambda \quad \text{and}\quad R^{T}+|\nabla
^{T}f|^{2}-2\lambda f=c
\end{equation*}
for a constant $c$, so that
\begin{equation*}
\Delta _{B}f-|\nabla ^{T}f|^{2}+2\lambda f=2n\lambda -c.
\end{equation*}
The soliton is \emph{normalized} if $c=0$, as in \eqref{eq:pre-normal}; then
\begin{equation*}
\Delta _{B}f-|\nabla ^{T}f|^{2}+2\lambda f=2n\lambda .
\end{equation*}
In the quotient-flow and local-splitting arguments below we take
$\lambda=\frac12$; Subsection~\ref{sec:end} states its estimates with
$\lambda>0$ explicit. Appendix~\ref{app:reduced} uses the flow normalization
$\partial_t h=-2\mathrm{Ric}_h$. In setting (S), write $f_Z$ for the
function on the chosen quotient characterized by $f=\pi^*f_Z$.
The \emph{$f$-volume} is
\begin{equation*}
\operatorname{Vol}_{f}^{T}(M):=\int_Z e^{-f_Z}\,d\mathrm{vol}_{h},
\qquad h=h_Z(-1).
\end{equation*}
Thus the notation $M$ records the chosen Sasakian lift, while the
integration takes place on its smooth quotient at the normalized time
$-1$. No factor $(4\pi)^{-n}$ is included in this definition.
\end{definition}

The $f$-volume of a normalized shrinking soliton does not depend on the
choice of the normalized potential and is invariant under quotient isometries
(compare \cite[Lemma~2.1]{Nab} and \cite{CN}).

\subsection{Asymptotic Solitons and Transverse Type I Limits}

We first work in setting (S); all reduced lengths and volumes below
are the quotient objects of Appendix~\ref{app:reduced}.

\begin{theorem}
\label{thm:asymp}In setting (S), let $(Z,h(t))$, $t<0$, be a
complete K\"{a}hler--Ricci flow with
\begin{equation*}
|\mathrm{Rm}(h(t))|\leq C/|t|
\end{equation*}
which is $\kappa $-noncollapsed. Fix $x\in Z$ and let $\tau
_{i}^{-}\downarrow 0$, $\tau _{i}^{+}\uparrow \infty $. Subsequences of
\begin{equation*}
(Z,(\tau _{i}^{\pm })^{-1}h(\tau _{i}^{\pm }t),x)
\end{equation*}
converge smoothly on compact subsets of space-time, $t<0$, to the flows of complete noncollapsed gradient shrinking
K\"{a}hler--Ricci solitons, normalized at $t=-1$. If
the two limiting metrics are holomorphically isometric, the original
quotient flow is also a shrinking soliton, isometric to them. The conclusion
concerns the smooth quotients; by Proposition~\ref{prop:contactization},
the limiting solitons admit local Sasakian contactizations.
\end{theorem}

\begin{proof}
This is Naber's asymptotic-soliton theorem \cite[Theorem~2.1]{Nab} for the
complete quotient flow on $Z$. Quotient noncollapsing, the curvature bounds
and Shi's
estimates give pointed smooth limits $S^{\pm }$, which are K\"{a}hler since
the complex structures are parallel. By Proposition~\ref{prop:singular} and
scale invariance, the singular reduced lengths converge to functions $\ell
^{T,\pm }$ whose reduced volumes are the constants $c_{0}$ and $c_{\infty }$
of \eqref{eq:endpointvolumes}. Hence the defect $\delta $ in \eqref{eq:magic}
and \eqref{eq:magic2}
vanishes, i.e.\ both hold with equality, and the equality case of
Perelman's monotonicity \cite[Section~7]{Per1} (see also
the proof of \cite[Theorem~2.1]{Nab}), applied on the smooth quotient limit,
gives
\begin{equation*}
\mathrm{Ric}^{T}+\operatorname{Hess}^{T}\ell ^{T,\pm
}=\frac{1}{2\bar{\tau}}g^{T},
\end{equation*}
with the normalization of Definition~\ref{def:soliton}. If $S^{+}$ and $S^{-}
$ are isometric, then $c_{0}=c_{\infty }$ by the independence of the $f
$-volume of the normalized potential, so $\widetilde{V}^{T}$ is constant and
the flow is itself a shrinking soliton.
\end{proof}

\begin{remark}
\label{rem:careful} As in \cite{Nab}, if $(M,g(t))$ is already a transverse
shrinking soliton it need not be transversely isometric to $S^{-}$, and the
singular reduced transverse length from $(x,0)$ need not be a soliton
potential. For a shrinking quotient, the backward asymptotic soliton
$S^{+}$ is isometric to the original quotient \cite[Remark~2.7 and
Theorem~3.2]{Nab}. The proof of Theorem~\ref{thm:kappa} below invokes
Naber's noncollapsing theorem directly.
\end{remark}

\begin{definition}
\label{def:typeI} Let a Sasaki--Ricci flow $(M,g(t))$ have a finite
singular time $0<T_0<\infty$ at the endpoint of its maximal interval
$[0,T_0)$. It has a \emph{transverse Type I singularity} at $T_0$ if
there is $C>0$ with
\begin{equation*}
\big|{\mathrm{Rm}^{T}}[g(t)]\big|\leq \frac{C}{|T_{0}-t|}.
\end{equation*}
\end{definition}

\begin{theorem}
\label{thm:typeI}In setting (S), let $(Z,h(t))$, $0\leq t<T_{0}<\infty $, be
a complete K\"{a}hler--Ricci flow with
\begin{equation*}
|\mathrm{Rm}(h(t))|\leq C/(T_{0}-t).
\end{equation*}
Assume explicitly that the quotient flow is uniformly $\kappa $-noncollapsed
at the curvature-admissible scales used in the dilations. Fix $x\in Z$, let
$t_{i}\uparrow T_{0}$, and set $\tau _{i}=T_{0}-t_{i}$. A subsequence of
\begin{equation*}
(Z,h_i(t),x),\qquad h_i(t)=\tau_i^{-1}h(T_0+\tau_i t),\qquad
-T_0/\tau_i\leq t<0,
\end{equation*}
converges smoothly on every compact space-time subset of the limiting
interval $(-\infty,0)$ to a complete normalized gradient shrinking
K\"{a}hler--Ricci soliton with bounded curvature. The conclusion allows
a flat limit at a basepoint that does not detect the singularity.
\end{theorem}

\begin{proof}
Apply the smooth Type I tangent-flow theorem of Naber
\cite[Theorem~1.4]{Nab} (see also Enders--M\"{u}ller--Topping \cite{EMT}) on
$Z$,
using the stated noncollapsing hypothesis for compactness. The parallel
complex structures pass to the limit. A local contactization is supplied by
Proposition~\ref{prop:contactization}. We have retained noncollapsing as an
explicit sufficient hypothesis instead of deducing a uniform reduced-volume
lower bound from pointwise bounds for the reduced length: such a deduction
also needs volume control. This extra hypothesis is automatic for the
bounded-curvature gradient quotient shrinkers used in
\cite[Theorem~\IIref{T51}{4.3}]{PartII}, by Theorem~\ref{thm:kappa} applied on
the quotient.
\end{proof}

\subsection{Noncollapsing and the Gradient Theorem}

\label{sec:noncollapsing}

\begin{theorem}
\label{thm:perelmanNC} Let $(M,g(t))$, $t\in \lbrack 0,T]$, be a
Sasaki--Ricci flow whose transverse metrics $h(t)$, normalized by
\eqref{eq:marked-flow}, are as in setting (S), and let
$(x,T)\in M\times \{T\}$,
$0<\bar{\tau}\leq T$ and $V>0$. There is $\kappa =\kappa (n,V)>0$ such that if
$0<r\leq \sqrt{\bar{\tau}}$, $|\mathrm{Rm}^{T}|\leq r^{-2}$ on the transverse
parabolic tube $P^{T}(x,T,r)$, and $V_{(x,T)}^{T}(\bar{\tau})\geq V$, then
\begin{equation*}
\operatorname{Vol}^{T}\!\big(B_{r}^{T}(x)\big)\geq \kappa r^{2n},
\end{equation*}
where $\operatorname{Vol}^{T}$ is the quotient measure \eqref{eq:analyticspace}.
\end{theorem}

\begin{proof}
Let $\bar x=\pi(x)$. By the definitions of Appendix~\ref{app:reduced},
$V_{(x,T)}^T(\bar\tau)$ is exactly the normalized Perelman reduced volume
of $(Z,h(t))$ based at $(\bar x,T)$, and the image of the transverse
parabolic tube is
$B_{h(T)}(\bar x,r)\times(T-r^2,T]$.
The curvature condition is precisely $|\mathrm{Rm}_h|\le r^{-2}$ there.
Apply Perelman's quotient no-local-collapsing theorem in real dimension
$m=2n$, in the explicit formulation \cite[Theorem~3.1]{Nab}
(see also \cite[Chapter~8]{MT}). Its hypotheses are completeness of
this smooth flow, $r\le\sqrt{\bar\tau}$, the displayed parabolic
curvature bound, and the lower reduced-volume bound $V$.
It gives $\operatorname{Vol}_{h(T)}B_{h(T)}(\bar x,r)\ge\kappa(2n,V)r^{2n}$.
The quotient measure \eqref{eq:analyticspace} identifies this volume
with the stated $\operatorname{Vol}^T$. If setting (S) uses a finite cover,
all these objects are formed on its specified quotient with the chosen
lift of $x$, as required by the definition of setting (S).
\end{proof}

\begin{theorem}
\label{thm:kappa} Let $(M,\xi ,\eta ,\Phi ,g,f)$ be a normalized complete
gradient shrinking Sasaki--Ricci soliton with bounded transverse curvature in setting
(S). Then there is
\begin{equation*}
\kappa =\kappa \big(n,\operatorname{Vol}_{f}^{T}(M)\big)>0
\end{equation*}
such that the associated Sasaki--Ricci flow of Proposition~\ref{lem:selfsim}
is transversely $\kappa $-noncollapsed.
\end{theorem}

\begin{proof}
The basic potential descends to $f_Z$ on the complete smooth quotient,
where
\[
 \mathrm{Ric}_{h}+\operatorname{Hess}_{h}f_Z=\tfrac12h,
 \qquad R_h+|\nabla f_Z|^2=f_Z.
\]
Thus $(Z,h,f_Z)$ is a normalized complete Riemannian gradient shrinker
of real dimension $2n$ with bounded full curvature. Its $f$-volume is
positive and finite: Lemma~\ref{lem:growth}, after the fixed normalization
change, gives quadratic growth of $f_Z$, whereas a lower Ricci bound
from bounded curvature gives at most exponential volume growth.
The associated quotient Ricci flow is the one in
Proposition~\ref{lem:selfsim}. Naber's noncollapsing theorem
\cite[Theorem~1.3]{Nab} therefore applies directly, with
\[
 \kappa=\kappa\!\left(2n,\int_Ze^{-f_Z}d\mathrm{vol}_h\right)>0.
\]
This theorem concerns precisely the parabolic-ball condition in
Definition~\ref{def:controlled}; it does not assume noncollapsing in
advance. By \eqref{eq:analyticspace}, its volume and our quotient
$f$-volume agree. Pulling back the ball notation gives the assertion.
In particular this proof does not use Theorem~\ref{thm:asymp}, whose
hypotheses already include noncollapsing.
\end{proof}

\begin{corollary}
\label{cor:controlled}Let $(M,\xi ,\eta ,\Phi ,g,f)$ be a normalized complete
gradient shrinking Sasaki--Ricci soliton in setting (S) with
\begin{equation*}
|\mathrm{Rm}^{T}|\leq C.
\end{equation*}
Its associated flow $g^{T}(t)=(-t)\Psi _{t}^{\ast }g^{T}$ of
Proposition~\ref{lem:selfsim} satisfies
\begin{equation*}
|\mathrm{Rm}^{T}(g^{T}(t))|\leq C/|t|
\end{equation*}
for $t<0$ and is transversely $\kappa $-noncollapsed with $\kappa =\kappa
(n,\operatorname{Vol}_{f}^{T}(M))$. That is, it is transversely $(C,\kappa
)$-controlled.
In particular Theorem~\ref{thm:asymp} applies to it.
\end{corollary}

\begin{proof}
Curvature norms scale by $(-t)^{-1}$ under $g^{T}(t)=(-t)\Psi _{t}^{\ast
}g^{T}$, and the noncollapsing is Theorem~\ref{thm:kappa}.
\end{proof}

\begin{theorem}
\label{thm:gradient}Suppose setting (S) holds. Let $M$ denote the
original total space and $q:\widetilde M\to M$ the chosen finite
covering. Assume that a complete real transverse soliton vector field $X$
on $M$ lifts to a field on $\widetilde M$ that projects under
$\pi:\widetilde M\to Z$ to a vector field on the chosen smooth quotient.
Assume bounded transverse curvature and
\begin{equation*}
\mathrm{Ric}^{T}+\tfrac{1}{2}\mathcal{L}_{X}g^{T}=\lambda g^{T},\qquad
\lambda >0.
\end{equation*}
Then there is a global real basic function $f$ on the original total space
with
\begin{equation*}
\mathrm{Ric}^{T}+\operatorname{Hess}^{T}f=\lambda g^{T}.
\end{equation*}
Its transverse $(1,0)$ gradient is holomorphic.
\end{theorem}

\begin{proof}
The lift of $X$ to the finite covering is complete, and its projection
to $Z$ is complete because its integral curves are the projections
of the lifted integral curves. Apply Naber's gradient-structure theorem
\cite[Theorem~1.2]{Nab} to the complete bounded-curvature smooth
quotient $(Z,h_Z)$. Let $f_Z$ be
the resulting potential and put $\widetilde f=\pi^*f_Z$ on the chosen
total space. Write $d$ for the degree of $q$ and define on the original
total space
\[
f(p)=\frac1d\sum_{y\in q^{-1}(p)}\widetilde f(y).
\]
On an evenly covered neighborhood this is the average over the smooth
local inverse branches of $q$. It is therefore smooth, and the local
formulas agree on overlaps because the branches are merely permuted.
Each branch pulls the lifted transverse metric and Ricci tensor back
to the original ones. Every summand satisfies that equation. Since the equation is
affine in the potential and the weights $1/d$ sum to one, so does their
average $f$.
Each summand is basic, hence $\xi f=0$. This argument works for any
finite covering and does not require a deck group acting transitively
on its fibers. The Hessian is $J$-invariant, proving transverse
holomorphicity. Normalize the additive constant after averaging.
Finally, $X-\nabla^T f$ is transversely Killing by subtraction of the
two real soliton equations.
\end{proof}

\subsection{Naber-Type Splitting of the Leaf Space}

The following statements are the Sasakian formulation, in setting (S), of
Naber's results \cite[Propositions~4.1 and~4.2 and Corollary~4.1]{Nab} on the
leaf space. Together with
Theorem~\ref{thm:splitting}, they describe the global splitting on the smooth
quotient
used in the proof of \cite[Theorem~\IIref{T51}{4.3}]{PartII}. For irregular
foliations they are replaced by Corollary~\ref{prop:conditional}.

All pointed convergence in the rest of this subsection takes place on the
complete smooth quotient $Z$ in setting (S), or on the quotient of the
specified finite cover. Thus no degenerate tensor on a Sasakian total space
is treated as a Riemannian metric.

\begin{theorem}
\label{P54}Let $(M^{5},\xi ,\eta ,\Phi ,g,f)$ be a normalized complete gradient shrinking
Sasaki--Ricci soliton with $\lambda=\tfrac12$ in setting (S). Suppose that its smooth quotient
$(Z,h,J,f_{Z})$ is noncompact and $|\mathrm{Rm}(h)|\leq C$, and let
$h(t)=(-t)\Psi _{t}^{\ast }h$, $t<0$, be the quotient flow induced by the
ancient flow of Proposition~\ref{lem:selfsim}, so that $h(-1)=h$. For any $q_{i}\in
Z$ with $f_{Z}(q_{i})\rightarrow \infty $, a subsequence of $(Z,h(t),q_{i})$
converges smoothly on compact space-time sets, $t<0$, to a complete
noncollapsed ancient K\"{a}hler--Ricci flow
$(Z_{\infty },h_{\infty }(t),q_{\infty })$ with $|\mathrm{Rm}(h_{\infty }(t))|\leq C/|t|$, whose time $-1$ slice is the pointed limit
of $(Z,h,q_{i})$. Its universal cover splits as
\begin{equation*}
(\widehat{Z}_{\infty },\widehat{h}_{\infty
}(t))=(\mathbb{C},h_{\mathrm{E}})\times (N,k(t)),
\end{equation*}
where $N$ is complete and $k(t)$ is a noncollapsed ancient Ricci flow on a surface with
$|\mathrm{Rm}(k(t))|\leq C/|t|$. In particular, the limit
$(Z_{\infty },h_{\infty }(t))$ itself
splits locally in this fashion.
\end{theorem}

\begin{proof}
Corollary~\ref{cor:controlled} applies on $Z$. On compact time intervals in
$(-\infty ,0)$, its curvature bound and noncollapsing give an
injectivity-radius bound at the basepoints, and Shi's estimates and
Hamilton's compactness theorem give the pointed flow limit. Escape is
uniform on such intervals: from
\begin{equation*}
|\nabla f_{Z}|^{2}=f_{Z}-R
\end{equation*}
with $0\leq R\leq C_{0}$ one gets
\begin{equation*}
f_{Z}(F_{s}(q))-C_{0}\geq e^{s}(f_{Z}(q)-C_{0})
\end{equation*}
for $s\geq 0$ along the flow $F_{s}$ of $\nabla f_{Z}$, with a comparable
bound for bounded negative $s$, and $\Psi _{t}=F_{-\log (-t)}$. Hence the
normalized potentials give a parallel unit field in the limit, as in
Theorem~\ref{thm:splitting}, and the parallel complex structure gives a second
one.
De Rham splitting on the simply connected cover yields the flat complex
factor, which persists in time by the product evolution and uniqueness
\cite{CZ}.

We make the transfer of noncollapsing to the surface factor explicit. Let
$\kappa >0$ be a noncollapsing constant for the limit flow. The lifted flow
on its universal cover is also $\kappa $-noncollapsed: at any fixed time the
covering map sends each metric ball onto the corresponding ball downstairs,
by lifting minimizing geodesics. Consequently a curvature bound on the
parabolic neighborhood with that fixed final-time ball descends to the
corresponding neighborhood downstairs, while the volume of the upstairs ball
is at least that of the downstairs ball.

Now fix $q\in N$, $t_{0}<0$ and $r>0$ such that
\begin{equation*}
|\mathrm{Rm}(k(s))|_{k(s)}\leq r^{-2}\quad \text{on }B_{k(t_{0})}(q,r)\times
(t_{0}-r^{2},t_{0}].
\end{equation*}
The flat factor contributes no curvature, so the product flow satisfies the
same bound on its parabolic neighborhood based on $B_{\widehat{h}_{\infty
}(t_{0})}((0,q),r)$. At the final time $t_{0}$,
\begin{equation*}
B_{\widehat{h}_{\infty }(t_{0})}((0,q),r)\subset
B_{h_{\mathrm{E}}}(0,r)\times B_{k(t_{0})}(q,r).
\end{equation*}
Product noncollapsing and the product volume form therefore give
\begin{equation*}
\kappa r^{4}\leq \operatorname{Vol}_{\widehat{h}_{\infty
}(t_{0})}B_{\widehat{h}_{\infty }(t_{0})}((0,q),r)\leq \pi
r^{2}\operatorname{Vol}_{k(t_{0})}B_{k(t_{0})}(q,r).
\end{equation*}
Thus
\begin{equation*}
\operatorname{Vol}_{k(t_{0})}B_{k(t_{0})}(q,r)\geq (\kappa /\pi )r^{2}
\end{equation*}
at every such parabolic scale, proving that $(N,k(t))$ is $\kappa /\pi
$-noncollapsed in the parabolic sense of Definition~\ref{def:controlled}.
\end{proof}

\begin{theorem}
\label{thm:naber-split}Assume the hypotheses of Theorem~\ref{P54}; in
particular, $M$ is complete and five-dimensional and its chosen smooth
quotient $(Z,h,J,f_Z)$ is noncompact. Every sequence $q_i\in Z$ with
$f_Z(q_i)\to\infty$ has a subsequence, still denoted $q_i$, and two
auxiliary sequences $q_j^{\pm}\in Z$ with $f_Z(q_j^{\pm})\to\infty$ such that
\begin{equation*}
(Z,h,q_{i})\longrightarrow (Z_{\infty },h_{\infty },q_{\infty }),\qquad
(Z,h,q_{j}^{\pm })\longrightarrow (Z_{\infty }^{\pm },h_{\infty }^{\pm
},q_{\infty }^{\pm }).
\end{equation*}
The first limit is the time $-1$ slice of the flow in Theorem~\ref{P54}. The
other two are its normalized asymptotic shrinking
K\"{a}hler--Ricci solitons. Each
limit splits locally off a flat complex line, and the universal covers split
globally with complete surface factors. If the two pointed asymptotic
solitons are isometric, the first limit is isometric to them and is a
shrinking soliton. If $q_{i}$ belong to one gradient trajectory, the
auxiliary points may be chosen on that same trajectory.
\end{theorem}

\begin{proof}
Let $h(t)=(-t)\Psi _{t}^{\ast }h$ be the quotient flow; Theorem~\ref{P54}
gives smooth convergence to $h_{\infty }(t)$. Apply Naber's
asymptotic-soliton theorem (equivalently Theorem~\ref{thm:asymp}) to this
complete controlled limit, with $\tau _{j}^{-}\downarrow 0$ and $\tau
_{j}^{+}\uparrow \infty $. For each $j$ choose $i(j)$, increasing in $j$,
such that the two rescalings of the flows based at $q_{i(j)}$ approximate
those of $h_{\infty }$ within $1/j$ in $C^{j}$ on balls of radius $j$ and on
$[-j,-1/j]$, and such that $f_{Z}(q_{j}^{\pm })>j$ for $q_{j}^{\pm }=\Psi
_{-\tau _{j}^{\pm }}(q_{i(j)})$; this is possible because each $\Psi _{-\tau
_{j}^{\pm }}$ is a fixed proper diffeomorphism and $q_{i}$ escapes. The
exact pointed isometries $(Z,(\tau _{j}^{\pm })^{-1}h(-\tau _{j}^{\pm
}),q_{i(j)})\cong (Z,h,q_{j}^{\pm })$ give the auxiliary convergences. The
maps $\Psi _{t}$ preserve gradient trajectories and the complex structure,
Theorem~\ref{thm:splitting} gives the complex-line splittings, and the
rigidity clause of Naber's theorem identifies $h_{\infty }(t)$ with its
asymptotic solitons when the two pointed limits are isometric.
\end{proof}

\begin{corollary}
\label{C41app}Assume the hypotheses of Theorem~\ref{P54}, including
completeness and $\dim M=5$. Let $\gamma$ be an integral curve of
$\nabla^T f$ whose projection satisfies $f_Z\to\infty$ along its forward end.
There is an escaping sequence $q_{i}\in \pi (\gamma )$ such that $(Z,h,q_{i})
$ converges to a complete noncollapsed shrinking K\"{a}hler--Ricci soliton.
Its universal cover is
\begin{equation*}
(\mathbb{C},h_{\mathrm{E}})\times (N,k),
\end{equation*}
where $(N,k)$ is a complete shrinking K\"{a}hler--Ricci soliton of complex
dimension one; by the classification of complete two-dimensional gradient
shrinkers, $(N,k)$ is the flat plane or the round sphere of Gaussian
curvature $\lambda=\tfrac12$. The corresponding local transverse limit on
$M$ is its local contactization.
\end{corollary}

\begin{proof}
Apply \cite[Corollary~4.1]{Nab} to the complete smooth quotient
$(Z,h,f_{Z})$, which has bounded curvature and is noncollapsed by
Corollary~\ref{cor:controlled}; Naber's corollary gives an escaping sequence
on the prescribed trajectory
with a complete shrinking limit and a real-line splitting. The parallel
complex structure, Theorem~\ref{thm:splitting} and de Rham splitting give
the flat complex line, and the soliton equation restricted to the
complementary factor makes it a complete shrinking K\"{a}hler--Ricci soliton
of complex dimension one. Since $N$ is simply connected and oriented, the
classification recalled in the proof of Theorem~\ref{thm:splitting}(3)
leaves only the flat plane and the round sphere, whose Gaussian curvature
is then the soliton constant $\lambda=\tfrac12$.
\end{proof}

\subsection{Local Splitting}

We now treat Reeb foliations of arbitrary type.

\begin{theorem}
\label{thm:splitting}\label{lA1}Let $p_{i}\in M$ with $f(p_{i})\rightarrow
\infty $ on a complete five-dimensional gradient shrinking Sasaki--Ricci
soliton with bounded transverse curvature.

\begin{enumerate}
\item Every developed local transverse limit (Definition~\ref{def:developed})
at these basepoints has orthonormal parallel fields $V$, $JV$.
Consequently it splits locally as
\begin{equation}
(U,h_{\mathrm{E}})\times (\Sigma ,h_{\Sigma }),\qquad U\subset \mathbb{C},
\label{eq:local-split}
\end{equation}
with its product K\"{a}hler structure. The assertion also holds transversely
in a smooth pointed limit of full Sasakian structures. Note that it is not
an ambient Riemannian product splitting.

\item Suppose a nonflat product limit satisfies the local soliton equation.
That is, on a neighborhood containing its product chart, there is a smooth
function $F$ with
\begin{equation}
\mathrm{Ric}(h_{\lim })+\operatorname{Hess}F=\lambda h_{\lim }.  \label{LS}
\end{equation}
Then
\begin{equation*}
F(z,y)=\frac{\lambda }{2}|z|^{2}+\ell (z)+\psi (y),\qquad
\operatorname{Hess}_{h_{\Sigma }}\psi =(\lambda -K)h_{\Sigma },
\end{equation*}
where $\ell $ is affine and $K$ is the Gaussian curvature of $h_{\Sigma }$.
In particular,
\begin{equation}
dK=K\,d\psi .  \label{eq:dK}
\end{equation}

\item In the nonflat situation of part (2), if the surface factor,
together with the potential $\psi$, extends to a complete smooth connected
two-dimensional gradient shrinker with the same constant $\lambda $, then
its local metric has
\begin{equation*}
K\equiv \lambda .
\end{equation*}
It is locally the round sphere factor of $\mathbb{C}\times \mathbb{CP}^{1}$.
More generally, when a complete smooth global K\"{a}hler limit is available,
the parallel fields give a global splitting
\begin{equation*}
(\mathbb{C},h_{\mathrm{E}})\times (N,k)
\end{equation*}
of its simply connected cover.
\end{enumerate}
\end{theorem}

\begin{proof}
(1) Set $a_{i}=|\nabla ^{T}f|(p_{i})$ and $u_{i}=(f-f(p_{i}))/a_{i}$. The
Hamilton identity \eqref{eq:canonical} and the Lipschitz estimate for
$\sqrt{f}$ imply, on fixed-radius charts,
\begin{equation*}
a_{i}\rightarrow \infty ,\qquad \frac{f}{f(p_{i})}\rightarrow 1,\qquad
|\nabla ^{T}u_{i}|\rightarrow 1,\qquad
\operatorname{Hess}^{T}u_{i}=\frac{\lambda
h-\mathrm{Ric}^{T}}{a_{i}}\rightarrow 0.
\end{equation*}
Derivatives of the Hessian are controlled by
Proposition~\ref{prop:derivatives}. The normalized functions therefore
converge smoothly,
after passing to a subsequence, on a smaller transverse chart to a function
$u_{\infty }$ with unit parallel gradient $V$. Since $J$ is parallel, $JV$ is
also parallel, unit, and orthogonal to $V$. The two-dimensional real
distribution they span and its orthogonal complement are parallel and
$J$-invariant. Local de Rham decomposition proves \eqref{eq:local-split}. For
full Sasakian convergence the mixed ambient Hessian components of $u_{i}$
are controlled algebraically by $du_{i}$ and the Sasakian tensors; they are
not required to tend to zero.

(2) The mixed Hessian components of $F$ vanish, so $F$ separates into
functions of the two factors. Its Euclidean Hessian is $\lambda
h_{\mathrm{E}}$, giving the asserted quadratic plus affine term. In two
dimensions $\mathrm{Ric}=Kh_{\Sigma }$, giving the remaining equation. The
scalar
identity $d(2K)=2\,\mathrm{Ric}(\nabla \psi ,\cdot )$, the analogue of
\eqref{eq:pre-identities}, is exactly \eqref{eq:dK}.

(3) The classification of complete smooth two-dimensional gradient shrinking
Ricci solitons leaves the flat plane and the round sphere, up to allowable
quotients; see Hamilton \cite{Ham} and Bernstein--Mettler \cite{BM}, and also
\cite{PW,MW2}.
Nonflatness eliminates the plane and its flat quotients. In the oriented
smooth case no nontrivial finite orientation-preserving isometry group acts
freely on $\mathbb{S}^{2}$, so the complete model is the sphere. The soliton
constant fixes its Gaussian curvature at $\lambda $. The global splitting
assertion is the complete de Rham theorem applied on the simply connected
cover.
\end{proof}

\begin{remark}
On the quotient of a regular structure a complete smooth global K\"{a}hler
limit is available, and the global splitting of Theorem~\ref{thm:splitting}(3)
applies; see \cite[Theorem~\IIref{T51}{4.3}]{PartII} for the regular case and
\cite[Theorem~\IIref{thm:split-structure}{6.1}]{PartII} for global splittings
of the total space.
\end{remark}

\begin{corollary}
\label{cor:flatlimit}Under the hypotheses of Theorem~\ref{thm:curvature}, if
$f(p_{i})\rightarrow \infty $ and $R^{T}(p_{i})\rightarrow 0$, every
developed local limit at $p_{i}$ is transversely flat, and its Sasakian
contactization is locally Heisenberg.
\end{corollary}

\begin{proof}
Assume that $M$ is not transversely flat, so that $R^{T}>0$ and
$|\mathrm{Rm}^{T}|\leq CR^{T}$ by Theorem~\ref{thm:curvature}. No gradient
estimate for $\log R^{T}$ is needed; we use the parabolic strong maximum
principle on the associated flow. After a fixed $D$-homothety we may assume
$\lambda =\frac{1}{2}$. Let $g^{T}(t)=(-t)\Psi _{t}^{\ast }g^{T}$ be the flow
of Proposition~\ref{lem:selfsim}, so that $g^{T}(-1)=g^{T}$. By
Proposition~\ref{prop:derivatives} and the evolution equations, all space
and time derivatives of $\mathrm{Rm}^{T}(g^{T}(t))$ are bounded on
$M\times \lbrack -2,-1]$, and the metrics $g^{T}(t)$, $t\in \lbrack -2,-1]$,
are uniformly equivalent. Restricting them to transversals through $p_{i}$
in the charts of Definition~\ref{def:developed} and passing to a
subsequence, we obtain a local K\"{a}hler--Ricci flow $h_{\infty }(t)$,
$t\in \lbrack -2,-1]$, on a transverse ball $B$ centered at the basepoint
$o$, such that $h_{\infty }(-1)$ is the developed local transverse limit.
Its scalar curvature $R_{\infty }\geq 0$ satisfies
\begin{equation*}
\partial _{t}R_{\infty }=\Delta R_{\infty }+2|\mathrm{Ric}(h_{\infty
})|^{2}\geq \Delta R_{\infty },\qquad R_{\infty }(o,-1)=\lim_{i\rightarrow
\infty }R^{T}(p_{i})=0.
\end{equation*}
By the strong maximum principle for nonnegative supersolutions of the heat
equation, $R_{\infty }\equiv 0$ on $B\times \lbrack -2,-1]$. The inequality
$|\mathrm{Rm}^{T}|\leq CR^{T}$ is invariant under the parabolic rescalings
and diffeomorphisms defining $g^{T}(t)$, so it holds for every $g^{T}(t)$
with the same constant and passes to the limit. Hence
$\mathrm{Rm}(h_{\infty }(-1))\equiv 0$ on $B$, that is, the developed local
transverse limit is flat. On a flat transverse chart $d\eta =2\omega
_{\mathrm{E}}$, and
Proposition~\ref{prop:contactization} identifies its connection with the
Heisenberg one \eqref{eq:heisenberg} after a change of the local Reeb
coordinate. The globally flat branch is immediate.
\end{proof}

\begin{remark}
The normalized function $u_{i}$ used in the proof of Theorem
\ref{thm:splitting} satisfies a limiting equation with zero Hessian. It is not
a
limiting shrinking potential. The unnormalized potentials $f-f(p_{i})$
generally have unbounded gradients and need not converge. Therefore the
existence of a shrinking potential on the surface factor requires a separate
condition. This input, the resulting two-dimensional equation, and an
explicit irregular example whose surface factor is not round are discussed in
\cite[Theorem~\IIref{thm:irregular-local}{5.1} and
Example~\IIref{ex:counterexample}{5.2}]{PartII}.
\end{remark}

\subsection{Transverse Structure at the End}

\label{sec:end}

Let $(M^{5},\xi ,\eta ,\Phi ,g,f)$ be connected, complete and noncompact,
and satisfy the full real transverse equation
\begin{equation*}
Q+\operatorname{Hess}^{T}f=\lambda h,\qquad \xi f=0,\qquad h=g|_{D},\quad D=\ker
\eta ,\quad \lambda >0.
\end{equation*}
Write
\begin{equation*}
T=\mathrm{Rm}^{T},\quad Q=\mathrm{Ric}^{T},\quad s=R^{T}=R_{g}+4,\qquad \omega
^{T}=\tfrac{1}{2}d\eta .
\end{equation*}
The scalar relation is \eqref{eq:pre-ric} with $n=2$. The normalization of
Sections~\ref{sec:identities} and \ref{sec:curvature} is $\lambda =6$; here
we keep $\lambda $ explicit. Choose the compatible sign for the transverse
complex structure $J$ so $\omega ^{T}(U,V)=h(JU,V)$. Our curvature
components are $T_{ijkl}=\left\langle
R^{T}(e_{k},e_{l})e_{j},e_{i}\right\rangle $. All lower-case indices in this
section run from one to four.

We work on the exhaustion space $\mathcal{X}$ of Lemma~\ref{lem:reeb}. When
the Reeb closure is compact, infinity is measured on $M$ by distance to a
fixed compact orbit; when the action is proper $\mathbb{R}$, it is measured
on the complete quotient $Z$. The normalized potential is proper on
$\mathcal{X}$ (Lemma~\ref{lem:growth}). We use ``transverse
infinity'' to mean $f\rightarrow \infty $, not arbitrary
escape along an infinite Reeb line. The identities are
\begin{align}
s+\Delta _{B}f& =4\lambda , & \nabla ^{T}s& =2Q(\nabla ^{T}f,\cdot ),
\label{eq:trace} \\
s+|\nabla ^{T}f|^{2}& =2\lambda f, & L_{f}f& =4\lambda
-2\lambda f,  \label{eq:ham} \\
L_{f}s& =2\lambda s-2|Q|^{2}, & \sum_{i}T_{ijkl;i}&
=Q_{jl;k}-Q_{jk;l}=\sum_{i}T_{ijkl}f_{;i}.  \label{eq:scalarbianchi}
\end{align}

Assume $s$ is bounded above. By the identities and estimates of
Sections~\ref{sec:identities} and~\ref{sec:curvature}
(Theorem~\ref{thm:curvature} and,
for the derivative, Proposition~\ref{prop:gradT}), either $s=T=0$ or $s>0$
everywhere. In the nonflat case, 
\begin{equation}
sf\geq c>0,\quad |T|+|Q|\leq Cs\quad \text{and}\quad |\nabla
^{T}T|\leq Cs,  \label{eq:input}
\end{equation}
and
\begin{equation}
c(1+r)^{2}\leq f\leq C(1+r)^{2}\quad \text{and}\quad |\mathrm{Rm}_{g}|\leq
C(1+s).  \label{eq:input-growth}
\end{equation}
Here $r$ is distance to a fixed point of $\mathcal{X}$; replacing it by
distance to a fixed compact Reeb orbit closure changes it by a bounded
amount. Constants may depend on a fixed compact sublevel and the soliton,
not only on $\sup s$. Section~\ref{sec:curvature} proves these estimates in
the normalization $\lambda =6$; a fixed transverse $D$-homothety gives the displayed versions. Bounds for all higher transverse derivatives follow from
Proposition~\ref{prop:derivatives}.

Fix a constant $F\geq \sup s/\lambda $. By \eqref{eq:ham},
\begin{equation*}
|\nabla ^{T}f|^{2}=2\lambda f-s\geq \lambda f>0\qquad \text{on }\{f\geq F\};
\end{equation*}
in the proofs below $F$ may be enlarged finitely many times. On $\{f\geq F\}$
define the $J$-invariant real two-plane
\begin{equation*}
E=\operatorname{span}\{\nabla ^{T}f,J\nabla ^{T}f\}^{\perp }\subset D,
\end{equation*}
and let $P_{E}$ be the orthogonal projection of $\Lambda ^{2}D$ onto its
one-dimensional subspace $\Lambda ^{2}E$. The transverse curvature operator
$\mathcal{R}^{T}$ is normalized so that a surface of Gaussian curvature $K$
has curvature operator $K$ on its oriented area form.

\begin{theorem}
\label{thm:general}Assume, as above, that $s=R^{T}$ is bounded from above. In
the nonflat case one has, on $\{f\geq F\}$,
\begin{equation}
\Bigl\Vert\mathcal{R}^{T}-\frac{s}{2}P_{E}\Bigr\Vert_{\mathrm{op}}\leq
\frac{Cs}{\sqrt{f}}\leq C\sqrt{\frac{s}{f}}\leq \frac{C}{1+r}.
\label{eq:rankbounded}
\end{equation}
Equivalently, $\Vert s^{-1}\mathcal{R}^{T}-\frac{1}{2}P_{E}\Vert
_{\mathrm{op}}\leq Cf^{-1/2}$.
If in addition $s\rightarrow 0$ uniformly as $f\rightarrow \infty $, then
\begin{equation}
\frac{c}{f}\leq s\leq \frac{C}{f},\qquad |(\nabla ^{T})^{j}T|\leq
C_{j}f^{-1-j/2}\quad (j\geq 0),\qquad
\Bigl\Vert\mathcal{R}^{T}-\frac{s}{2}P_{E}\Bigr\Vert_{\mathrm{op}}\leq
Cf^{-2};  \label{eq:decay-f}
\end{equation}
equivalently, by \eqref{eq:input-growth},
\begin{align}
c(1+r)^{-2}\leq s& \leq C(1+r)^{-2}, & |(\nabla ^{T})^{j}T|& \leq
C_{j}(1+r)^{-2-j},  \label{eq:decay} \\
& & \Bigl\Vert\mathcal{R}^{T}-\frac{s}{2}P_{E}\Bigr\Vert_{\mathrm{op}}& \leq
C(1+r)^{-4}.  \label{eq:rankdecay}
\end{align}
In particular the transverse curvature operator is bounded below by
$-Cs(1+r)^{-1}\mathrm{Id}$, hence by $-C(1+r)^{-1}\mathrm{Id}$, in the
bounded-scalar case and by $-C(1+r)^{-4}\mathrm{Id}$ in the decaying branch.
Every developed local transverse limit
(Definition~\ref{def:developed}) along a sequence $p_{i}$ with
$f(p_{i})\rightarrow \infty $ consequently has nonnegative curvature
operator; its surface factor in \eqref{eq:local-split} has $K\geq 0$.
\end{theorem}

\begin{proof}
The Hamilton identity gives $|\nabla ^{T}f|^{2}=2\lambda f-s\geq \lambda f$
on $\{f\geq F\}$, by the choice of $F$. The flow of
$\nabla ^{T}f/|\nabla ^{T}f|^{2}$ changes $f$ at unit speed and is complete
over each finite interval of potential values by properness. It therefore
identifies the exterior $\{f\geq F\}$ with
$[F,\infty)\times \{f=F\}$ for large $F$; basicness makes the identification
Reeb-equivariant.

To prove \eqref{eq:rankbounded}, fix a point of $\{f\geq F\}$ and choose
\begin{equation*}
e_{4}=\frac{\nabla ^{T}f}{|\nabla ^{T}f|},\qquad
e_{3}=Je_{4},\qquad
e_{1}\in E\text{ a unit vector},\qquad
e_{2}=Je_{1}.
\end{equation*}
Since $J$ is orthogonal and skew and $E$ is the $J$-invariant orthogonal
complement of $\operatorname{span}\{e_{3},e_{4}\}$, this is an orthonormal basis
of $D$ with $E=\operatorname{span}\{e_{1},e_{2}\}$.

First, in this basis $f_{;i}=|\nabla ^{T}f|\,\delta _{i4}$, so the last
identity in \eqref{eq:scalarbianchi} reads
\begin{equation*}
|\nabla ^{T}f|\,T_{4jkl}
=Q_{jl;k}-Q_{jk;l}
\qquad \text{for all }j,k,l.
\end{equation*}
The right-hand side is the difference of two index permutations of
$\nabla ^{T}Q$, so squaring and summing over $j,k,l$ gives
$|\nabla ^{T}f|^{2}\sum_{j,k,l}T_{4jkl}^{2}
\leq 4|\nabla ^{T}Q|^{2}$.
Moreover, $Q$ is the contraction of $T$ over one pair of indices, hence
$\nabla ^{T}Q$ is the same contraction of $\nabla ^{T}T$, and the
Cauchy--Schwarz inequality over the four summands of the contraction
gives $|\nabla ^{T}Q|\leq 2|\nabla ^{T}T|$.
Together with $|\nabla ^{T}f|^{2}\geq \lambda f$, this yields
\begin{equation*}
\Bigl(\sum_{j,k,l}T_{4jkl}^{2}\Bigr)^{1/2}
\leq \frac{2|\nabla ^{T}Q|}{|\nabla ^{T}f|}
\leq \frac{4|\nabla ^{T}T|}{\sqrt{\lambda f}}.
\end{equation*}

Next, since $\nabla ^{T}J=0$, every curvature endomorphism $R^{T}(Z,W)$
commutes with $J$; as $J$ is orthogonal, this gives
\begin{equation*}
T(JX,JY,Z,W)=T(X,Y,Z,W).
\end{equation*}
Taking $X=e_{4}$ and replacing $Y$ by $-JY$, we obtain
\begin{equation*}
T(e_{3},Y,Z,W)=-T(e_{4},JY,Z,W).
\end{equation*}
As $J$ permutes the basis vectors up to sign
($Je_{1}=e_{2}$, $Je_{2}=-e_{1}$, $Je_{3}=-e_{4}$, $Je_{4}=e_{3}$),
the numbers $T_{3jkl}$ are, up to sign and a permutation of the index $j$,
the numbers $T_{4jkl}$. Hence
$\sum_{j,k,l}T_{3jkl}^{2}=\sum_{j,k,l}T_{4jkl}^{2}$.

Furthermore, the symmetries
$T_{ijkl}=-T_{jikl}=-T_{ijlk}=T_{klij}$
move any prescribed index to the first position, up to sign.
Hence the preceding estimates give
\begin{equation}
|T_{ijkl}|
\leq \varepsilon
:=\frac{4|\nabla ^{T}T|}{\sqrt{\lambda f}}
\qquad
\text{whenever }
\{i,j,k,l\}\cap \{3,4\}\neq \emptyset.
\label{eq:rankcomponent}
\end{equation}

To pass to the curvature operator, give $\Lambda ^{2}D$ the inner product
for which $\{e_{i}\wedge e_{j}\}_{i<j}$ is orthonormal.
The normalization of $\mathcal{R}^{T}$ fixed before the theorem means
\begin{equation*}
\langle \mathcal{R}^{T}(e_{i}\wedge e_{j}),e_{k}\wedge e_{l}\rangle
=T_{ijkl},
\end{equation*}
since
$T_{1212}=\langle R^{T}(e_{1},e_{2})e_{2},e_{1}\rangle$
is the sectional curvature of $\operatorname{span}\{e_{1},e_{2}\}$.
The matrix of $\mathcal{R}^{T}$ in this basis is a symmetric
$6\times 6$ matrix, and every entry except the one indexed by
$(e_{1}\wedge e_{2},e_{1}\wedge e_{2})$
involves an index $3$ or $4$.
Since $\Lambda ^{2}E$ is spanned by $e_{1}\wedge e_{2}$, we may write
\begin{equation*}
\mathcal{R}^{T}=kP_{E}+\mathcal{E},\qquad k=T_{1212},
\end{equation*}
where the matrix of $\mathcal{E}$ vanishes at
$(e_{1}\wedge e_{2},e_{1}\wedge e_{2})$
and its other $35$ entries are bounded by $\varepsilon$,
by \eqref{eq:rankcomponent}.
Bounding the operator norm by the Hilbert--Schmidt norm gives
$\Vert \mathcal{E}\Vert _{\mathrm{op}}\leq \sqrt{35}\,\varepsilon$,
and the five remaining diagonal entries give
$|\operatorname{tr}\mathcal{E}|\leq 5\varepsilon$.

Moreover, by the definition of the scalar curvature,
\begin{equation*}
s=\sum_{i,j}T_{ijij}
=2\sum_{i<j}T_{ijij}
=2\,\operatorname{tr}\mathcal{R}^{T}
=2k+2\,\operatorname{tr}\mathcal{E}.
\end{equation*}
Hence
$|k-\frac{s}{2}|=|\operatorname{tr}\mathcal{E}|
\leq 5\varepsilon$,
and
\begin{equation*}
\Bigl\Vert \mathcal{R}^{T}-\frac{s}{2}P_{E}\Bigr\Vert _{\mathrm{op}}
\leq
\Bigl|k-\frac{s}{2}\Bigr|
+\Vert \mathcal{E}\Vert _{\mathrm{op}}
\leq
(5+\sqrt{35})\,\varepsilon
\leq
\frac{44}{\sqrt{\lambda}}\,
\frac{|\nabla ^{T}T|}{\sqrt{f}}
\qquad \text{on }\{f\geq F\}.
\tag{$\ast$}
\label{eq:rank-pointwise}
\end{equation*}
This pointwise inequality uses only \eqref{eq:scalarbianchi}, the
K\"{a}hler symmetries and the choice of $F$; it does not use any
curvature bound.

We now apply Proposition~\ref{prop:gradT}, which gives
$|\nabla ^{T}T|\leq Cs$ on $M$.
(That proposition is stated in the normalization $\lambda=6$;
a fixed $D$-homothety $g^{T}\mapsto \tau g^{T}$ multiplies
$|\nabla ^{T}T|$ by $\tau^{-3/2}$ and $s$ by $\tau^{-1}$
and leaves $f$ unchanged, so the inequality holds at every scale
with another constant.)
Inserting it into \eqref{eq:rank-pointwise} gives the first inequality
of \eqref{eq:rankbounded}. Since $s\leq S_{0}:=\sup s$,
\begin{equation*}
\frac{s}{\sqrt{f}}
=\sqrt{s}\,\sqrt{\frac{s}{f}}
\leq \sqrt{S_{0}}\,\sqrt{\frac{s}{f}},
\qquad
\sqrt{\frac{s}{f}}
\leq \sqrt{\frac{S_{0}}{f}}
\leq \frac{C}{1+r},
\end{equation*}
the last step by \eqref{eq:input-growth}; these are the second and
third inequalities.
Dividing by $s>0$ gives the equivalent form
$\Vert s^{-1}\mathcal{R}^{T}-\frac{1}{2}P_{E}\Vert _{\mathrm{op}}
\leq Cf^{-1/2}$.
Finally, since $P_{E}\geq 0$ and $s>0$,
\begin{equation*}
\mathcal{R}^{T}
=\frac{s}{2}P_{E}
+\Bigl(\mathcal{R}^{T}-\frac{s}{2}P_{E}\Bigr)
\geq -\frac{Cs}{\sqrt{f}}\,\mathrm{Id}
\geq -\frac{Cs}{1+r}\,\mathrm{Id},
\end{equation*}
which is the first lower bound stated in the theorem.

We now assume that $s\rightarrow 0$ uniformly as $f\rightarrow \infty$
and prove \eqref{eq:decay-f}.
For $p>0$, by \eqref{eq:inversepowers-general} with $m=4$,
\begin{equation}
L_{f}(f^{-p})
=2p\lambda f^{-p}
+2p\lambda(p-1)f^{-p-1}
-p(p+1)sf^{-p-2}.
\label{eq:inversepowers}
\end{equation}
By $|Q|\leq Cs$, there is $c_{1}>0$ such that
$L_{f}s\geq 2\lambda s-c_{1}s^{2}$.
Choose $F$ larger if necessary so that on $\{f\geq F\}$,
\begin{equation*}
s<\frac{\lambda}{2c_{1}}
\quad \text{and}\quad
L_{f}(f^{-2})\geq 4\lambda f^{-2}.
\end{equation*}
Put $a=\lambda F/c_{1}$, $b=c_{1}/(2\lambda)$, and
$w=s-af^{-1}+ba^{2}f^{-2}$.
Equations~\eqref{eq:inversepowers} and
$L_{f}(f^{-1})\leq 2\lambda f^{-1}$ give
\begin{equation*}
L_{f}w
\geq 2\lambda w-c_{1}(s-af^{-1})(s+af^{-1})
\geq [2\lambda-c_{1}(s+af^{-1})]w.
\end{equation*}
The bracket is at least $\lambda/2$, $w<0$ on $f=F$, and
$w\rightarrow 0$ at infinity.
The exterior maximum principle proves $w\leq 0$, hence $s\leq a/f$.
Together with $sf\geq c$, this proves the scalar bounds and
$|T|\leq C/f$.

For the first derivative, \eqref{eq:deriv} and $|T|\leq C/f$ show that
$\Xi=|\nabla ^{T}T|$ satisfies
\begin{equation*}
L_{f}\Xi\geq (3\lambda-C_{1}/f)\Xi.
\end{equation*}
Also $\Xi\rightarrow 0$ by \eqref{eq:input}.
For
\begin{equation*}
\mathscr{A}=L_{f}-3\lambda+C_{1}/f
\quad \text{and}\quad
\varphi=f^{-3/2}-Bf^{-5/2},
\end{equation*}
the exact expansion is
\begin{align*}
\mathscr{A}\varphi
&=
\left(C_{1}+\tfrac{3}{2}\lambda-2\lambda B\right)f^{-5/2}
-B\left(C_{1}+\tfrac{15}{2}\lambda\right)f^{-7/2}\\
&\hspace{8mm}
-\tfrac{15}{4}sf^{-7/2}
+\tfrac{35}{4}Bsf^{-9/2}.
\end{align*}
Choose $B$ so the first coefficient is negative.
On a sufficiently large exterior, $\varphi>0$,
\begin{equation*}
\mathscr{A}\varphi\leq 0
\quad \text{and}\quad
3\lambda-C_{1}/f>0.
\end{equation*}
Comparing $\Xi$ with a sufficiently large multiple of $\varphi$,
first on the inner boundary and then at infinity, gives
$\Xi\leq Cf^{-3/2}$.
At a hypothetical positive maximum the norm is smooth, so Kato's
inequality suffices.

For higher derivatives we argue by induction on $j$, the cases $j=0,1$
having been proved.
Let $\Xi_{j}=|(\nabla ^{T})^{j}T|$, bounded a priori by
Proposition~\ref{prop:derivatives}, and assume
$\Xi_{i}\leq C_{i}f^{-1-i/2}$ for $i<j$.
In \eqref{eq:deriv-k}, the terms with a factor $T$ are bounded by
$Cf^{-1}\Xi_{j}$, and the remaining products by $Cf^{-(j+4)/2}$.
Hence, by Kato's inequality,
\begin{equation*}
\bigl(L_{f}-(j+2)\lambda+C/f\bigr)\Xi_{j}
\geq -Cf^{-(j+4)/2}.
\end{equation*}
Put $p=1+j/2$ and $\mathscr{A}=L_{f}-2p\lambda+C/f$.
Formula \eqref{eq:inversepowers} implies that, for a fixed sufficiently
large $\beta$ and on a fixed sufficiently large exterior,
\begin{equation*}
\mathscr{A}\bigl(Gf^{-p}-G\beta f^{-p-1}\bigr)
\leq -Cf^{-p-1}
\end{equation*}
when the constant $G$ is large enough.
The bracket is positive there; choose $G$ to dominate $\Xi_{j}$
on the inner boundary.
For any $\varepsilon>0$ add $\varepsilon f^{\alpha}$ with $\alpha>0$
to the barrier.
On enlarging the exterior, $\mathscr{A}(f^{\alpha})<0$, while
boundedness of $\Xi_{j}$ forces the difference to be negative at infinity.
A positive interior maximum is impossible.
Letting $\varepsilon\downarrow 0$ proves
$\Xi_{j}\leq C_{j}f^{-1-j/2}$, which is the induction step.

Returning to \eqref{eq:rankcomponent} and \eqref{eq:rank-pointwise},
the bound $|\nabla ^{T}T|\leq C_{1}f^{-3/2}$ just proved shows that
every component with an index $4$ or $3$ is $O(f^{-2})$, and that
the left-hand side of \eqref{eq:rank-pointwise} is $O(f^{-2})$.
Equation~\eqref{eq:input-growth} converts the bounds to
\eqref{eq:decay}--\eqref{eq:rankdecay}.

Finally, $s\geq 0$ and the error in \eqref{eq:rankbounded} tends to zero
on every fixed-radius chart escaping in potential.
Curvature converges smoothly on developed limits, which proves the
last assertion.
\end{proof}

\begin{corollary}
\label{prop:conditional}Let $(M^{5},\xi ,\eta ,\Phi ,g,f)$ be a connected
complete gradient shrinking Sasaki--Ricci soliton with
\begin{equation*}
R^{T}\leq C,
\end{equation*}
and let $(p_{i})$ be any sequence in $M$ with $f(p_{i})\rightarrow \infty $.
Then every developed local transverse limit at $(p_{i})$
(Definition~\ref{def:developed}) splits locally as
\begin{equation}
(U,h_{\mathrm{E}})\times (\Sigma ^{2},h_{\Sigma }),\qquad U\subset
\mathbb{C},\qquad 0\leq K_{\Sigma }\leq C.  \label{eq:conditionalproduct}
\end{equation}
If $R^{T}(p_{i})\rightarrow 0$, the limiting transverse metric is flat and
the limiting Sasakian metric is locally Heisenberg. In particular, all of
this holds for every sequence $(p_{i})$ on an integral curve of $\nabla
^{T}f$ which tends to infinity.
\end{corollary}

\begin{proof}
We check that each step uses only $f(p_{i})\rightarrow \infty $. Hypotheses
and conclusions are invariant under a fixed $D$-homothety, so we may assume
the normalization of Theorem~\ref{thm:curvature}, and we normalize $f$ by
\eqref{eq:canonical}. If $R^{T}\equiv 0$, then $\mathrm{Rm}^{T}\equiv 0$ by
Proposition~\ref{prop:flat}, every developed local limit is transversely
flat, and the assertions hold with $K_{\Sigma }\equiv 0$, the last one by
Corollary~\ref{cor:flatlimit}. Assume therefore $R^{T}>0$.

(i) \emph{Existence of the limits.} By Theorem~\ref{thm:curvature},
$|\mathrm{Rm}^{T}|\leq CR^{T}\leq C$. Hence the transverse curvature is
bounded, and Proposition~\ref{prop:unwrap} shows that $\sup_{M}|\nabla
^{k}\mathrm{Rm}_{g}|<\infty $ for every $k$ and that every sequence of
basepoints, in particular $(p_{i})$, admits developed local limits.

(ii) \emph{The charts escape in potential.} Since $f$ is basic, $|\nabla
_{g}\sqrt{f}|=|\nabla ^{T}\sqrt{f}|\leq \sqrt{\lambda /2}$ by \eqref{eq:Lf-f},
so $\sqrt{f}$ is $\sqrt{\lambda /2}$-Lipschitz for the
distance of $g$. The developed chart of radius $\rho \leq \rho _{0}$ at
$p_{i}$ is mapped by $\exp _{p_{i}}$ into the metric ball $B(p_{i},\rho
_{0})$. Therefore
\begin{equation*}
\sqrt{f}\geq \sqrt{f(p_{i})}-\sqrt{\lambda /2}\,\rho _{0}\longrightarrow
\infty
\end{equation*}
uniformly on these charts.

(iii) \emph{The local splitting.} Theorem~\ref{thm:splitting}(1) is stated
for an arbitrary sequence with $f(p_{i})\rightarrow \infty $, and its proof
uses only $a_{i}=|\nabla ^{T}f|(p_{i})\rightarrow \infty $, which follows
from $a_{i}^{2}=2\lambda f(p_{i})-R^{T}(p_{i})\geq 2\lambda f(p_{i})-C$. It
gives orthonormal parallel fields $V$ and $JV$ on the limit, where $V$ is
the limit of $\nabla ^{T}f/a_{i}$, and the local product
\eqref{eq:conditionalproduct}, whose surface factor is tangent to
$\{V,JV\}^{\perp }$.

(iv) \emph{The bounds for $K_{\Sigma }$.} By (ii) and \eqref{eq:rankbounded},
$\Vert \mathcal{R}^{T}-\frac{s}{2}P_{E}\Vert _{\mathrm{op}}\leq
Cs/\sqrt{f}\leq C/\sqrt{f}\rightarrow 0$ uniformly on the charts, and
$\frac{s}{2}P_{E}\geq 0$. The pulled-back curvature operators converge
smoothly, so the
curvature operator of the limit is nonnegative; this is the last assertion
of Theorem~\ref{thm:general}, which is also stated for arbitrary sequences.
For the local product \eqref{eq:conditionalproduct} the curvature operator
is $K_{\Sigma }$ on $\Lambda ^{2}T\Sigma $ and vanishes on the orthogonal
complement of $\Lambda ^{2}T\Sigma $. Hence $K_{\Sigma }\geq 0$. For the upper bound, the scalar curvature of this local product is
$2K_\Sigma$. Passing the original bound $R^T\leq C$ to the limit gives
$2K_\Sigma\leq C$, and hence the stated bound $K_\Sigma\leq C$
with the same scalar-curvature constant.

(v) \emph{Flat limits.} If $R^{T}(p_{i})\rightarrow 0$, then
Corollary~\ref{cor:flatlimit}, which is also stated for an arbitrary
sequence with $f(p_{i})\rightarrow \infty $ and $R^{T}(p_{i})\rightarrow 0$,
shows that the limiting transverse metric is flat and that the limiting
Sasakian structure is locally Heisenberg.

None of these steps uses that the $p_{i}$ lie on one integral curve of $\nabla
^{T}f$. This proves the corollary, and in particular its special case
for sequences on an integral curve which tends to infinity.
\end{proof}

\begin{remark}
A gradient trajectory matters for the following finer statements in
setting (S). Along a projected gradient trajectory, moving basepoints
can be traded for times of the quotient flow at a fixed basepoint:
if $q_i=\Psi_{-\tau_i}(x_0)$ with $\tau_i\downarrow0$, then
$(Z,h,q_i)$ is isometric to $(Z,\tau_i^{-1}h(-\tau_i),x_0)$, as in the
proof of Theorem~\ref{thm:naber-split}. These are parabolic dilations
of a single ancient flow about the time $0$ at the fixed point $x_0$;
the scale factors $\tau_i^{-1}$ tend to infinity.
This is what produces the shrinking soliton limits of
Corollary~\ref{C41app}, and, in \cite[Theorem~\IIref{T51}{4.3}]{PartII},
the dichotomy $R^{T}\rightarrow 0$ or $R^{T}\rightarrow 2\lambda $ along
the trajectory together with the uniqueness of the limit. For an arbitrary
sequence with $f(p_{i})\rightarrow \infty $, different subsequences may
have different limits.

The surface factor in Corollary~\ref{prop:conditional} need not satisfy
a local soliton equation, and even when it does, it need not be round: a nonregular quasi-regular or irregular foliation does not in general
supply a complete smooth global surface factor to which the
classification of smooth surface shrinkers applies.
The examples in \cite[Remark~\IIref{rem:weighted-example}{4.3} and
Example~\IIref{ex:counterexample}{5.2}]{PartII} have limits along escaping
gradient trajectories whose surface factors have nonconstant curvature.
Under the regular cylindrical-trajectory hypotheses of
\cite[Theorem~\IIref{T51}{4.3}]{PartII}, the pointed Sasakian structures
converge to a complete Sasakian structure whose quotient splits globally;
see also \cite[Theorem~\IIref{thm:regular}{4.5} and the paragraph following
Remark~\IIref{rem:cylmodel}{4.2}]{PartII}.
\end{remark}

\begin{remark}
\label{rem:rank-improved}The proof of \eqref{eq:rankbounded} reduces the
rank-one estimate to the
pointwise inequality \eqref{eq:rank-pointwise}, so every derivative bound
for $T$ gives a rank-one estimate. The bound $|\nabla ^{T}T|\leq C\sqrt{s}$
of \eqref{eq:gradTsqrt} alone gives only
$C\sqrt{s/f}$; Proposition~\ref{prop:gradT} gives $Cs/\sqrt{f}$; and in the
decaying branch $|\nabla ^{T}T|\leq Cf^{-3/2}$ gives $Cf^{-2}$. The
improvement matters where $s$ is small. The weaker bound \eqref{eq:gradTsqrt} gives only $\Vert
s^{-1}\mathcal{R}^{T}-\frac{1}{2}P_{E}\Vert _{\mathrm{op}}\leq C(sf)^{-1/2}$,
which is bounded because $sf\geq c$ but need not tend to zero, whereas
\eqref{eq:rankbounded} gives $Cf^{-1/2}$. Thus the normalized curvature
operator $s^{-1}\mathcal{R}^{T}$ is asymptotically of rank one at a rate
which does not depend on the size of $s$, and the negative part of
$\mathcal{R}^{T}$ is small relative to $s$: $\mathcal{R}^{T}\geq
-Cs(1+r)^{-1}\mathrm{Id}$.
\end{remark}

\appendix

\section{Transverse Reduced Geometry for Controlled Sasaki--Ricci Flows}
\label{app:reduced}

The material of this appendix is standard; we include it for the
reader's convenience and to fix the normalizations. Throughout this
appendix we use the single cover and quotient fixed in setting (S). We write $M$ for the chosen total space and
$\pi:M\to Z$ for its time-independent quotient map. The fixed
complex manifold $Z$ carries a complete smooth
K\"{a}hler--Ricci flow $h(t)$ with
$\partial_t h=-2\mathrm{Ric}_{h}$. Its real dimension is $m=2n$.
For the ancient results $t<0$ and $|\mathrm{Rm}(h(t))|\le C/|t|$.
For the finite-time statements it suffices to work on a compact time
interval on which this quotient flow has bounded curvature.
All distances, curve variations, exponential maps, cut loci, and
integrals in this appendix are first defined on this same $Z$.
Superscript $T$ denotes their quotient expression or basic pullback;
$\nabla^T$, $\Delta_B$, and $R^T$ therefore
mean $\nabla_h$, $\Delta_h$, and $R_h$ in the calculations below.
This is Perelman's normalization \cite[Section~7]{Per1}; the complete
bounded-curvature theory is given in \cite[Chapters~6--7]{MT} (see also
\cite{CaZh}) and the ancient estimates in \cite[Section~2]{Nab}.

\subsection{The Transverse \texorpdfstring{$\mathcal{L}$}{L}-Geodesic}

Fix a regular base time $T$ and a point $\bar x\in Z$. When a basepoint
$x\in M$ is specified, $\bar x=\pi(x)$; similarly $\bar y=\pi(y)$.
Set $h_\tau=h(T-\tau)$, so $\partial_\tau h_\tau=2\mathrm{Ric}_{h_\tau}$.

\begin{definition}
An admissible quotient curve on $[0,\bar\tau]$ is a continuous map
$\gamma:[0,\bar\tau]\to Z$ for which
$\alpha(s)=\gamma(s^2)$ is piecewise smooth on $[0,\sqrt{\bar\tau}]$.
We also use its $W^{1,2}$ completion in the variable $s$ in minimization.
Zero-speed curves, in particular constant curves, are allowed.
Write $X=d\gamma/d\tau$ and $A=d\alpha/ds=2sX$.
The velocity $X$ may have the integrable $O(\tau^{-1/2})$ singularity
at the initial endpoint.
\end{definition}

\begin{definition}
For an admissible curve define
\[
\begin{aligned}
 \mathcal L^T(\gamma)&=\int_0^{\bar\tau}\sqrt\tau
       \bigl(R_h(\gamma(\tau),T-\tau)+|X|_{h_\tau}^2\bigr)d\tau\\
 &=\int_0^{\sqrt{\bar\tau}}
       \left(\tfrac12|A|_{h(T-s^2)}^2+2s^2R_h(\alpha(s),T-s^2)\right)ds.
\end{aligned}
\]
Its minimizing-action function (the quotient $\mathcal L$-distance) is
\[
 L_{(\bar x,T)}(\bar y,\bar\tau)
 =\inf_{\gamma(0)=\bar x,\,\gamma(\bar\tau)=\bar y}\mathcal L^T(\gamma).
\]
On the total space the notation
${}^T\mathcal L_{(x,T)}^{\bar\tau}(y)$ means, by definition,
$L_{(\pi(x),T)}(\pi(y),\bar\tau)$.
\end{definition}

This pullback is basic in both endpoints. Equivalently, one may take
horizontal lifts starting at $x$, with their final point free anywhere
on the Reeb fiber $\pi^{-1}(\bar y)$. Projection preserves the integrand,
and every quotient curve has a horizontal lift with that initial point;
the principal circle or proper-line bundle in setting (S) has a smooth
invariant connection, so local trivializations along its compact image
continue the lift to the endpoint. Hence the two infima agree.
Fixing the final point in $M$ would define a different, constrained
problem and is not part of the definition. Variations below are variations
on $Z$. Although each curve can be lifted horizontally, the variation
field of those lifts need not be horizontal; only its projection
represents the quotient variation field.

\begin{proposition}
\label{prop1}Critical quotient curves for $\mathcal L^T$ with fixed
quotient endpoints satisfy
\begin{equation}
 E(X):=D_\tau X-\tfrac12\nabla R^T
       +\frac{X}{2\tau}+2(\mathrm{Ric}^T)^\#X=0.
 \label{el}
\end{equation}
Here $D_\tau$ is covariant differentiation along $\gamma$ using $h_\tau$.
\end{proposition}

\begin{proof}
For a quotient variation with field $Y$, differentiation of the action
and integration by parts give
\[
 \delta_Y\mathcal L^T
 =2\sqrt\tau\langle X,Y\rangle\big|_0^{\bar\tau}
       -2\int_0^{\bar\tau}\sqrt\tau\langle E(X),Y\rangle d\tau.
\]
The coefficient $2\mathrm{Ric}^T$ in \eqref{el} comes from
$\partial_\tau h_\tau=2\mathrm{Ric}^T$. The formula is obtained first
on $[\varepsilon,\bar\tau]$ and then by $s=\sqrt\tau$ at zero.
For fixed initial point $Y(0)=0$, its initial boundary term vanishes.
Arbitrary compactly supported quotient variations prove \eqref{el},
and the same identity proves its sufficiency for stationarity.
\end{proof}

\begin{definition}
A transverse $\mathcal L$-geodesic is a quotient solution of \eqref{el}
which is smooth in $s=\sqrt\tau$ up to $s=0$, or a horizontal lift of
that solution with the leafwise endpoint convention just specified.
\end{definition}

The geodesic equation is regular at the initial time in the $s$ variable:
\[
 D_s A=2s^2\nabla R^T-4s(\mathrm{Ric}^T)^\#A.
\]
Thus the initial data are $\alpha(0)=\bar x$ and $A(0)=v\in T_{\bar x}Z$.
The linearized equation must also use this regular initial-value problem.
Let $\mathcal C_\tau=\partial_\tau\nabla^{h_\tau}$; it is the tensor
\[
 \langle\mathcal C_\tau(U,V),W\rangle
 =(\nabla_U\mathrm{Ric}^T)(V,W)
 +(\nabla_V\mathrm{Ric}^T)(U,W)
 -(\nabla_W\mathrm{Ric}^T)(U,V).
\]

\begin{lemma}
The field $Y$ of a variation through quotient $\mathcal L$-geodesics
satisfies
\begin{equation}
 \begin{split}
 0={}&D_\tau^2Y+R^{\nabla}(Y,X)X-\mathcal C_\tau(Y,X)
       -\tfrac12\nabla_Y\nabla R^T+\frac{D_\tau Y}{2\tau}\\
    &+2\bigl(\nabla_Y(\mathrm{Ric}^T)^\#\bigr)X
       +2(\mathrm{Ric}^T)^\#D_\tau Y .
 \end{split}
 \label{jac}
\end{equation}
Here $R^\nabla$ is the curvature endomorphism in \eqref{eqn-SR1}.
\end{lemma}

\begin{proof}
Torsion-freeness gives $\nabla_YX=D_\tau Y$ on the quotient variation.
For its time-dependent connection the commutator is
\[
 \nabla_YD_\tau X=D_\tau\nabla_YX+R^\nabla(Y,X)X-
 \mathcal C_\tau(Y,X).
\]
Differentiating \eqref{el} gives \eqref{jac}.
The displayed formula for $\mathcal C_\tau$ follows by differentiating
the Koszul formula. In particular the evolving-connection term is retained.
\end{proof}

\begin{definition}
A transverse $\mathcal L$-Jacobi field is a solution of \eqref{jac}
whose expression along $\alpha(s)$ is smooth at $s=0$.
Its initial data are $Y(0)$ and $D_sY(0)$.
\end{definition}

\begin{definition}
Let $\mathcal D_x\subset T_{\bar x}Z\times(0,\infty)$ be the open set of
pairs $(v,\tau)$ for which the solution with $\alpha(0)=\bar x$,
$A(0)=v$ exists through $\tau$. Define
\[
 \mathcal L^T\exp_x(v,\tau)=(\gamma_v(\tau),T-\tau)\in Z\times(-\infty,T),
 \qquad
 \widetilde{\mathcal L^T}(v,\tau)=\mathcal L^T(\gamma_v|_{[0,\tau]}).
\]
Write $\mathcal D_x^\tau$, $\mathcal L^T\exp_x^\tau$, and
$\widetilde{\mathcal L_\tau^T}$ for their time slices.
Identifying $T_{\bar x}Z$ with the contact space $\ker\eta_x\subset T_xM$ by
$d\pi$ does not
change the target: the exponential map takes values in $Z$, not in
an open subset of the $(2n+1)$-dimensional total space.
\end{definition}

Smooth dependence for the regular $s$-equation makes these maps smooth.

\begin{proposition}
\label{prop2}Every transverse $\mathcal L$-Jacobi field is obtained by
a smooth variation through quotient $\mathcal L$-geodesics.
\end{proposition}

\begin{proof}
Choose a curve of initial points with derivative $Y(0)$ and a smooth
family of initial vectors whose covariant derivative is $D_sY(0)$.
Solve the regular $s$-equation with these initial data. Its variation
field has the prescribed initial data and satisfies the same linearized
regular equation as $Y$. Uniqueness for that linear ordinary differential
equation proves the assertion. No condition on $D_\tau Y(0)$ is imposed,
since that derivative can be singular.
\end{proof}

\begin{corollary}
\label{cor1}For $v\in\mathcal D_x^\tau$ and $w\in T_{\bar x}Z$ there is a
unique Jacobi field with $Y(0)=0$ and $D_sY(0)=w$, and
\[
 (d\mathcal L^T\exp_x^\tau)_v w=Y(\tau).
\]
\end{corollary}
\begin{proof}
Differentiate the geodesics with initial vectors $v+uw$ at $u=0$.
\end{proof}

\begin{lemma}
\label{lem1}For $v\in\mathcal D_x^\tau$ and $w\in T_{\bar x}Z$,
\[
 d\widetilde{\mathcal L_\tau^T}|_v(w)
 =2\sqrt\tau\left\langle X(\tau),
          (d\mathcal L^T\exp_x^\tau)_vw\right\rangle_{h_\tau}.
\]
\end{lemma}
\begin{proof}
Apply the first variation formula to $\gamma_{v+uw}$. Its integral
term vanishes and its initial point is fixed, leaving the final term.
\end{proof}

\begin{definition}
The injectivity set ${}^T\widetilde{\mathcal U}_x$ consists of pairs
$(v,\tau)$ for which $\gamma_v$ uniquely minimizes $\mathcal L^T$ to
its endpoint and $d\mathcal L^T\exp_x^\tau$ is nonsingular.
Its image ${}^T\mathcal U_x\subset Z\times(-\infty,T)$ is the
quotient injectivity domain; write ${}^T\mathcal U_x(\tau)$ for its
spatial slice. The corresponding basic domain in $M$ is the preimage
under $\pi$, not the image of a $2n$-dimensional exponential map into $M$.
\end{definition}

\begin{proposition}
The map $\mathcal L^T\exp_x:{}^T\widetilde{\mathcal U}_x\to
{}^T\mathcal U_x$ is a diffeomorphism between open sets, and
\[
 L_{(\bar x,T)}(\mathcal L^T\exp_x^\tau(v),\tau)
       =\widetilde{\mathcal L_\tau^T}(v).
\]
\end{proposition}
\begin{proof}
The inverse function theorem applies at every nonconjugate endpoint.
Uniqueness of the minimizing geodesic makes these local inverses agree.
Openness and persistence of minimizing uniqueness follow from the
usual quotient cut-locus argument: a converging sequence of competing
minimizers has a minimizing limit, and uniqueness and nonsingularity
then place it in the same inverse chart. The compactness needed here
is proved in the next lemma. For completeness, its energy bound also
controls the initial velocities of such minimizing geodesics. On their
common compact spatial set and a fixed regular time slab, the regular
$s$-equation and $\partial_s h(T-s^2)=4s\mathrm{Ric}$ imply
$|d|A|/ds|\le C(1+|A|)$. A pointwise bound at one time follows from the
uniform $L^2$ energy bound, and Gronwall's inequality propagates it over
the entire interval. The geodesic equation then gives smooth convergence,
after a subsequence, including convergence of the initial velocities.
Thus the competing initial data really do enter the same inverse chart.
See \cite[Sections~6.2--6.3 and
Proposition~7.5]{MT} for this argument for a smooth time-dependent
Riemannian metric. The equality is the definition of the minimizing value.
\end{proof}

\begin{corollary}
\label{coro1}On ${}^T\mathcal U_x(\tau)$ one has
\[
 \nabla^T L_{(\bar x,T)}(\cdot,\tau)=2\sqrt\tau X(\tau),
\]
where $\gamma$ is the unique minimizing quotient geodesic.
The identical formula holds for the basic pullback.
\end{corollary}
\begin{proof}
Combine Lemma~\ref{lem1} with the inverse exponential chart.
\end{proof}

\begin{lemma}
\label{lem:quotient-minimizer}
Let $[T-\bar\tau,T]$ be a regular closed time slab of the complete
quotient flow on which the curvature is bounded, with $\bar\tau>0$.
For every $\bar x,\bar y\in Z$ there is a minimizing quotient
$\mathcal L$-geodesic from $(\bar x,T)$ to $(\bar y,T-\bar\tau)$.
In particular, every $T<0$ and $\bar\tau>0$ are allowed for a
$C$-controlled ancient flow.
\end{lemma}

\begin{proof}
On the finite slab $[T-\bar\tau,T]$ curvature is uniformly bounded,
so the metrics are uniformly equivalent to the complete metric $h(T)$,
and $R_h$ is uniformly bounded. Set $b=\sqrt{\bar\tau}$.
A minimizing sequence for the $s$-form of the action has uniformly
bounded $\int_0^b|\alpha'|_{h(T)}^2ds$. Its curves therefore have uniformly
bounded $h(T)$ length and lie in one compact ball about $\bar x$.
Cauchy--Schwarz gives uniform one-half H\"older continuity.
A subsequence converges uniformly and weakly in $W^{1,2}$ in finitely
many coordinate charts. The kinetic integrand is uniformly positive
quadratic in the velocity; uniform convergence of its smooth metric
coefficients and weak lower semicontinuity give the lower-semicontinuity
of its integral. The scalar-curvature term converges uniformly under the
integral. The limiting curve thus minimizes the action with the prescribed
endpoints. The Euler--Lagrange equation of this regular Lagrangian gives
a smooth solution in $s$, including the endpoints by ordinary differential
equation regularity. Returning to $\tau=s^2$ gives \eqref{el}.
This is ordinary $W^{1,2}$ compactness in $s$, not an assertion of
unweighted $W^{1,2}$ compactness in $\tau$.
\end{proof}

\begin{remark}
The quotient $\mathcal L$-distance $L_{(\bar x,T)}$ is locally Lipschitz
on its positive-time domain, smooth on ${}^T\mathcal U_x$, and locally
semiconcave in the spatial
variable away from $\tau=0$. Its cut locus has measure zero.
These are the complete bounded-curvature results of
\cite[Proposition~7.5, Theorem~7.13, and Lemma~7.15]{MT}:
one fixes a final segment of a minimizing curve to construct smooth
upper barriers, and applies the second-variation bound to those segments.
The pulled-back distance is consequently basic and locally Lipschitz on
$M$. No assertion here identifies its cut locus with that of a
fixed-endpoint horizontal sub-Riemannian problem (for which see \cite{ABB}).
\end{remark}

\begin{proposition}
\label{prop:LapL}Let $(Z,h(t))$ be the complete quotient flow in
setting (S), with $|\mathrm{Rm}(h(t))|\le C/|t|$ and normalization
\eqref{eq:marked-flow}. In this proposition and its variational formulas,
$x$ and $y$ denote points of $Z$; the same formulas on $M$ mean their
basic pullbacks under $\pi$. We suppress the bars on these quotient
endpoints and use ${}^T\mathcal L$ for $L_{(x,T)}$. For
$y\in {}^{T}\mathcal{U}_{x}(\bar{\tau})$ let $\widetilde{\gamma }$ be
the unique minimizing transverse $\mathcal{L}$-geodesic from
$(x,T)$ to $(y,T-\bar{\tau})$. Then
\begin{equation*}
\Delta _{B}\left( ^{T}\mathcal{L}_{(x,T)}^{\bar{\tau}}\right)
(y)\leq \frac{2n}{\sqrt{\bar{\tau}}}-2\sqrt{\bar{\tau}}R^{T}(y)
-\frac{1}{\bar{\tau}}K^{T}[\widetilde{\gamma }],
\end{equation*}
where $K^{T}[\widetilde{\gamma }]:=
\int_{0}^{\bar{\tau}}\tau ^{\frac{3}{2}}H(X(\tau ))\,d\tau $ and
\begin{equation*}
H(X):= -\frac{\partial R^{T}}{\partial \tau }-\frac{R^{T}}{\tau }
-2\langle \nabla ^{T}R^{T},X\rangle +2\mathrm{Ric}^{T}(X,X)
\end{equation*}
is Perelman's Harnack quantity, computed with respect to
$h_{\tau }$.
\end{proposition}

\begin{proof}
All computations are on $Z$, with its actual Riemannian metric $h_\tau$.
The distance was defined there before being pulled back, so ordinary
quotient variations have no horizontality constraint. We give
Perelman's computation \cite[Section~7]{Per1} (see also
\cite[Sections~6.4 and~7.2]{MT}) in the normalization
$\partial_\tau h_\tau=2\mathrm{Ric}_{h_\tau}$ and real dimension $2n$. We
indicate the three steps. Below, $\nabla _{X}$ denotes covariant
differentiation along $\gamma $ with respect to the Levi-Civita connection
of $h_{\tau }$, and $\langle \cdot ,\cdot \rangle
=h_{\tau }$.

(i) \emph{Second variation.} Let $\widetilde{\Phi }(\tau ,u)$ be a
variation of the minimizing quotient curve through admissible quotient
curves with fixed initial point. The ordinary quotient fields are
$\widetilde{X}=\widetilde{\Phi }_{\ast
}\partial _{\tau }$ and $\widetilde{Y}=\widetilde{\Phi }_{\ast }\partial _{u}$,
and put $X=\widetilde{X}(\cdot ,0)$, $Y=\widetilde{Y}(\cdot ,0)$, so that
$Y(0)=0$. Since $\partial _{\tau }h_{\tau }=2\mathrm{Ric}^{T}$, for vector
fields $A,B$ along $\gamma $
\begin{equation}
\frac{d}{d\tau }\langle A,B\rangle =2\mathrm{Ric}^{T}(A,B)+\langle \nabla
_{X}A,B\rangle +\langle A,\nabla _{X}B\rangle ,  \label{form1}
\end{equation}
and the variation of the Christoffel symbols is given by
\[
\begin{aligned}
\langle (\partial _{\tau }\Gamma )(A,B),W\rangle
 &= (\nabla_{A}\mathrm{Ric}^{T})(B,W)\\
 &\quad +(\nabla _{B}\mathrm{Ric}^{T})(A,W)
       -(\nabla_{W}\mathrm{Ric}^{T})(A,B).
\end{aligned}
\]
Using \eqref{form1}, $\nabla _{X}Y=\nabla
_{Y}X$, the Euler--Lagrange equation \eqref{el} and this formula when
commuting $\nabla _{Y}$ with $\nabla _{X}$, one obtains, exactly as in
\cite[Section~7]{Per1},
\begin{equation}
\begin{split}
&\frac{d^{2}}{du^{2}}\mathcal{L}^{T}(\widetilde{\gamma }_{u})\Big|_{u=0}
=2\sqrt{\bar{\tau}}\langle \nabla _{Y}\widetilde{Y},X\rangle
(\bar{\tau})\\
&\quad +\int_{0}^{\bar{\tau}}\sqrt{\tau
}\Big[\operatorname{Hess}^{T}R^{T}(Y,Y)-2\,\mathrm{Rm}^{T}(X,Y,X,Y) \\
&\qquad\qquad +2|\nabla _{X}Y|^{2}+2(\nabla _{X}\mathrm{Ric}^{T})(Y,Y)
 -4(\nabla_{Y}\mathrm{Ric}^{T})(X,Y)\Big]d\tau ,
\end{split}
\label{form2}
\end{equation}
where $\mathrm{Rm}^{T}(X,Y,X,Y)$ is the sectional curvature form in the
convention \eqref{eqn-SR1}.

(ii) \emph{Hessian estimate along adapted fields.} A field
$V(\tau )\in T_{\gamma(\tau)}Z$ along $\gamma$ is called \emph{adapted} if
\begin{equation}
\nabla _{X}V=-\mathrm{Ric}^{T}(V,\cdot )^{\#}+\frac{1}{2\tau }V.
\label{adap1}
\end{equation}
We choose $Y$ adapted with the required terminal value. In the variable
$s=\sqrt\tau$ it vanishes smoothly at $s=0$ and is realized by an
ordinary quotient curve variation, for example using the pointwise
Riemannian exponential. It need not be a Jacobi field; the varied curves
are competitors, not necessarily geodesics. Let
$q(u)=\widetilde{\Phi }(\bar{\tau},u)$. Since
$^{T}\mathcal{L}_{(x,T)}^{\bar{\tau}}(q(u))\leq
\mathcal{L}^{T}(\widetilde{\gamma }_{u})$ with equality at $u=0$,
\begin{equation}
\frac{d^{2}}{du^{2}}\,{}^{T}\mathcal{L}_{(x,T)}^{\bar{\tau}}(q(u))\Big|_{u=0}\leq
\frac{d^{2}}{du^{2}}\mathcal{L}^{T}(\widetilde{\gamma
}_{u})\Big|_{u=0},  \label{ineq1}
\end{equation}
and, by Corollary~\ref{coro1}, the boundary term in \eqref{form2} is exactly
the first-order term $\langle \nabla _{Y}\widetilde{Y},\nabla
\,{}^{T}\mathcal{L}\rangle $ in the Hessian. By \eqref{form1} and
\eqref{adap1},
$\frac{d}{d\tau }|Y|^{2}=|Y|^{2}/\tau $, hence
\begin{equation}
|Y(\tau )|^{2}=\frac{\tau }{\bar{\tau}}|Y(\bar{\tau})|^{2}.
\label{form3}
\end{equation}
Moreover $2|\nabla _{X}Y|^{2}=2|\mathrm{Ric}^{T}(Y,\cdot )|^{2}-\frac{2}{\tau
}\mathrm{Ric}^{T}(Y,Y)+\frac{|Y|^{2}}{2\tau ^{2}}$ and
\begin{equation*}
\frac{d}{d\tau }\mathrm{Ric}^{T}(Y,Y)=\frac{\partial
\mathrm{Ric}^{T}}{\partial \tau }(Y,Y)+(\nabla
_{X}\mathrm{Ric}^{T})(Y,Y)-2|\mathrm{Ric}^{T}(Y,\cdot )|^{2}+\frac{1}{\tau
}\mathrm{Ric}^{T}(Y,Y).
\end{equation*}
Inserting these identities into \eqref{form2}, one finds that the integrand
equals
\begin{equation*}
-\sqrt{\tau }H(X,Y)-2\frac{d}{d\tau }\big(\sqrt{\tau
}\,\mathrm{Ric}^{T}(Y,Y)\big)+\frac{|Y|^{2}}{2\tau ^{3/2}},
\end{equation*}
where
\begin{equation*}
\begin{split}
H(X,Y)& := -\operatorname{Hess}^{T}R^{T}(Y,Y)+2\,\mathrm{Rm}^{T}(X,Y,X,Y)\\
&\quad -4(\nabla _{X}\mathrm{Ric}^{T})(Y,Y)
       +4(\nabla _{Y}\mathrm{Ric}^{T})(X,Y) \\
&\quad -2\frac{\partial \mathrm{Ric}^{T}}{\partial \tau
}(Y,Y)+2|\mathrm{Ric}^{T}(Y,\cdot )|^{2}-\frac{1}{\tau }\mathrm{Ric}^{T}(Y,Y).
\end{split}
\end{equation*}
Integrating, and using \eqref{ineq1}, \eqref{form3} and $Y(0)=0$, we
obtain, for $|Y(\bar{\tau})|=1$,
\begin{equation}
\operatorname{Hess}^{T}\left( ^{T}\mathcal{L}_{(x,T)}^{\bar{\tau}}\right)
(Y(\bar{\tau}),Y(\bar{\tau}))\leq
\frac{1}{\sqrt{\bar{\tau}}}-2\sqrt{\bar{\tau}}\,\mathrm{Ric}^{T}(Y,Y)(\bar{\tau})-\int_{0}^{\bar{\tau}}\sqrt{\tau
}H(X,Y)\,d\tau .
\label{ineq2}
\end{equation}

(iii) \emph{Trace.} Let $\{Y_{j}(\bar{\tau})\}_{j=1}^{2n}$ be an
orthonormal basis of $T_yZ$ at time $T-\bar\tau$, extended to adapted
fields along $\gamma $. By \eqref{form1} and \eqref{adap1},
$\langle Y_{j},Y_{k}\rangle (\tau )=\frac{\tau }{\bar{\tau}}\delta
_{jk}$, so $W_{j}=\sqrt{\bar{\tau}/\tau }\,Y_{j}$ are orthonormal.
Summing \eqref{ineq2} over $j$ and using $\sum_{j}H(X,W_{j})=H(X)$, which
follows from the contracted second Bianchi identity,
$\sum_{j}\frac{\partial \mathrm{Ric}^{T}}{\partial \tau
}(W_{j},W_{j})=\frac{\partial R^{T}}{\partial \tau }+2|\mathrm{Ric}^{T}|^{2}$
and
$\frac{\partial R^{T}}{\partial \tau }=-\Delta
_{B}R^{T}-2|\mathrm{Ric}^{T}|^{2}$, we obtain
\begin{equation*}
\Delta _{B}\left( ^{T}\mathcal{L}_{(x,T)}^{\bar{\tau}}\right) (y)\leq
\frac{2n}{\sqrt{\bar{\tau}}}-2\sqrt{\bar{\tau}}R^{T}(y)-\int_{0}^{\bar{\tau}}\frac{\tau
^{3/2}}{\bar{\tau}}H(X)\,d\tau
=\frac{2n}{\sqrt{\bar{\tau}}}-2\sqrt{\bar{\tau}}R^{T}(y)-\frac{1}{\bar{\tau}}K^{T}[\widetilde{\gamma
}].
\end{equation*}
With the normalization \eqref{eq:marked-flow} this is consistent with
\eqref{eq:e3} below.
\end{proof}

\subsection{The Reduced Transverse Length}

We make the quotient measure explicit and compare it with ambient tube
volume. The next comparison allows arbitrary Reeb type solely to explain
why the global reduced-volume theory is formulated on the smooth quotient
in setting (S); it does not extend that theory to irregular orbit spaces.

By Lemma~\ref{lem:reeb}, applied at a fixed time, either the closure $H$ of
the Reeb flow in the isometry group is a compact torus, or the Reeb flow is
a closed, free and proper $\mathbb{R}$-action whose quotient $Z=M/\mathbb{R}$
is a smooth complete K\"{a}hler manifold.

The same group acts isometrically at every time: the Reeb flow preserves
each transverse metric and $\eta $, and its limits in the initial isometry
group preserve the later smooth tensors as well. In particular, the compact
torus $H$ is fixed along the flow. Every compact Sasakian manifold belongs
to the compact-torus case. This includes irregular structures.

Throughout, $m=2n$ denotes the real dimension of the transverse space. In the
compact-torus case, write
\begin{equation*}
P_{x}=Hx=\overline{\operatorname{orb}_{\xi }x},\qquad \mathcal{Q}=M/H.
\end{equation*}
We do not identify $\mathcal{Q}$ with an $m$-dimensional orbifold. The
orbits $P_{x}$ may have different dimensions. The quotient is used as a
proper metric space, not as a substitute for a transverse coordinate chart.

For comparison with the main text, the full tube measure is
\begin{equation}
d\mu_t^{\mathrm{tube}}=\eta\wedge\frac{(\omega_t^T)^n}{n!}.
\label{eq:measure}
\end{equation}
It is the volume of $h(t)+\eta^2$, which need not be Sasakian at every time.
Replacing $\eta$ by $\eta+\beta$ for basic $\beta$ leaves this measure
unchanged. In setting (S), however, we henceforth use
\begin{equation}
(\mathcal{X},d\mu_t^T)=(Z,d\mathrm{vol}_{h_Z(t)}).  \label{eq:analyticspace}
\end{equation}
Thus integrals and test functions in global statements are on $Z$; their
pullbacks are basic functions on the chosen total space. In the proper-line
case no integration over the infinite Reeb fiber is performed. The quotient
identities are
\begin{equation}
\partial_t d\mu_t^T=-R^T d\mu_t^T,\qquad \int_Z v\Delta_Bu\,d\mu_t^T
=-\int_Z\langle\nabla^Tv,\nabla^Tu\rangle d\mu_t^T.  \label{eq:measureevolve}
\end{equation}
The integration-by-parts identity is stated for smooth compactly
supported quotient test functions; noncompactly supported versions below
are justified by the indicated integrability and cutoff arguments.
All scaling in these global statements is ordinary quotient parabolic
scaling. It does not assert convergence of contact bundles or Reeb fields.

For contrast, if $P_x$ is a compact $k$-dimensional Reeb orbit closure, a
Gaussian tube diagnostic at a fixed smooth time has the asymptotic
\begin{equation}
\begin{aligned} \mathcal V^{\mathrm{tube}}_{x,T}(\tau) &:=
(4\pi\tau)^{-m/2}\int_M
e^{-d_{g(T)}(P_x,y)^2/(4\tau)}d\mathrm{vol}_{g(T)}(y)\\
&=\operatorname{Vol}_{g(T),k}(P_x)(4\pi\tau)^{(1-k)/2}(1+o(1)). \end{aligned}
\label{12}
\end{equation}
Here one assumes a Gaussian-integrable volume tail, as holds under the
bounded ambient curvature hypotheses used in the main text. This is a
tubular-neighborhood calculation, not a definition of the normalized reduced
volume in setting (S). It illustrates why the full tube measure
\eqref{eq:measure} is not used
for the normalized reduced volume when the orbit closures have dimension
$k>1$; no global irregular reduced-volume theorem is asserted.

We now introduce the reduced transverse length, which is invariant under
parabolic scaling.

\begin{definition}
The reduced transverse length based at $(x,T)$ is
\begin{equation*}
{\ell
^{T}}_{(x,T)}(y,\bar{\tau}):=\frac{{}^{T}\mathcal{L}_{(x,T)}^{\bar{\tau}}(y)}{2\sqrt{\bar{\tau}}},
\end{equation*}
and $K^{T}$ is the transverse Harnack integral of Proposition~\ref{prop:LapL}.
We write $\ell _{x,\bar{\tau}}^{T}:=\ell _{(x,T)}^{T}(\cdot ,\bar{\tau})$
when $T$ is fixed.
\end{definition}

\begin{lemma}
\label{lem:ell} Let $q$ lie in the transverse injectivity domain
${}^{T}\mathcal{U}_{x}(\bar{\tau})$ and let $\tilde{\gamma}_{q}$ be the unique
minimizing transverse $\mathcal{L}$-geodesic from $(x,T)$ to
$(q,T-\bar{\tau})$. Then at $(q,T-\bar{\tau})$
\begin{align}
\frac{\partial {\ell ^{T}}}{\partial \bar{\tau}}& ={R^{T}}-\frac{{\ell
^{T}}}{\bar{\tau}}+\frac{1}{2\bar{\tau}^{3/2}}{K^{T}}[\tilde{\gamma}_{q}],
\label{eq:e1} \\
|{\nabla ^{T}}{\ell ^{T}}|^{2}& =\frac{{\ell
^{T}}}{\bar{\tau}}-{R^{T}}-\frac{1}{\bar{\tau}^{3/2}}{K^{T}}[\tilde{\gamma}_{q}],
\label{eq:e2} \\
\Delta _{B}{\ell ^{T}}& \leq
\frac{n}{\bar{\tau}}-{R^{T}}-\frac{1}{2\bar{\tau}^{3/2}}{K^{T}}[\tilde{\gamma}_{q}].
\label{eq:e3}
\end{align}
Consequently, with $\delta \geq 0$ the defect in \eqref{eq:e3},
\begin{align}
\frac{\partial {\ell ^{T}}}{\partial \bar{\tau}}-\Delta _{B}{\ell
^{T}}+|{\nabla ^{T}}{\ell ^{T}}|^{2}-{R^{T}}+\frac{n}{\bar{\tau}}& =\delta \
\geq \
0,  \label{eq:magic} \\
2\Delta _{B}{\ell ^{T}}-|{\nabla ^{T}}{\ell ^{T}}|^{2}+{R^{T}}+\frac{{\ell
^{T}}-2n}{\bar{\tau}}& =-2\delta \ \leq \ 0.  \label{eq:magic2}
\end{align}
The displayed identities and the function $\delta$ are pointwise on the
injectivity domain. Across the cut locus, \eqref{eq:magic} and
\eqref{eq:magic2} hold as distributional inequalities; the nonnegative
defect there is understood distributionally and need not be a smooth
function.
\end{lemma}

\begin{proof}
We argue as in \cite[Chapter~7]{MT} and write
$\mathcal{L}^{T}={}^{T}\mathcal{L}_{(x,T)}^{\bar{\tau}}$ and
$K^{T}=K^{T}[\tilde{\gamma}_{q}]$. On the injectivity domain $\nabla
^{T}\mathcal{L}^{T}=2\sqrt{\bar{\tau}}X(\bar{\tau})$ (Corollary~\ref{coro1})
and $\partial
_{\bar{\tau}}\mathcal{L}^{T}=\sqrt{\bar{\tau}}(R^{T}-|X|^{2})$, while
integrating
\begin{equation*}
\frac{d}{d\tau }[\tau ^{3/2}(R^{T}+|X|^{2})]=\frac{1}{2}\sqrt{\tau
}(R^{T}+|X|^{2})-\tau ^{3/2}H(X)
\end{equation*}
along $\tilde{\gamma}_{q}$ gives
\begin{equation*}
\bar{\tau}^{3/2}(R^{T}+|X|^{2})(\bar{\tau})=\frac{1}{2}\mathcal{L}^{T}-K^{T}.
\end{equation*}
This yields \eqref{eq:e1} and \eqref{eq:e2}, and \eqref{eq:e3} is
Proposition~\ref{prop:LapL} divided by $2\sqrt{\bar{\tau}}$. Then
\eqref{eq:magic} is \eqref{eq:e1}$-$\eqref{eq:e3}$+$\eqref{eq:e2} and
\eqref{eq:magic2} is $2\times $\eqref{eq:e3}$-$\eqref{eq:e2}; the
distributional inequalities use more than local Lipschitz regularity:
fixing an initial part of a minimizing quotient geodesic gives smooth
upper barriers for its final segment. The second-variation bound yields
local spatial semiconcavity. The singular part of the distributional
Hessian is nonpositive, so the Laplacian upper bound, and hence
\eqref{eq:magic} and \eqref{eq:magic2}, extend across the cut locus.
This is the quotient argument of \cite[Theorem~7.13]{MT}; the cut locus
has measure zero and the first-order identities hold almost everywhere.
\end{proof}

Throughout the rest of this subsection the complete quotient flow
$(Z,h(t))$ is $C$-controlled and $(\bar x,T)\in Z\times(-\infty,0)$
is fixed. Symbols $x,y$ in transverse notation mean their projections
when a quotient distance or integral is taken. No noncollapsing
hypothesis is imposed in the estimates that follow.

\begin{proposition}
\label{prop:quadratic} For a controlled ancient transverse flow and every
$T<0$, there are $A,\mu $, independent of $T$, such that
\begin{align}
\ell ^{T}(y,\tau )& \geq -\mu ,\qquad |\ell ^{T}(x,\tau )|\leq \mu ,
\label{eq:basebound} \\
A^{-1}\left( 1+\frac{d_{T-\tau }^{T}(x,y)}{\sqrt{\tau }}\right) ^{2}-A& \leq
\ell ^{T}(y,\tau )\leq A\left( 1+\frac{d_{T-\tau }^{T}(x,y)}{\sqrt{\tau
}}\right) ^{2},  \label{eq:twosided} \\
|{\nabla ^{T}}\ell ^{T}|(y,\tau )& \leq \frac{A}{\sqrt{\tau }}\left(
1+\frac{d_{T-\tau }^{T}(x,y)}{\sqrt{\tau }}\right) ,  \label{eq:gradbound} \\
|\partial _{\tau }\ell ^{T}|(y,\tau )& \leq \frac{A}{\tau }\left(
1+\frac{d_{T-\tau }^{T}(x,y)}{\sqrt{\tau }}\right) ^{2}.  \label{eq:timebound}
\end{align}
Derivative assertions are understood almost everywhere and in the
corresponding Lipschitz sense.
\end{proposition}

\begin{proof}
These are the complete-quotient estimates of
\cite[Propositions~2.1--2.3 and Lemmas~2.5--2.6]{Nab}, in real dimension
$2n$. The proofs of these estimates use curvature control and
completeness, not the noncollapsing part of Naber's controlled-flow
terminology. We specify the inputs to avoid a circular use below.
At time $T-\sigma$, $|\mathrm{Rm}|\le C/|T-\sigma|\le C/\sigma$;
Shi's estimates give the corresponding spatial and time derivative bounds for
the scalar curvature.
Thus $|R_h|\le C_1/\sigma$, which bounds the scalar part of the action
below by $-2C_1\sqrt\tau$. The constant quotient curve at $\bar x$
is an admissible competitor and gives the corresponding upper bound
there. This proves \eqref{eq:basebound}.
The Harnack integral estimate along minimizing quotient geodesics is
\[
 |K^T|\le A\sqrt\tau(1+\mu+\ell^T).
\]
It follows by inserting the curvature and derivative bounds in $H(X)$,
using Cauchy--Schwarz for the $\nabla R$ term, and controlling
$\int\sqrt\sigma|X|^2$ by the minimizing action and its scalar part.
Equation~\eqref{eq:e2} now gives
$|\nabla\ell^T|^2\le A\tau^{-1}(1+\mu+\ell^T)$.
Integration along $h(T-\tau)$-minimizing geodesics from $\bar x$
proves the upper quadratic bound and \eqref{eq:gradbound}.
For the lower bound one also uses the index-form estimate
\[
 \int_\beta\mathrm{Ric}_{h(T-s)}(\dot\beta,\dot\beta)
       \le A/\sqrt{|T-s|}
\]
for every unit-speed minimizing quotient geodesic $\beta$.
Indeed, cut off parallel normal fields over lengths $\sqrt{|T-s|}$
at its endpoints; the index form controls the middle integral and
$|\mathrm{Ric}|\le C_2/|T-s|$ controls the end intervals.
The resulting distance distortion estimate, compared along a minimizing
$\mathcal L$-geodesic, is exactly the lower-bound argument in
\cite[Proposition~2.3]{Nab} and gives
$\ell^T\ge A^{-1}(1+d_{h(T-\tau)}(\bar x,\bar y)/\sqrt\tau)^2-A$.
Finally \eqref{eq:e1} and the Harnack-integral estimate imply
\eqref{eq:timebound}. Every constant in these arguments depends only
on $n,C$; the bounds extend almost everywhere and in the local
Lipschitz sense across the quotient cut locus.
\end{proof}

\begin{lemma}
\label{lem:BG}If $\mathrm{Ric}_{h_{Z}}\geq -(m-1)b^{2}h_{Z}$, then
\begin{equation}
\frac{\mu ^{T}(B_{Z}(x,R))}{\mu ^{T}(B_{Z}(x,r))}\leq \frac{\mu
_{m,-b^{2}}(R)}{\mu _{m,-b^{2}}(r)},\qquad 0<r<R.  \label{eq:BG}
\end{equation}
Here $\mu _{m,-b^{2}}(r)$ denotes the volume of a geodesic ball of radius
$r$ in the $m$-dimensional simply connected space form of constant
sectional curvature $-b^{2}$. In particular the quotient volume has at most exponential growth.
\end{lemma}

\begin{proof}
This is ordinary Bishop--Gromov comparison on the complete smooth
$m$-dimensional quotient $Z$. For a regular lift, multiplying both volumes by
the same fixed fiber period gives the corresponding tube ratio. No
comparison theorem for an arbitrary irregular orbit space is needed or
asserted in this appendix.
\end{proof}

\subsection{The Reduced Transverse Volume and the Singular Reduced
Transverse Length}

Setting (S), with the quotient measure \eqref{eq:analyticspace}, is assumed
throughout this subsection. Define
\begin{equation}
u_{L}=(4\pi \tau )^{-m/2}e^{-\ell ^{T}},\qquad V_{(x,T)}^{T}(\tau
)=\int_{Z}u_{L}\,d\mu _{T-\tau }^{T}.  \label{eq:reducedvolume}
\end{equation}
Under the parabolic scaling $h_{c}(t)=c\,h(c^{-1}t)$ one has
$\mathcal{L}_{c}^{T}=\sqrt{c}\,\mathcal{L}^{T}$ at corresponding times, so
both $\ell ^{T}$ and \eqref{eq:reducedvolume} are invariant under
transverse parabolic scaling.

\begin{definition}
\label{def:redvol} For a measurable subset $\mathcal A\subset
Z\times(-\infty,0)$ define its spatial time slice by
$\mathcal A_{\bar\tau}:=\{y\in Z:(y,T-\bar\tau)\in\mathcal A\}$ and put
\begin{equation*}
{V^{T}}_{\mathcal{A}}(\bar{\tau}):=\int_{\mathcal{A}_{\bar{\tau}}}(4\pi
\bar{\tau})^{-n}e^{-{\ell ^{T}}_{x,\bar{\tau}}}\,d\mu^{T}_{T-\bar\tau}.
\end{equation*}
When $\mathcal{A}_{\bar{\tau}}=Z$ for all $\bar{\tau}$ we write
$V^{T}(\bar{\tau})$ and call it the reduced transverse volume.
\end{definition}

\begin{proposition}
\label{prop:volmonotone}Let the complete quotient flow in setting (S)
be $C$-controlled.

\begin{enumerate}
\item $(\partial _{\bar{\tau}}\ell ^{T}+\frac{n}{\bar{\tau}}-R^{T})e^{-\ell
^{T}}\in L^{1}(Z,d\mu ^{T})$, and \eqref{eq:magic} holds in the
distributional sense against test functions $\psi e^{-\ell ^{T}}$ with $\psi
$ and $|\nabla ^{T}\psi |$ bounded. In particular $\int_{Z}(\Delta _{B}\ell
^{T}-|\nabla ^{T}\ell ^{T}|^{2})e^{-\ell ^{T}}d\mu ^{T}=0$ in this extended
sense.

\item $V_{(x,T)}^{T}(\tau )$ is finite and nonincreasing in $\tau $, and
$V_{(x,T)}^{T}(\tau )\rightarrow 1$ as $\tau \downarrow 0$.
\end{enumerate}
\end{proposition}

\begin{proof}
Work on $Z$ with $m=2n$ and put
$u=(4\pi\tau)^{-m/2}e^{-\ell^T}$. Proposition~\ref{prop:quadratic}
and Lemma~\ref{lem:BG} imply Gaussian decay of $u$, and integrability
of $u|\nabla\ell^T|$, $u|\nabla\ell^T|^2$, and
$u|\partial_\tau\ell^T|$ at every positive time, locally uniformly
on compact time intervals. The distributional inequality
\eqref{eq:magic}, by the chain rule for locally semiconcave functions,
is
\[
 (\partial_\tau-\Delta_h+R_h)u\le0.
\]
Use compactly supported spatial cutoffs $\phi_R$ with
$0\le\phi_R\le1$, $\phi_R\to1$, and $|\nabla\phi_R|\le C/R$.
The time-slice integration by parts is interpreted distributionally;
for example,
\[
 \int_Z\phi_R(\Delta_h\ell^T-|\nabla\ell^T|^2)e^{-\ell^T}
       d\mathrm{vol}_{h}
 =-\int_Z\langle\nabla\phi_R,\nabla\ell^T\rangle e^{-\ell^T}
       d\mathrm{vol}_{h}\longrightarrow0.
\]
The same bound applies with an extra factor $\psi$ of bounded value
and gradient. Smooth approximation on compact sets justifies these
tests against the weak inequalities, and Gaussian integrability makes
the cutoff errors vanish. This proves (1), including its stated extended
pairing. Since
$\partial_\tau d\mathrm{vol}_{h(T-\tau)}=R_hd\mathrm{vol}_{h(T-\tau)}$,
integration of the inequality for $u$ on any compact positive-time
interval gives $V(\tau_2)\le V(\tau_1)$ for $\tau_2>\tau_1$.
At a regular basepoint the parabolically rescaled quotient converges
locally to Euclidean space and $\ell^T$ converges to
$|z|^2/4$. The quadratic lower bound and the volume comparison control
the Gaussian tails in this rescaling. Hence dominated convergence
gives $\lim_{\tau\downarrow0}V(\tau)=
(4\pi)^{-m/2}\int_{\mathbb R^m}e^{-|z|^2/4}dz=1$.
This uses the quotient volume, not the full tube measure
\eqref{eq:measure}.
\end{proof}

\begin{proposition}
\label{prop:singular} Let $T_{i}\uparrow 0$ and fix $x$. After a
subsequence, the functions
\begin{equation*}
l_{i}(y,t)=\ell _{(x,T_{i})}^{T}(y,T_{i}-t),\qquad t<T_{i},
\end{equation*}
converge locally uniformly on $Z\times(-\infty,0)$ to a locally
Lipschitz function $l(y,t)$, and their pullbacks converge locally
uniformly to its basic pullback on $M$. Their spatial gradients converge
strongly in $L_{\mathrm{loc}}^{p}$ for every finite $p$. With $\tau =-t$, the
limit
satisfies the Hamilton--Jacobi identity
$2\partial_{\tau}l+|\nabla^{T}l|^{2}=R^{T}-l/\tau$ (twice \eqref{eq:e1} plus
\eqref{eq:e2}), \eqref{eq:magic}, \eqref{eq:magic2}, and the bounds of
Proposition~\ref{prop:quadratic}. Its reduced volume is finite and
nonincreasing at positive
scales.
\end{proposition}

\begin{proof}
Proposition~\ref{prop:quadratic} is uniform in $T_i$ and gives local
space-time equicontinuity on $Z\times(-\infty,0)$. A diagonal
Arzel\`a--Ascoli argument gives a locally uniform limit.
For completeness, the strong gradient passage can be made on a
compact set $K\times[a,b]$, $b<0$. For large $i$ every minimizing
curve to this set has a terminal time segment of a fixed positive
length wholly inside a compact subinterval of $(-\infty,0)$.
On that segment the curvature and all its needed derivatives are
uniformly bounded, and the terminal velocity is bounded by
\eqref{eq:gradbound} and the first-variation formula. Its regular
geodesic equation bounds the velocity on a shorter fixed segment.
Hold the preceding part of the minimizing curve fixed and vary only
this final segment. Its second variation supplies upper Hessian
barriers with one constant on $K\times[a,b]$, independent of $i$;
therefore the functions are uniformly semiconcave in space there.
Local uniform convergence of uniformly semiconcave functions gives
almost-everywhere convergence of spatial gradients at differentiability
points. The uniform Lipschitz bounds then imply strong local
$L^p$ convergence for every finite $p$.
The time derivatives converge weakly. Passing to the almost-everywhere
Hamilton--Jacobi identity with these strong spatial gradients proves
the claimed identity distributionally and hence almost everywhere.
The weak forms of \eqref{eq:magic} and \eqref{eq:magic2} pass to the
limit in the same way; all other bounds follow from uniform convergence.
This is the singular reduced-length construction of
\cite[Proposition~2.5]{Nab}, now for the specified quotient object.
For fixed negative time the measure is $d\mathrm{vol}_{h(t)}$ for all
$i$. The uniform quadratic lower bound and Bishop--Gromov comparison
give an integrable Gaussian majorant, so the reduced volumes converge.
Their finiteness and monotonicity pass to the limit from
Proposition~\ref{prop:volmonotone}.
\end{proof}

\begin{definition}
For a subsequential limit $l(y,t)$ in Proposition~\ref{prop:singular},
set $\ell^T(y,\tau):=l(y,-\tau)$ for $\tau>0$. We call this function
\emph{a singular reduced transverse length} from $(x,0)$, and define
its \emph{singular reduced transverse volume} by
\begin{equation*}
\widetilde V^T(\tau)=\int_Z(4\pi\tau)^{-n}e^{-\ell^T(y,\tau)}
\,d\mu^T_{-\tau}(y).
\end{equation*}
A choice of subsequential limit is fixed whenever these quantities are
used; uniqueness is not asserted.
\end{definition}

\begin{remark}
By Proposition~\ref{prop:volmonotone}, each regular reduced volume is at
most one, so $0<\widetilde{V}^{T}\leq 1$, and $\widetilde{V}^{T}$ is
nonincreasing. We write
\begin{equation}
c_{0}=\lim_{\tau \downarrow 0}\widetilde{V}^{T}(\tau )\in (0,1],\qquad
c_{\infty }=\lim_{\tau \rightarrow \infty }\widetilde{V}^{T}(\tau )\in
\lbrack 0,1].  \label{eq:endpointvolumes}
\end{equation}
If the quotient flow is $\kappa $-noncollapsed, then \eqref{eq:basebound}
and \eqref{eq:twosided}, integrated on a curvature-admissible quotient ball
of radius $a(m,C)\sqrt{\tau }$, give $c_{\infty }\geq c(m,C)\kappa >0$.
\end{remark}

\vskip12pt

\noindent{\bf Ackowledgements} Shu-Cheng Chang is partially supported by the Startup Foundation
for Advanced Talents of the Shanghai Institute for Mathematics and
Interdisciplinary Sciences (No.~2302-SRFP-2024-0049). Chien
Lin is partially supported by the NSTC (No.~114-2115-M-003-002-MY2).
Hongbing Qiu is partially supported by the NSFC
(No.~12471050).

\vskip12pt

\noindent{\bf Declaration of AI use.} The authors used ChatGPT (OpenAI) 
during the
preparation of the manuscript for language editing and exploratory calculations.
All mathematical arguments and computations were independently checked by the
authors.

\end{document}